\documentclass[11pt,a4paper]{article}
\usepackage[margin=27mm,headheight=15pt]{geometry}

\usepackage[titles]{tocloft}

\usepackage{titlesec}

\titleformat*{\section}{\large\bfseries}
\titleformat*{\subsection}{\normalsize\bfseries}
\titleformat*{\subsubsection}{\normalsize\bfseries}

\usepackage[T1]{fontenc}
\usepackage{lmodern,microtype,amsmath,amssymb,amsthm,mathtools,mathrsfs}
\usepackage{booktabs,tabularx,array,graphicx,enumitem,fancyhdr}
\usepackage{tikz}\usetikzlibrary{arrows.meta,positioning,calc,fit}
\usepackage{fvextra}
\usepackage[dvipsnames]{xcolor}
\usepackage[colorlinks=true,linkcolor=MidnightBlue,citecolor=MidnightBlue,urlcolor=MidnightBlue]{hyperref}
\usepackage[nameinlink,noabbrev]{cleveref}
\hypersetup{pdftitle={Zero area of the quartic Fibonacci Julia set},pdfauthor={},pdfsubject={Analytic criterion and computer-assisted quartic application}}
\numberwithin{equation}{section}
\newtheorem{theorem}{Theorem}[section]
\newtheorem{proposition}[theorem]{Proposition}
\newtheorem{lemma}[theorem]{Lemma}
\newtheorem{corollary}[theorem]{Corollary}
\theoremstyle{definition}\newtheorem{definition}[theorem]{Definition}
\theoremstyle{remark}\newtheorem{remark}[theorem]{Remark}
\newcommand{\C}{\mathbb C}\newcommand{\R}{\mathbb R}\newcommand{\D}{\mathbb D}
\newcommand{\area}{\operatorname{area}}\newcommand{\Dom}{\operatorname{Dom}}
\newcommand{\CritVal}{\operatorname{CritVal}}\newcommand{\Pc}{\mathcal P}
\newcommand{\eps}{\varepsilon}\newcommand{\wh}[1]{\widehat{#1}}
\newcommand{\cl}[1]{\overline{#1}}
\newcommand{\cR}{\mathscr R}\newcommand{\cF}{\mathscr F}\newcommand{\cS}{\mathscr S}
\newcommand{\cB}{\mathscr B}\newcommand{\cX}{\mathscr X}\newcommand{\cE}{\mathscr E}
\newcommand{\HH}{\widehat H}
\newcommand{\disk}[1]{\mathbb D_{#1}}\newcommand{\Comp}{\operatorname{Comp}}
\newcommand{\Ak}{\mathcal A}\newcommand{\Jk}{\mathcal J}
\setlist{leftmargin=*,itemsep=3pt,topsep=5pt}
\title{Zero area of the canonical quartic\\Fibonacci Julia set\\[5pt]
\large An analytic criterion for even degrees $\ell\ge4$\\and a computer-assisted quartic application}
\author{}\date{12 November 2026}
\definecolor{diagramblue}{RGB}{36,86,123}
\definecolor{diagramteal}{RGB}{27,117,109}
\definecolor{diagramorange}{RGB}{173,92,23}
\tikzset{
  manuscript diagram/.style={font=\small, line cap=round, line join=round,
    >=Latex},
  map arrow/.style={->, line width=.65pt, shorten <=2pt, shorten >=2pt},
  map label/.style={font=\scriptsize, inner sep=2pt, align=center},
  cell box/.style={draw=diagramblue, line width=.6pt, rounded corners=2pt,
    fill=diagramblue!3, minimum height=.8cm, inner sep=5pt, align=center}
}

\title{Zero area of quartic
Fibonacci Julia sets}
\author{Runze Zhang}
\date{\vspace{-5ex}}
\begin{document}
\maketitle

\thispagestyle{plain}
\begin{abstract}
We give a computer-assisted proof that the Julia set of the degree-four Fibonacci renormalization 2-periodic point has zero planar Lebesgue measure. As a corollary, the degree four Fibonacci polynomial $z^4+c_{\mathrm{Fib}}$ has zero measure Julia set.
\end{abstract}
\vspace{3pt}

\setcounter{tocdepth}{1}
\begingroup\small\setlength{\parskip}{0pt}
\tableofcontents
\endgroup

% BEGIN sections/01_introduction.tex
\section{Introduction}\label{sec:main}

Fibonacci maps are real unimodal maps with a recurrent critical point
whose closest-return combinatorics is governed by the Fibonacci
sequence. For their holomorphic extensions, determining whether the
Julia set has zero or positive planar Lebesgue measure is a
long-standing problem. The quadratic Fibonacci Julia set has zero
area, as follows from Lyubich's theorem on quadratic polynomials
which are not infinitely renormalizable and have no irrationally
indifferent periodic orbit \cite{LyubichMeasure,LM}. To our knowledge,
before the present work, degree two was the only critical degree
for which zero area had been established in the Fibonacci family.
The quadratic argument uses the geometry of Yoccoz puzzles.
Its zero-area estimates do not extend directly to higher critical
degree, where the geometry of the critical nest is qualitatively
different \cite{BKNS,LM}.

Avila and Lyubich developed a general framework for studying the
measure and dimension of Julia sets through renormalization
\cite{ALDimension}. For Feigenbaum maps with periodic renormalization
combinatorics, they established a trichotomy between the
\emph{lean}, \emph{balanced}, and \emph{black-hole} regimes.
In the lean regime the Julia set has Hausdorff dimension strictly
less than two; in the balanced regime it has Hausdorff dimension
two but zero area; and in the black-hole regime it has positive
area. Their approach relates these alternatives to the competition
between orbits entering deeper critical neighborhoods and orbits
escaping from the corresponding return dynamics.

Fibonacci maps are not renormalizable in the usual
Douady--Hubbard sense of polynomial-like renormalization
\cite{DH}. They are, however, infinitely renormalizable in
Lyubich's generalized sense, as developed for Fibonacci dynamics
by Lyubich and Milnor \cite{LM}. In this setting a renormalization
has two source domains, carrying a critical branch and a univalent
branch, with a common range. This structure makes the
landing--escape approach available in the Fibonacci setting as
well; Avila and Lyubich explicitly discuss this extension of
their methods in \cite{ALDimension}.

For the classical quadratic period-doubling Feigenbaum map,
Dudko and Sutherland implemented this approach through an explicit
recursive area estimate and a rigorous numerical verification
\cite{DS}. They proved that its Julia set has Hausdorff dimension
strictly less than two, and hence zero area. Their argument combines
a decomposition into suitable inverse copies with quantitative
distortion bounds, reducing the conclusion to a finite numerically verifiable inequality.
Our proof follows this strategy and adapts it to the two-branch
geometry of Fibonacci renormalization.

We consider the \emph{canonical quartic Fibonacci map}: a normalized
representative of the period-two orbit of generalized Fibonacci
renormalization in critical degree four. Its critical branch is
given by the canonical solution of the Cvitanovi\'c--Fibonacci
equation (for the existence and uniqueness of the solution, see the work of H. Epstein \cite{E} and Nowicki-van Strien \cite{NS})
\[
 f(z)=-\alpha^{-2}f\bigl(\alpha f(\alpha z)\bigr),
 \qquad f(0)=1,\qquad 0<\alpha<1,
\]
with $f(z)=F(z^4)$ and $F'(0)<0$. 
The precise normalization, branch domains, and associated
generalized polynomial-like map are given in
Section~\ref{sec:preliminaries}.

\begin{theorem}[The quartic Fibonacci Julia set has zero area]
\label{thm:main}
Let $D_{\rm can}:U^0\sqcup U^1\to V$ be the normalized canonical
Fibonacci map of critical degree four, with the positive lateral
sign in \eqref{eq:physical-D}. Then
\[
 \area K(D_{\rm can})=\area J(D_{\rm can})=0.
\]
\end{theorem}

The analytic part of the proof applies to canonical Fibonacci maps
of every even critical degree $\ell\ge4$. We use the nested critical
neighborhoods
\[
 W_j=(-\alpha^2)^jV
\]
and study the sets of points which eventually enter these
neighborhoods. Let $a_j$ denote the corresponding normalized
landing areas, as defined in \eqref{eq:L-definition}.
The main objective is to prove that $a_j$ tends to zero.
This implies zero area of the filled Julia set: the points whose
orbits accumulate at the critical point must enter arbitrarily
deep critical neighborhoods, while the remaining surviving points
form a null set. The precise argument is given in
Lemma~\ref{lem:null-detection}.

For the numerical part, the criterion we verify is Theorem~\ref{thm:supported}.
It is essentially a comparison between an upper bound
for landing area and a lower bound for loss area weighted by
distortion. To describe it, let $E$ be a certified subset of an
appropriate normalized loss set, and let $C_t^+$ cover the
normalized landing set $L_t^+$, up to a null set. The loss set
records points which eventually cease returning to the relevant
critical neighborhood, subject to the avoidance conditions
specified in Section~\ref{sec:pairing}. Define
\[
 g_E(z)=\int_E \overline C(z,w)^{-2}\,dA(w),
\]
where $\overline C(z,w)$ bounds the derivative ratios of the
univalent inverse branches used in the landing--loss decomposition.
Thus $g_E$ measures the amount of certified loss available to
control a landing copy at the position $z$.

With the admissibility conditions and normalizations of
Theorem~\ref{thm:supported}, it is enough to establish
\[
 a_n\le A_n^+<\gamma,
 \qquad
 \mathop{\rm ess\,inf}_{z\in C_t^+}g_E(z)\ge\gamma.
\]
The corresponding lower-half-plane estimates follow by complex
conjugation. These inequalities yield a uniform recurrence
\[
 a_{k+n-1}\le\theta a_k,
 \qquad
 \theta=\frac{A_n^+}{\gamma}<1,
\]
for all sufficiently large $k$. The uniformity is essential:
the same support $C_t^+$ contains every deeper landing set,
and the same certified loss set remains admissible.
Iteration therefore gives exponential decay of the landing areas.

For the quartic map, we verify this criterion with
$(n,m_0,t)=(3,2,5)$. The certified bounds give
\[
 \theta=\frac{2128573712302080}{2160534643703221}<1,
 \qquad
 a_{k+2}\le\theta a_k\quad(k\ge5).
\]
The complete hypotheses and exact numerical bounds are stated in
Proposition~\ref{prop:certificate-intro}. Together with
Theorem~\ref{thm:supported}, they prove
Theorem~\ref{thm:main}.

Two features of the verification are worth emphasizing.
First, the numerical approximation of the renormalization map is
separated from its rigorous validation and canonical identification.
A finite Newton iteration produces an approximate solution.
Interval estimates then validate an exact analytic solution,
including its Taylor tail. An independent construction of proper
branch domains and an all-depth return argument identify this
solution with the canonical Fibonacci map; see
Sections~\ref{sec:validation} and~\ref{sec:identification}.

Second, both the landing and loss certificates combine
Koebe estimates for inverse branches with direct forward propagation
of whole cells. For landing, successive refinements replace only
cells retained by the preceding stage. For loss, independently
certified witnesses are combined by an exact geometric union,
so that overlapping regions contribute only once.
The resulting landing cover and loss union are then used to
bound the distortion integral. This separates the dynamical
certification of the sets from the subsequent kernel estimate;
see Sections~\ref{sec:algorithms} and~\ref{sec:certificate}.

\begin{remark}[Related results on the measure of Julia sets]
The Fibonacci problem belongs to a broader study of zero- and
positive-area Julia sets. Beyond Lyubich's theorem for non-renormalizable quadratic polynomials,
Petersen and Zakeri proved zero area for quadratic Siegel
polynomials whose rotation numbers satisfy
$\log a_n=O(\sqrt n)$, where $a_n$ are the continued-fraction
entries \cite{PZ}. Graczyk and Smirnov obtained zero-area and
dimension estimates under suitable summability conditions
on derivatives along critical orbits \cite{GS}.

Positive-area Julia sets occur even for quadratic polynomials.
Buff and Ch\'eritat constructed examples with a Siegel disk,
with a Cremer point, and with infinite satellite renormalization
\cite{BC}. Avila and Lyubich subsequently constructed
positive-area Feigenbaum Julia sets, including examples with
stationary primitive combinatorics \cite{ALArea}.
For rational maps, Yang obtained nowhere-dense positive-area
Julia sets without irrationally indifferent periodic points,
Herman rings, or polynomial-like renormalizations \cite{Yang}.

Within the Fibonacci family, sufficiently high even critical
degree is expected to produce positive-area Julia sets;
see the discussion in \cite{ALArea}. The earlier positive-area
claim of Nowicki and van Strien \cite{NS} contains a gap,
as noted in the arXiv record, and does not settle this question.
There is also a related real measure problem: wild Cantor
attractors were established for sufficiently flat Fibonacci
maps by Bruin, Keller, Nowicki and van Strien \cite{BKNS},
and explicit computer-assisted criteria were developed by
Dudko and Gaidashev \cite{DG}. These results concern the
one-dimensional measure of an attracting basin and do not
by themselves determine the planar area of the complex
Julia set.
\end{remark}

The paper is organized as follows. Part~I develops the analytic
criterion: it recalls the canonical Fibonacci map, establishes
univalent continuation for the relevant inverse copies,
constructs the disjoint landing--loss decomposition, and proves
the uniform area recurrence. Part~II treats the quartic map,
from analytic validation and canonical identification to
domain enclosures, landing and loss certification, and the
final kernel estimate. Section~\ref{sec:reproducibility}
describes the verification code and gives an ordered procedure
for reproducing the certificate.

\paragraph{AI Acknowledgement.}The main content of the manuscript is generated by ChatGPT 5.6 Sol and ChatGPT 6 Astra under a multi-turn human-AI interaction. The author takes full responsibility for the results in the paper.

\part*{\Large{Part I. Analytic criterion for even critical degree $\ell\ge4$}}
\addcontentsline{toc}{part}{Part I. Analytic criterion}

% BEGIN sections/02_preliminaries.tex
\section{Preliminaries: canonical maps, domains and puzzles}\label{sec:preliminaries}
Sections~\ref{sec:preliminaries}--\ref{sec:criterion} concern any even
critical degree $\ell\ge4$. The numerical application, starting in
Section~\ref{sec:validation}, has $\ell=4$. We collect the known dynamical
inputs here, so that the physical map and the domains used later are fixed
before the area argument begins.

\subsection{Normalization and physical iteration}
Let $f$ be the canonical real-symmetric solution, in the critical-maximum
normalization,
\begin{equation}\label{eq:functional}
 f(z)=-\alpha^{-2}f\bigl(\alpha f(\alpha z)\bigr),\qquad
 f(0)=1,\qquad 0<\alpha<1,
\end{equation}
with $f(z)=F(z^\ell)$ and $F'(0)<0$ near zero. Put
\begin{equation}\label{eq:canonical-auxiliary}
 q_\alpha(z)=f(\alpha z),\quad G(z)=f(\alpha^2z),\quad
 \HH(x)=F(\alpha^\ell x),\quad q_\alpha(z)=\HH(z^\ell).
\end{equation}
A holomorphic map is \emph{proper} if inverse images of compact sets are
compact. A two-branch generalized polynomial-like map
$D:U^0\sqcup U^1\to V$ has bounded quasidisk domains, disjoint source
closures compactly contained in $V$, and proper onto branches. Here the
central degree is $\ell$, the other branch is univalent, and the only
critical point is $0$; its orbit does not escape.

All iterations of $D$ below are \emph{physical}: at each step the point
must lie in one of its declared source domains. More precisely,
\begin{gather*}
 \Dom(D^0)=V,\qquad
 \Dom(D^{k+1})=\{z\in U^0\sqcup U^1:Dz\in\Dom(D^k)\},\\
 K(D)=\bigcap_{k\ge0}\Dom(D^k),\qquad J(D)=\partial K(D).
\end{gather*}
The notation $D^{-k}(A)$ always denotes the inverse image inside
$\Dom(D^k)$. A continued analytic formula outside a physical source is
not thereby an iterate of $D$.

\subsection{Canonical geometry used in the proof}
The following facts are from Buff~\cite{Buff}.

\paragraph{P1. The slit preimage.}
The inverse germ $\iota=\HH^{-1}$ at $1$ extends univalently to
\begin{equation}\label{eq:epstein-slit}
 \Delta_\alpha=\C\setminus
 ((-\infty,-\alpha^{-2}]\cup[\alpha^{-4},\infty)).
\end{equation}
It is anti-Herglotz: $\Im\iota(w)<0$ when $\Im w>0$. Define
\begin{equation}\label{eq:epstein-source}
 E_\alpha=\{z:z^\ell\in\iota(\Delta_\alpha)\}.
\end{equation}
This is Buff's bounded domain $W_\alpha$; we rename it to distinguish it
from our homothetic tower $W_j$. The map
$q_\alpha:E_\alpha\to\Delta_\alpha$ is proper of degree $\ell$, with
sole critical point $0$ \cite[Section 3.3, Theorem 1; Section 4.3]{Buff}.

\paragraph{P2. The attracting basin and the range.}
Let $\rho$ be the first positive zero of $q_\alpha$. It is an attracting
fixed point of $G$. Its immediate basin $B$ is a quasidisk,
$\overline B=K(G)\Subset E_\alpha/\alpha$, and
$G:B\to B$ is proper of degree $\ell$. The physical range is
\begin{equation}\label{eq:canonical-V}
 V=\alpha B,\qquad \alpha E_\alpha\subset V\Subset E_\alpha.
\end{equation}
The maximal domains of $f$ and $q_\alpha$ are $V$ and $B$, respectively;
on $B$ one has
\begin{equation}\label{eq:canonical-linearization}
 q_\alpha\circ G=-\alpha^2q_\alpha.
\end{equation}
Also $q_\alpha:V\to B$ is proper of degree $\ell$. The domains $B,V$
are invariant under conjugation and multiplication by $\ell$th roots
of unity \cite[Section 4.3, Theorem B and Lemmas 6--7]{Buff}.

\paragraph{P3. The two physical branches.}
Set $A_{-1}=V$, $A_j=(G|_B)^{-(j+1)}(V)$ for $j\ge0$, and
\begin{equation}\label{eq:canonical-sources}
 Z=\alpha V=\alpha^2B,\qquad U^0=\alpha^2A_0,\qquad
 U^1=\operatorname{Comp}_{1}(q_\alpha|_V)^{-1}(Z).
\end{equation}
The domains $Z,U^1$ are the immediate basins of a two-cycle for
$q_\alpha$; they are disjoint, $U^0\Subset Z$, $U^1\Subset V$,
$q_\alpha(Z)=U^1$, and $q_\alpha:U^1\to Z$ is biholomorphic.
The resulting physical map is
\begin{equation}\label{eq:physical-D}
 D(z)=\begin{cases}f(z),&z\in U^0,\\
 \alpha^{-1}f(\alpha z),&z\in U^1.
 \end{cases}
\end{equation}
Its critical value is $1\in U^1$. We use this \emph{positive lateral
sign} throughout \cite[Section 4.4, Theorem C, Steps 1--4]{Buff}.

\paragraph{P4. Renormalization germs and uniqueness.}
Two generalized polynomial-like maps define the same germ if their filled
sets agree and the maps agree on that set \cite[Section 2.4, Definition 7]{Buff}.
For each even $\ell\ge4$ there is a unique real Fibonacci renormalization
cycle of period two. Buff states this as Theorem 3 in Section 4.2 and
explicitly attributes it to van Strien--Nowicki, Theorem 7.1
\cite{NS}. We use this cycle theorem through Buff's formulation, after
constructing a map in its stated class. The critical-maximum normalization
and the sign $D^2(0)=-\alpha<0$ select the phase used here.
The separate positive-area claim in the Nowicki--van Strien preprint has
a known gap; that claim is not an input to this paper.

\begin{figure}[htbp]
 \centering\includegraphics[width=0.6\textwidth]{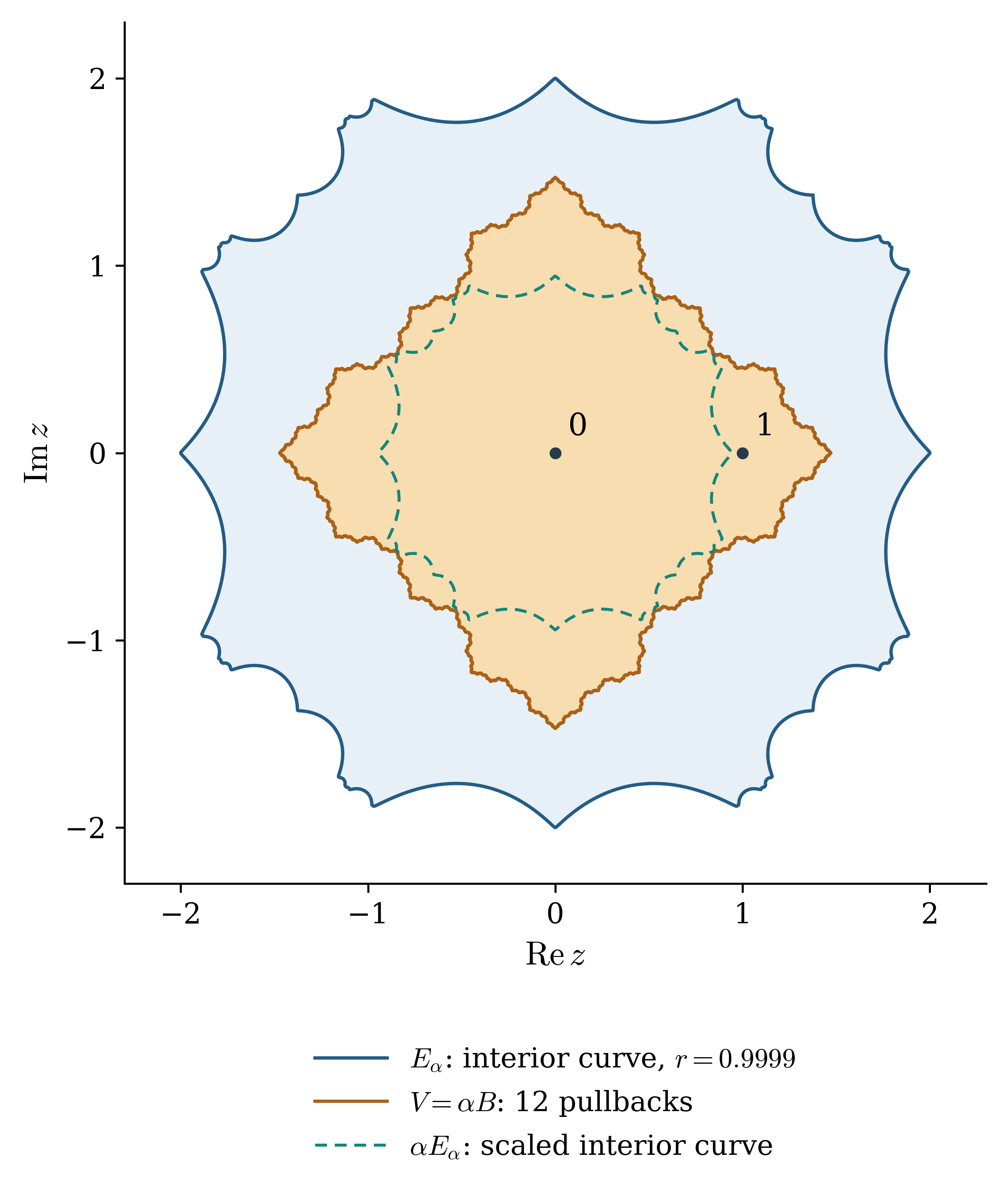}
 \caption{Numerical illustration of domains $E_\alpha$ and
 $V=\alpha B$ for the quartic solution. For the blue interior curve,
 $r=0.9999$ is the radius of an auxiliary circle $|\zeta|=r$ in a
 unit-disk coordinate on $\Delta_\alpha$. The orange
 contour approximates $V$ by twelve pullbacks; the dashed green curve
 is the $\alpha$-scaled blue curve. Rigorous enclosures are constructed
 separately in Section~\ref{sec:geometry}.}\label{fig:buff-domains}
\end{figure}
\begin{figure}[htbp]
 \centering\includegraphics[width=\textwidth]{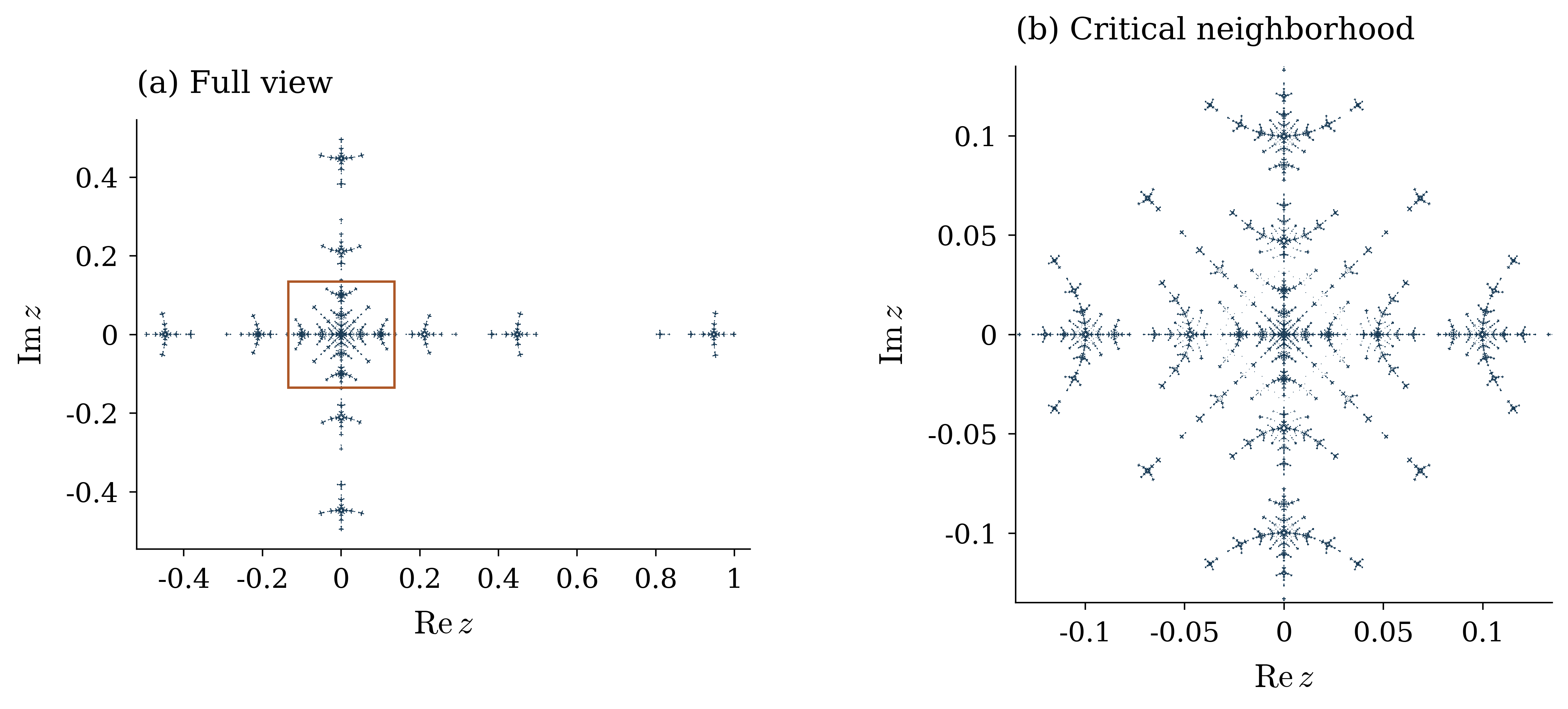}
 \caption{Inverse-iteration illustration of the quartic Fibonacci Julia
 set, with an enlarged critical neighborhood on the right. The plotting
 procedure and its physical iteration depths are specified in
 Section~\ref{sec:illustrations}.}\label{fig:julia}
\end{figure}

\begin{figure}[htbp]
\centering
\begin{tikzpicture}[manuscript diagram]
  % A single range domain: both mapping arrows return to this same V.
  \draw[diagramblue, fill=diagramblue!3, line width=.75pt]
    (0,0) ellipse[x radius=5.2cm,y radius=1.9cm];
  \draw[gray!65,dashed,fill=gray!7,line width=.55pt]
    (-1.6,0) ellipse[x radius=2.2cm,y radius=1.28cm];
  \draw[diagramblue,fill=diagramblue!7,line width=.65pt]
    (-1.6,0) ellipse[x radius=1.35cm,y radius=.83cm];
  \draw[diagramblue,fill=diagramblue!17,line width=.65pt]
    (-1.6,0) ellipse[x radius=.72cm,y radius=.43cm];
  \draw[diagramorange,fill=diagramorange!12,line width=.65pt]
    (2.55,0) ellipse[x radius=.78cm,y radius=.95cm];
  \node at (-4.35,.62) {$V$};
  \node at (-2.92,.91) {$\alpha V$};
  \node at (-2.53,.40) {$U^0$};
  \node at (-1.60,.24) {$W_1$};
  \node at (2.55,.52) {$U^1$};
  \fill (-1.6,-.07) circle[radius=1.1pt];
  \node[font=\scriptsize,below right=1pt] at (-1.6,-.07) {$0$};
  \fill (2.55,-.07) circle[radius=1.1pt];
  \node[font=\scriptsize,below right=2pt] at (2.55,-.07) {$1$};
  \draw[map arrow,diagramblue] (-1.60,.85)
    .. controls (-2.6,3.10) and (1.4,3.20) .. (2.80,1.60);
  \node[map label,text=diagramblue] at (.05,2.84)
    {$f$\quad (degree $\ell$)};
  \draw[map arrow,diagramorange] (3.34,-.05)
    .. controls (5.5,-.6) and (4.2,-3.05) .. (1.50,-1.82);
  \node[map label,text=diagramorange] at (3.35,-2.47)
    {$\alpha^{-1}q_\alpha$\quad (degree $1$)};
\end{tikzpicture}
\caption{Schematic geometry of the physical map. Both sources are compactly contained in $V$, with $W_1\Subset U^0\Subset\alpha V$ and $U^1\cap\alpha V=\varnothing$. Both arrows map onto the same outer domain $V$, with degrees $\ell$ and one.}\label{fig:physical}
\end{figure}
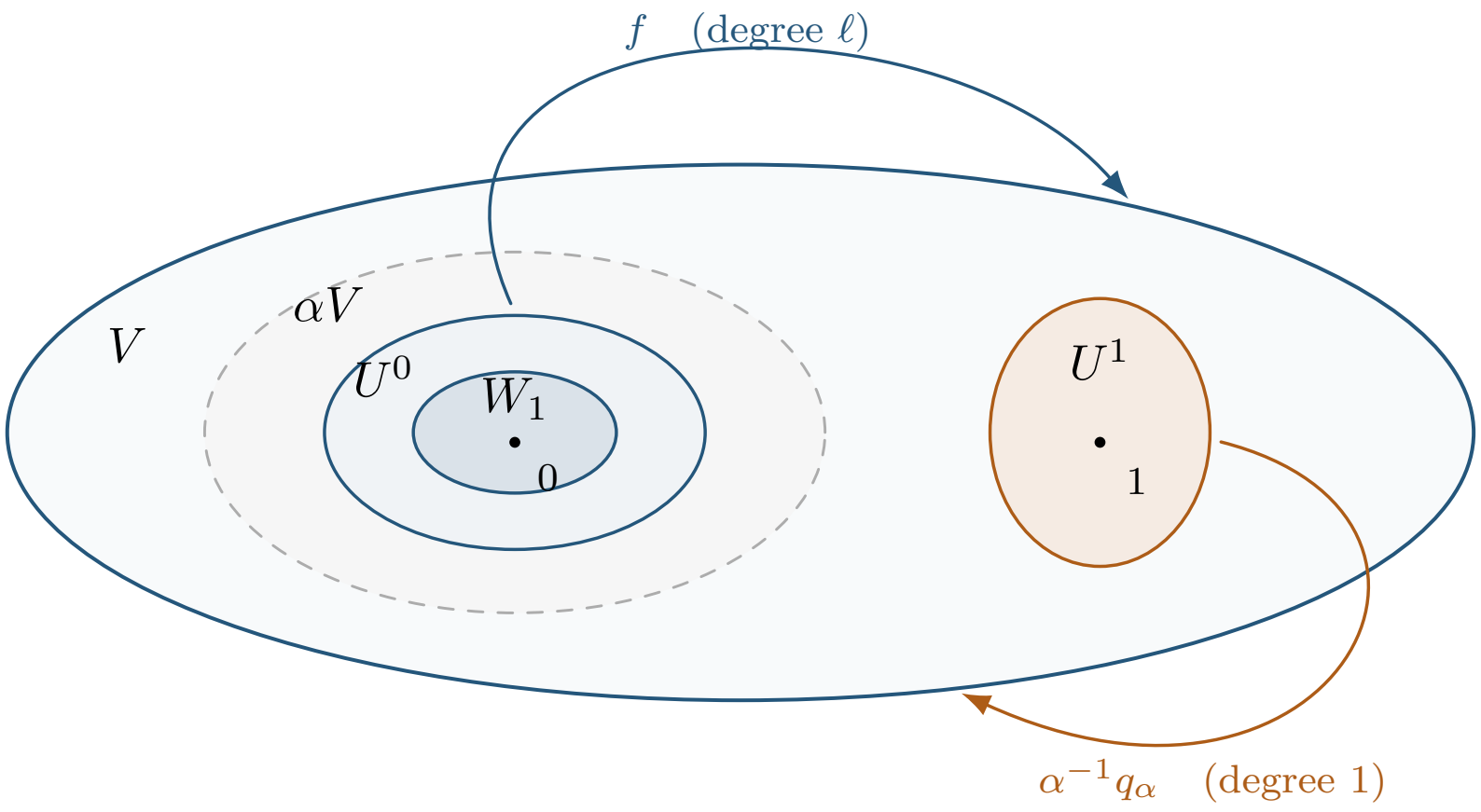
\subsection{Puzzles and Fibonacci combinatorics}
Write $c_r=D^r(0)$ and $\Pc(D)=\overline{\{c_r:r\ge1\}}$. For
$x\in\Dom(D^{d+1})$ let
\[
 P_d(x)=\operatorname{Comp}_x\Dom(D^{d+1}),\quad
 C_d=P_d(0),\qquad C_{-1}=V.
\]
These finite puzzle pieces satisfy the Markov property: physical
pullbacks are puzzle pieces, and two pieces are disjoint or nested
according to their depth. A deeper piece is compactly contained in the
shallower piece whenever they meet and their depths differ.
Use the Fibonacci times
\begin{equation}\label{eq:fib-times}
 S_0=1,\qquad S_1=2,\qquad S_{r+1}=S_r+S_{r-1}.
\end{equation}
The critical itinerary is characterized by
\begin{equation}\label{eq:fibonacci-combinatorics}
 c_{S_r}\in C_{S_{r+1}-3}\setminus C_{S_{r+1}-2}
 \qquad(r\ge0),
\end{equation}
with $C_{-1}=V$ at the first index
\cite[Section 4.1, Proposition 7; Section 4.4, Step 4]{Buff}.
An itinerary specifies every central or lateral branch and, when an
inverse is taken, every inverse-sheet choice. These labels remain part
of the data even when the formulas are analytically continued.

% END sections/02_preliminaries.tex

% BEGIN sections/02_dynamics.tex
\section{The homothetic tower and its intrinsic slit}\label{sec:dynamics}
\subsection{The two towers and the index shift}

Let $U_0^j=U^j$ and $V_0=V$. For $r\ge1$, the consecutive
puzzle ranges and sources are
\[
 V_r=C_{S_{r+1}-3},\qquad U_r^0=C_{S_{r+2}-3},\qquad
 U_r^1=P_{S_{r+1}+S_{r-1}-3}(c_{S_r}).
\]
These are the domains in Buff's return construction
\cite[Section 4.1, Lemma 2]{Buff}. His consecutive return at level $r\ge1$ is
\begin{equation}\label{eq:consecutive-return}
 f_r:U_r^0\sqcup U_r^1\longrightarrow V_r=U_{r-1}^0,
 \qquad
 f_r|_{U_r^0}=D^{\circ S_r},\quad
 f_r|_{U_r^1}=D^{\circ S_{r-1}}.
\end{equation}
These ranges are not homothetic.  In particular, neither $f_r$ nor $V_r$ is
identified with Buff's domain $D_r$.

For the present argument set
\begin{equation}\label{eq:W-tower}
 h(z)=-\alpha^2z,\qquad h_m=h^{\circ m},\qquad
 W_m=h_m(V),\qquad m\ge0.
\end{equation}
When taking absolute values, $h_m$ also denotes its real linear multiplier;
thus $|h_m|=\alpha^{2m}$.
Set $C_{-1}=V$; for $d\ge0$, $C_d$ is the component of $\Dom(D^{\circ(d+1)})$ containing $0$. Buff's critical-piece formula \cite[Section 4.4, Lemmas 9--10]{Buff} gives
\begin{equation}\label{eq:W-puzzle}
 W_m=C_{S_{2m}-2},\qquad W_{m+1}\Subset W_m,\qquad
 \bigcap_{m\ge0}W_m=\{0\}.
\end{equation}
Define the scaled two-branch map
\begin{equation}\label{eq:R-m}
 \cR_m=h_m\circ D\circ h_m^{-1}:
 h_mU^0\sqcup h_mU^1\longrightarrow W_m.
\end{equation}
For $m\ge1$, Buff's iterate identities restrict to
\begin{equation}\label{eq:physical-times}
 \cR_m|_{h_mU^0}=D^{\circ S_{2m}},
 \qquad
 \cR_m|_{h_mU^1}=D^{\circ S_{2m-1}}.
\end{equation}
The equalities hold only on the displayed sources.

Put $K_m=K(\cR_m)=h_mK(D)$.  For a target $W_m$, $m\ge1$, define
\begin{equation}\label{eq:parent}
 \boxed{\ \cF_m:=\cR_{m-1},\qquad K_m^{\rm par}:=K(\cF_m)=K_{m-1}.\ }
\end{equation}
Write $J(D):=\partial K(D)$.

\begin{center}
\small
\begin{tabularx}{\textwidth}{@{}>{\bfseries}lXXX@{}}
\toprule
Object & Target/range & Physical times & Role for a copy of $W_m$\\
\midrule
$f_r$ & $V_r$ & $S_r,S_{r-1}$ & consecutive canonical tower\\
$\cF_m=\cR_{m-1}$ & $W_{m-1}$ & $S_{2m-2},S_{2m-3}$ $(m\ge2)$ & parent clock/core; $\cF_1=D$\\
$\cR_m$ & $W_m$ & $S_{2m},S_{2m-1}$ & last return to target\\
\bottomrule
\end{tabularx}
\end{center}

Thus the ``shift'' has two distinct sources.  First, the parent of the target
$W_m$ is indexed $m-1$.  Second, normalizing a landing from $W_n$ uses
$h_{n-1}^{-1}$, and consequently sends $W_{n+m}$ to $W_{m+1}$; this produces
the factor $\widetilde\eta_{m+1}$ in the recurrence.

\subsection{The physical return clock}

\begin{lemma}[Setwise first returns]\label{lem:first-return}
For every $m\ge1$, the two sources in \eqref{eq:physical-times} are full
components of the indicated physical pullbacks of $W_m$, and
\begin{align}
 D^{\circ r}(h_mU^0)\cap W_m&=\varnothing
       &&(0<r<S_{2m}),\label{eq:first-central}\\
 D^{\circ r}(h_mU^1)\cap W_m&=\varnothing
       &&(0<r<S_{2m-1}).\label{eq:first-lateral}
\end{align}
They are therefore two principal components of the complete first-return map
to $W_m$.  Completeness of this two-component list is not asserted.
\end{lemma}

\begin{proof}
At $m=1$, Buff's formulas give
$D^{\circ3}:hU^0\to W_1$ and $D^{\circ2}:hU^1\to W_1$ as proper onto maps.
Here $W_1=C_1\Subset U^0$, $D(W_1)=U^1$, and the first images of both sources lie in $U^1$, disjoint from $W_1$.
The only other intermediate image, $D^2(hU^0)$, is a depth-two piece and
contains $c_2=-\alpha\in C_0\setminus C_1$.  The puzzle Markov alternative
therefore excludes an intermediate intersection with $C_1$.  Properness
makes each displayed source an entire pullback component.

Assume the result at level $m$.  In $\cR_m$-coordinates the next target is
$W_{m+1}=h_m(W_1)$.  Conjugating the base return diagram gives the two
macro-words
\[
 0,1,0\qquad\text{and}\qquad0,1.
\]
The base diagram excludes a return at a proper macro-prefix, while the
induction hypothesis excludes a physical return inside a macro-step.  The
physical lengths are
\[
 2S_{2m}+S_{2m-1}=S_{2m+2},
 \qquad
 S_{2m}+S_{2m-1}=S_{2m+1}.
\]
Compositions of the proper branches are proper, so the same open-and-closed
component argument proves the claim at $m+1$.
\end{proof}

For $x\in K_m$, let $a_j(x)$ be the branch label of $\cR_m^j(x)$.
Set $t_m(0)=S_{2m}$, $t_m(1)=S_{2m-1}$ for $m\ge1$, and
$t_0(0)=t_0(1)=1$. Its cumulative physical times are
\[
 T_{m,0}(x)=0,\qquad T_{m,q}(x)=\sum_{j=0}^{q-1}t_m(a_j(x)).
\]
An \emph{admissible $\cR_m$-word} is a finite composition of the two branches
in \eqref{eq:R-m} for which every intermediate state lies in the declared
source.  Its \emph{physical clock time} is the sum of the corresponding
lengths in \eqref{eq:physical-times}.

\begin{corollary}[Synchronization]\label{cor:synchronization}
If $x\in K_m$, then the physical visits of the orbit of $x$ to $W_m$ occur
exactly at endpoints of admissible $\cR_m$-blocks.  Equivalently,
\[
 \{k\ge0:D^kx\in W_m\}=\{T_{m,q}(x):q\ge0\},\qquad
 D^{T_{m,q}(x)}x=\cR_m^q(x).
\]
All these physical iterates are defined.  In particular the visits of the critical
orbit to $W_m$ are
\[
 h_mc_0,h_mc_1,h_mc_2,\ldots,
 \qquad
 \Pc(D)\cap W_m=h_m\Pc(D).
\]
Moreover, for $m\ge1$,
\begin{equation}\label{eq:core-inclusion}
 K_m\subset K_{m-1}\cap W_m.
\end{equation}
\end{corollary}

\begin{proof}
Concatenate the blocks of Lemma~\ref{lem:first-return}.  Equations
\eqref{eq:first-central}--\eqref{eq:first-lateral} exclude every intermediate
visit.  The conjugacy in \eqref{eq:R-m} gives the displayed critical orbit;
taking closures gives the postcritical identity.  Finally, after conjugating
by $h_{m-1}^{-1}$, an orbit in $K(\cR_1)$ is an infinite concatenation of the
physical blocks $D^{\circ3}|_{hU^0}$ and $D^{\circ2}|_{hU^1}$.  All
intermediate iterates are defined, so the point belongs to $K(D)$.  Scaling
back proves \eqref{eq:core-inclusion}.
\end{proof}

\subsection{Nearest postcritical endpoints}

The first seven physical labels are $0,1,0,0,1,0,1$ and
\[
 c_1=1,\ c_2=-\alpha,\ c_3=-\alpha^2,\ c_4=f(\alpha^2),\
 c_5=\alpha^3,\ c_6=f(\alpha^3),\ c_7=-\alpha f(\alpha^2).
\]
Write
\begin{equation}\label{eq:xi}
 \xi=c_4=f(\alpha^2)>0,\qquad c_7=-\alpha\xi.
\end{equation}
Let $I_m=W_m\cap\R=(a_m,b_m)$.  Define intrinsically
\begin{equation}\label{eq:intrinsic-endpoints}
 p_m=\max\bigl(\Pc(D)\cap(-\infty,a_m)\bigr),\qquad
 q_m=\min\bigl(\Pc(D)\cap(b_m,\infty)\bigr)
\end{equation}
and
\begin{equation}\label{eq:intrinsic-slit}
 \boxed{\ \Omega_m=\C\setminus
 \bigl(( -\infty,p_m]\cup[q_m,\infty)\bigr).\ }
\end{equation}
This definition refers to the actual quasidisk $W_m$, not to a round
surrogate.

\begin{lemma}[Intrinsic endpoint identity]\label{lem:endpoints}
For every $m\ge1$, the extrema in \eqref{eq:intrinsic-endpoints} exist and
\begin{equation}\label{eq:endpoint-pair}
 \{p_m,q_m\}=\{h_{m-1}c_7,h_{m-1}c_4\}.
\end{equation}
Consequently
\begin{equation}\label{eq:stationary-slit}
 \Omega_m=h_m\Omega_*,\qquad
 \Omega_*=\C\setminus
 \bigl(( -\infty,-\xi/\alpha^2]\cup[\xi/\alpha,\infty)\bigr),
 \qquad V\Subset\Omega_*.
\end{equation}
\end{lemma}

\begin{proof}
Let $V\cap\R=(-r,r)$.  Since $c_1=1\in U^1\Subset V$, one has $r>1$.
Let $\rho$ be the first positive zero of $f_\alpha(z)=f(\alpha z)$.
The central real interval contains $\pm\alpha$, whereas the disjoint
lateral interval contains $\xi$ and $1$, so $\xi>\alpha$.
Since $q_\alpha(0)=1$, $q_\alpha(1)=-\alpha^2$, and
$q_\alpha(\alpha)=\xi>0$, monotonicity gives $\alpha<\rho<1$.
Buff's Lemma 6 gives $f(\alpha^2\rho)=\rho$; strict decrease of $f$
then gives $\xi=f(\alpha^2)<\rho$. Thus
\[
 0<\alpha<\xi<\rho<1,\qquad f(\alpha^2\rho)=\rho.
\]
The functional equation at $z=\alpha$ gives $f(\alpha\xi)=\alpha^4$,
and at $z=\xi$ it gives
\[
 f(\xi)=-\alpha^{-2}f(\alpha^5)<-\rho/\alpha^2,
\]
because $\alpha^5<\alpha^2\rho$ and $f$ decreases on the positive
real Epstein interval. Also $\xi>\alpha\rho$.
The immediate invariant basin $\alpha^{-1}V$ of $z\mapsto f(\alpha^2z)$ is
compactly contained in the Epstein source whose real trace is
$(-\rho/\alpha^2,\rho/\alpha^2)$.  At $t=\xi/\alpha^2$ the image is
$f(\xi)<-\rho/\alpha^2$, so $t\notin\overline{\alpha^{-1}V}$. Indeed, if $t$ belonged to that closure, continuity and forward invariance would put its image in the closure of the basin, which is contained in the displayed source. The strict image inequality excludes this. Hence
\begin{equation}\label{eq:r-bound}
 1<r<\xi/\alpha.
\end{equation}

The real Fibonacci order in Dudko--Gaidashev
\cite[Lemma 3.1 and its proof]{DG}, in the critical-maximum normalization
used here, gives
\[
 \Pc(D)\subset[-\alpha,-\alpha\xi]\,
       \cup[-\alpha^2,\alpha^2]\,\cup[\xi,1].
\]
This is the three-cluster consequence of the five postcritical intervals
in that proof; their critical-minimum orbit coordinates change by
$x_j/s=-c_j$, where $s=-x_1>0$. Both $-\alpha\xi=c_7$ and $\xi=c_4$ belong to the orbit.
Since $r>1$ and $r<\xi/\alpha$, the middle cluster is contained in $W_1$
and
\[
 -\alpha\xi< -\alpha^2r<\alpha^2r<\xi.
\]
The two displayed orbit points are therefore exactly the nearest
postcritical points outside $W_1$.

Corollary~\ref{cor:synchronization} transports the order inside $W_{m-1}$ by
$h_{m-1}$; postcritical points outside $W_{m-1}$ are farther from $W_m$.
This proves \eqref{eq:endpoint-pair}.  Dividing the unordered pair by $h_m$
gives $\{-\xi/\alpha^2,\xi/\alpha\}$, proving the slit identity.  Finally,
\eqref{eq:r-bound} puts the real trace of $V$ strictly inside that of
$\Omega_*$; all nonreal points belong to the slit plane. Since $V$ is a real-symmetric Jordan domain, $\overline V\cap\R=[-r,r]$. Compactness yields
$V\Subset\Omega_*$.
\end{proof}

% END sections/02_dynamics.tex

% BEGIN sections/03_branches.tex
\section{Physical pullbacks and labelled inverse branches}\label{sec:branches}

The convention for $D^{-k}$ is the physical inverse-image convention
of Section~\ref{sec:preliminaries}. Set
\[
 H^+=\{z:\operatorname{Im}z>0\},\qquad
 H^-=\{z:\operatorname{Im}z<0\},\qquad
 P_m^\sigma=W_m\cap H^\sigma.
\]

\subsection{Full puzzle pullbacks}

For a physical iterate let $\Dom(D^{\circ k})$ be the open set on which all
first $k$ iterates follow physical branches.  We use the following standard
consequence of puzzle construction.

\begin{lemma}[Full-hull Markov alternative]\label{lem:markov}
Let $T$ be a finite critical puzzle piece, in particular some $W_m$, and let
$A$ be a component of
\[
 \Dom(D^{\circ k})\cap D^{-k}(T).
\]
Then $D^{\circ r}(A)$ is the corresponding full component of the remaining
pullback for $0\le r\le k$.  These are finite puzzle pieces: if $T$ has depth $d$, their depths are $d+k-r$.  Two such pieces
are disjoint or nested in the direction determined by depth; an intersecting deeper piece is compactly contained in the shallower one.
In particular, if $A$ contains $0$, then $A\subset T$.  Finally,
\begin{equation}\label{eq:postcritical-boundary}
 \Pc(D)\cap\partial T=\varnothing.
\end{equation}
\end{lemma}

\begin{proof}
On a fixed physical itinerary, properness at each step shows that a pullback
component maps properly and onto its terminal component.  The itinerary is a
puzzle cylinder; pulling back increases its depth, so Markov property
gives the asserted nesting direction.  A deeper component and $T$ which both
contain $0$ must therefore be nested as $A\subset T$.

The filled set $K(D)$ is compact and every point of it follows, at each finite
depth, an itinerary through the interiors of the two physical sources.
Consequently it lies in the interior of its finite cylinder and on no finite
puzzle boundary.  Since $\Pc(D)\subset K(D)$, \eqref{eq:postcritical-boundary}
follows.
\end{proof}

\begin{definition}[Tagged copy and full hull]\label{def:tagged}
A tagged copy $Q$ of $P_m^\sigma$ under $D^{\circ k}$ is a component of the physical inverse image $D^{-k}(P_m^\sigma)$ for which
\[
 D^{\circ k}:Q\longrightarrow P_m^\sigma
\]
is biholomorphic.  Its data include the entire central/lateral word and, at
each central inverse, the chosen root sector.  Its \emph{full hull} $\wh Q$ is
the component of
\[
 \Dom(D^{\circ k})\cap D^{-k}(W_m)
\]
containing $Q$.  Put $\wh Q_r=D^{\circ r}(\wh Q)$ and
$Q_r=D^{\circ r}(Q)$.
\end{definition}

\begin{definition}[Primitive and parent-core-separated]\label{def:eligible}
The tagged copy is \emph{primitive} if
\begin{equation}\label{eq:primitive}
 \wh Q_r\cap W_m=\varnothing\qquad(0\le r<k).
\end{equation}
If it is not primitive, put
\[
 r_*(Q)=\max\{r<k:\wh Q_r\cap W_m\ne\varnothing\},
 \qquad B(Q)=\wh Q_{r_*(Q)}.
\]
By Lemma~\ref{lem:markov}, $B(Q)\subset W_m$.  The copy is
\emph{parent-core-separated} if
\begin{equation}\label{eq:separated}
 B(Q)\cap K_m^{\rm par}=\varnothing,
 \qquad K_m^{\rm par}=K_{m-1}.
\end{equation}
We call a primitive or parent-core-separated copy \emph{eligible}.
\end{definition}

The parent index $m-1$ in \eqref{eq:separated} is essential: a finite
admissible inverse of $0$ under $\cR_{m-1}$ belongs to $K_{m-1}$,
but need not belong to $K_m$. The pairing contradiction therefore
requires separation from the parent core.

\begin{lemma}[Synchronization and tail stability]\label{lem:copy-sync}
Let $1\le n\le m$, and let $\wh Q\subset W_{n-1}$ be the full hull of a
tagged copy of $P_m^\sigma$ under $D^{\circ k}$.  If
$\wh Q\cap K_{n-1}\ne\varnothing$, then $k$ is the physical time of an
admissible $\cR_{n-1}$-word and that word agrees with $D^{\circ k}$ on
$\wh Q$.

If $Q$ is parent-core-separated at level $m$, then
\begin{equation}\label{eq:separated-hull}
 \wh Q\cap K_{m-1}=\varnothing.
\end{equation}
Every forward tail of an eligible copy, with the remaining labels and root
sectors, is again eligible for the same target.
\end{lemma}

\begin{proof}
Choose $x\in\wh Q\cap K_{n-1}$.  Since its terminal point lies in
$W_m\subset W_{n-1}$, Corollary~\ref{cor:synchronization} says that $k$ is a
$\cR_{n-1}$ clock time.  The admissible induced pullback of $W_m$ through $x$
and $\wh Q$ are full components of the same physical pullback and hence agree.

For \eqref{eq:separated-hull}, suppose instead that
$x\in\wh Q\cap K_{m-1}$.  Since $K_{m-1}\subset W_{m-1}$, full-hull nesting
first gives $\wh Q\subset W_{m-1}$, so the preceding synchronization statement
applies both at time $k$ and at $r_*(Q)$.  Invariance at parent-clock times
then gives $D^{\circ r_*(Q)}(x)\in B(Q)\cap K_{m-1}$, contrary to
\eqref{eq:separated}.  Finally, a tail taken after $r_*(Q)$ is primitive; a
tail taken at or before $r_*(Q)$ has the same last-return hull $B(Q)$ and is
separated.  Primitive tails are primitive directly from
\eqref{eq:primitive}.
\end{proof}

\subsection{The two labelled Epstein surfaces}

Let $f_\alpha(z)=f(\alpha z)$. It follows from H. Epstein's work \cite{E} that $f_\alpha$ is a degree-$\ell$ proper map over a Jordan domain
$E_1=E_\alpha$:
\[
 f_\alpha:E_1\longrightarrow
 \Delta_0:=\C\setminus
 \bigl(( -\infty,-\alpha^{-2}]\cup[\alpha^{-4},\infty)\bigr),
\]
ramified only at $0$.  Put $E_0=\alpha E_1$ and define
\begin{equation}\label{eq:analytic-branches}
 F_0=f:E_0\to\Delta_0,
 \qquad
 F_1(z)=\alpha^{-1}f(\alpha z):E_1\to
 \Delta_1:=\alpha^{-1}\Delta_0.
\end{equation}
Both analytic maps have degree $\ell$ and sole critical point $0$; their
restrictions to $U^0,U^1$ are the physical formulas in
\eqref{eq:physical-D}.  If $\rho$ is the first positive zero of $f_\alpha$,
then
\[
 \overline{E_0}\cap\R=[-\rho,\rho],\qquad
 \overline{E_1}\cap\R=[-\rho/\alpha,\rho/\alpha],\qquad
 \cl{E_0\cup E_1}\Subset\Delta_0\cap\Delta_1.
\]
The trace formulas are those in \cite[Section 4.3, proof of Lemma 6]{Buff}.
Boundedness and $\rho<1$ give the compact containment. Since $p_m,q_m$
are in $[-1,1]$, also $\Omega_m\subset\Delta_0\cap\Delta_1$ for $m\ge1$.

For a label word ${\bf a}=(a_0,\ldots,a_{k-1})\in\{0,1\}^k$, define its
compatibility domain by requiring
\[
 z_0\in E_{a_0},\qquad z_{r+1}=F_{a_r}(z_r),\qquad
 z_r\in E_{a_r}\quad(1\le r<k),
\]
and let
\begin{equation}\label{eq:word-map}
 G_{\bf a}=F_{a_{k-1}}\circ\cdots\circ F_{a_0}
\end{equation}
there.  Labels are retained outside the physical domains; in general
$G_{\bf a}$ is not the physical iterate $D^{\circ k}$ there.

\begin{lemma}[Proper word and crossing criterion]\label{lem:crossing}
For each $m\ge1$,
\[
 G_{\bf a}:G_{\bf a}^{-1}(\Omega_m)\longrightarrow\Omega_m
\]
is finite and proper, where the inverse image is relative to the compatibility
domain.  At a compatible tuple,
\begin{equation}\label{eq:chain-rule}
 G_{\bf a}'(z_0)=\prod_{r=0}^{k-1}F_{a_r}'(z_r).
\end{equation}
Thus every critical value is the terminal value of a composable suffix which
starts at the labelled critical point $0$, and every such value is real.

Let a tagged copy select the inverse sheet
$\psi_\sigma:H^\sigma\to C_\sigma$ agreeing with the physical inverse
on $P_m^\sigma$. For $b\in(p_m,q_m)$ the limit
\[
 \zeta_b=\lim_{z\to b,\ z\in H^\sigma}\psi_\sigma(z),\qquad
 z_0(b)=\zeta_b,\quad z_{r+1}(b)=F_{a_r}(z_r(b))
\]
exists and gives a compatible tuple ending at $z_k(b)=b$.
The inverse extends univalently to $\Omega_m$ if and only if
\[
 G_{\bf a}'(\zeta_b)\ne0
 \quad\text{for every }b\in\CritVal(G_{\bf a})\cap(p_m,q_m),
 \quad\text{equivalently}\quad\prod_{r=0}^{k-1}z_r(b)\ne0.
\]
This concerns the selected point of the fiber; the other points need not
be regular. Failure gives $b\in(p_m,q_m)$ and $0\le r<k$ with $z_r(b)=0$.
\end{lemma}

\begin{proof}
For a compact $L_k\Subset\Omega_m$, pull it backward recursively by the
proper maps $F_{a_{k-1}},\ldots,F_{a_0}$.  Compact containment of both sources
in both ranges makes every intermediate inverse compact in the required
source.  This proves properness; finiteness follows from the degree bound
$\ell^k$.  Equation \eqref{eq:chain-rule} gives the critical-value statement.

Over either half-plane there are no critical values.  The selected component
is therefore a proper unramified cover of a simply connected domain, hence a
biholomorphism.  The local normal form of a finite proper holomorphic map at a fiber point is $u\mapsto u^d$ in suitable coordinates. Applied over a sufficiently small upper or lower half-disk, this shows that the selected sheet has a unique limiting point over every point of $(p_m,q_m)$.  A regular fiber extends by the inverse function theorem; a
critical fiber cannot be crossed by a holomorphic right inverse because
differentiating $G_{\bf a}\circ\psi=\mathrm{id}$ would give $0=1$.
After the finitely many word-critical values are deleted, the word map is a
finite covering.  Lollipop loops with stems in $H^\sigma$ generate the
fundamental group; a regular continuation across each puncture kills the
corresponding monodromy.  The sheet is therefore single-valued on
$\Omega_m$, and the right-inverse identity makes it injective.
\end{proof}

\subsection{Recognition of a physical critical word}

Let $\iota_r\in\{0,1\}$ be the physical label of $c_r$.  The real Fibonacci
orbit order gives
\begin{equation}\label{eq:orbit-cores}
 c_r\in U^0\iff c_r\in[-\alpha,\alpha^2],
 \qquad
 c_r\in U^1\iff c_r\in[\xi,1].
\end{equation}

\begin{lemma}[Sheet recognition]\label{lem:recognition}
Suppose a labelled evaluation
\[
 y_0=0,\qquad y_{r+1}=F_{a_r}(y_r)\quad(0\le r<k)
\]
is compatible through time $k$ and
\begin{equation}\label{eq:recognition-window}
 y_k\in(-\xi/\alpha,\xi/\alpha).
\end{equation}
Then $a_r=\iota_r$ for $0\le r<k$ and $y_r=c_r$ for $0\le r\le k$.
\end{lemma}

\begin{proof}
Let $r$ be the first mismatch, so $y_r=c_r$.  If $\iota_r=0$ but $a_r=1$,
then $|y_r|\le\alpha$ and monotonicity on the real Epstein branch gives
\[
 y_{r+1}=\alpha^{-1}f(\alpha|y_r|)\ge\xi/\alpha>\rho.
\]
If $k=r+1$, this contradicts \eqref{eq:recognition-window}.
If $k\ge r+2$, the next central label is impossible because
$y_{r+1}\notin E_0$. If the next lateral label is defined, then
$u=\alpha y_{r+1}\in[\xi,\rho)$, and
\[
 y_{r+2}=\alpha^{-1}f(u)\le\alpha^{-1}f(\xi)<-\rho/\alpha^3.
\]
For $k=r+2$ this violates the terminal window; for $k>r+2$ it violates
compatibility, since $y_{r+2}\notin E_0\cup E_1$. The additional step is
needed precisely when the mismatch is before the terminal step.

If $\iota_r=1$ but $a_r=0$, then $y_r\in[\xi,1]$.  For $y_r\ge\rho$ the
central formula is outside its source; otherwise
\[
 y_{r+1}=f(y_r)\le f(\xi)<-\rho/\alpha^2,
\]
again outside both sources and the terminal window.  Both mismatches are
impossible.
\end{proof}

% END sections/03_branches.tex

% BEGIN sections/04_extension.tex
\section{Univalent extension on the homothetic tower}\label{sec:extension}

\begin{proposition}[Intrinsic-slit continuation for the $W_m$-domains]
\label{prop:main}
Let $m\ge1$, $\sigma\in\{+,-\}$, and let $Q$ be an eligible tagged copy of
$P_m^\sigma$ under $D^{\circ k}$ in the sense of
Definitions~\ref{def:tagged}--\ref{def:eligible}.  Then
\[
 \phi_Q=(D^{\circ k}|_Q)^{-1}:P_m^\sigma\longrightarrow Q
\]
extends to a univalent holomorphic map
\begin{equation}\label{eq:prop-extension}
 \Phi_Q:\Omega_m\longrightarrow\C.
\end{equation}
After normalization by $h_m$, every level uses the fixed slit $\Omega_*$ in
\eqref{eq:stationary-slit}.
\end{proposition}

\begin{proof}
Fix $m$.  Suppose that the assertion fails and choose an eligible
counterexample of minimum physical length $k$.  Let ${\bf a}$ be its complete
label/root-sector data.  Lemma~\ref{lem:crossing} produces a real crossing
value $b\in(p_m,q_m)$ and a compatible limiting tuple
$(z_0^*,\ldots,z_k^*)$, with $z_k^*=b$, in which some $z_s^*=0$.

\emph{The collision is at the initial coordinate.}
If $s>0$, the forward tail $Q_s$ is eligible by
Lemma~\ref{lem:copy-sync}.  Its selected inverse tends to $0$ over $b$.  A
holomorphic extension would satisfy a right-inverse identity for the suffix
word, whose first factor is critical at $0$; differentiation at $b$ would give
$1=0$.  Hence the tail is a shorter counterexample, contradicting minimality.
Therefore
\begin{equation}\label{eq:initial-collision}
 z_0^*=0,\qquad G_{\bf a}(0)=b.
\end{equation}

\emph{The labelled word is physical.}
The interval $(p_m,q_m)$ lies in the recognition window
$(-\xi/\alpha,\xi/\alpha)$: its endpoints are postcritical and hence lie in
$[-1,1]$, while $\xi>\alpha$.  Lemma~\ref{lem:recognition} applied to
\eqref{eq:initial-collision} gives
\begin{equation}\label{eq:physical-collision}
 b=D^{\circ k}(0)=c_k.
\end{equation}
This is the step that excludes a ghost critical value created by analytic
continuation of the lateral formula.

\emph{The collision lies in the target.}
By the intrinsic definition of $p_m,q_m$, no postcritical point lies in either
gap $(p_m,a_m)$ or $(b_m,q_m)$.  Hence
$b\in\cl{W_m}\cap\R$.  Equation \eqref{eq:postcritical-boundary} and
\eqref{eq:physical-collision} upgrade this to $b\in W_m\cap\R$.  Approaching
$b$ through $P_m^\sigma$, the inverse values remain in $Q$ and converge to
$0$; thus
\begin{equation}\label{eq:zero-in-closure}
 0\in\cl Q.
\end{equation}

\emph{The last physical return contradicts eligibility.}
The physical critical itinerary stays in the interiors of its finite puzzle
cylinders, so $D^{\circ k}$ is defined on a neighborhood of $0$.  Let $A$ be
the component of
\[
 \Dom(D^{\circ k})\cap D^{-k}(W_m)
\]
which contains $0$.  A small disk about $0$ lies in $A$, and
\eqref{eq:zero-in-closure} gives points of $Q$ in that disk.  Component
uniqueness therefore gives $A=\wh Q$.  Lemma~\ref{lem:markov} gives
$\wh Q\subset W_m$, so $Q$ is not primitive.

Let $t=r_*(Q)$. For every $r<k$, the full hull
$D^{\circ r}(\widehat Q)$ contains $c_r$. By the Markov alternative it
meets $W_m$ exactly when $c_r\in W_m$; in that case it is contained in
$W_m$. Thus $t$ is the last physical visit of the critical orbit to
$W_m$ before time $k$. If $t=0$, then $c_t=0\in K_m$. If $t>0$,
synchronization gives
\[
 c_t\in h_m\Pc(D)\subset K_m\subset K_{m-1}.
\]
In either case $c_t\in K_{m-1}$. As $c_t\in D^{\circ t}(\widehat Q)=B(Q)$, this contradicts
parent-core separation. The counterexample is neither primitive nor
separated, a contradiction.

Lemma~\ref{lem:crossing} now supplies the single-valued section on all of
$\Omega_m$; its right-inverse identity makes it univalent.
\end{proof}

\begin{remark}[What Proposition~\ref{prop:main} does and does not say]
The result applies to the strong parent-core-separated class used in the
pairing below.  It neither asserts that arbitrary half-pullbacks form a Markov
family nor identifies this class with every separated-copy convention in the
literature.  It also does not require the two principal branches in
Lemma~\ref{lem:first-return} to exhaust the complete first-return map.
\end{remark}

\begin{corollary}[Extension over the full target]\label{cor:full-copy}
For an eligible positive-length tag, its extended inverse maps $W_m$
biholomorphically onto its full physical hull. Consequently that hull,
and every positive-length forward tail hull, does not contain $0$.
There are exactly two half-tags for each such full inverse component.
\end{corollary}
\begin{proof}
The physical hull is a component of the labelled preimage of $W_m$:
it is open there, and properness of its physical map makes it relatively
closed. The connected image of the extended inverse meets this component,
so remains in it. Its image is open and relatively closed in the hull;
closedness follows from the right-inverse identity and continuity on
$W_m$. Thus it is the entire hull. A hull containing $0$ would make
the positive-length physical iterate critical at $0$, contradicting
biholomorphy. Apply the same argument to eligible tails.
\end{proof}

% END sections/04_extension.tex

% BEGIN sections/05_pairing.tex
\section{Landing sets and disjoint loss copies}\label{sec:pairing}

Use $W_1$ as the base range.  Define
\begin{equation}\label{eq:landing-global}
 \widetilde X_r=
 \{z\in W_1:D^{\circ k}(z)\in W_r
      \text{ for some defined }k\ge0\},
 \qquad
 \widetilde\eta_r=
 \frac{\area(\widetilde X_r)}{\area(W_1)}.
\end{equation}
Then $\widetilde X_{r+1}\subset\widetilde X_r$. Define the normalized
half-plane landing sets and their areas by
\begin{equation}\label{eq:L-definition}
 \begin{gathered}
 L_r^\tau=h^{-1}(\widetilde X_r)\cap P_0^\tau
 =\{z\in P_0^\tau:D^k(hz)\in W_r\text{ for some defined }k\ge0\},\\
 a_r=\area(L_r^+)=\area(L_r^-).
 \end{gathered}
\end{equation}
Thus $L_{r+1}^\tau\subset L_r^\tau$ and
\[
 \area(\widetilde X_r)=2\alpha^4a_r,\qquad
 a_r=\widetilde\eta_r\area(P_0^+).
\]
The Jacobian $\alpha^4$ is independent of the critical degree. Numerical
upper bounds are denoted $A_r^+$ to distinguish them from $a_r$.

For $n\ge2$, $m\ge1$, set
\begin{equation}\label{eq:Xnm}
 X_{n,m}=
 \{z\in W_n:\cR_{n-1}^{\circ q}(z)\in W_{n+m}
       \text{ for some defined }q\ge0\}.
\end{equation}
Because $h_{n-1}^{-1}W_n=W_1$ and
$h_{n-1}^{-1}W_{n+m}=W_{m+1}$,
\begin{equation}\label{eq:stationary-X}
 X_{n,m}=h_{n-1}\widetilde X_{m+1},
 \qquad
 \frac{\area(X_{n,m})}{\area(W_n)}=\widetilde\eta_{m+1}.
\end{equation}
This is the precise source of the $m+1$ shift.

Define the nonreturn/loss sets, counting only defined physical iterates,
\begin{align}
 Y_n&=\{z\in W_n:D^{\circ k}(z)\notin W_n
                 \text{ for every defined }k\ge1\},\label{eq:Yn}\\
 \Sigma_n&=\{z\in W_n:\cR_{n-1}^{\circ q}(z)\in Y_n
                 \text{ for some defined }q\ge0\},\label{eq:Sigma-n}\\
 \Sigma_{n,m}&=\{z\in\Sigma_n:
       \cR_{n-1}^{\circ q}(z)\notin W_{n+m}
                 \text{ for every defined }q\ge0\}.
 \label{eq:Sigma-nm}
\end{align}
Then
\begin{equation}\label{eq:disjoint-monotone}
 X_{n,m}\cap\Sigma_{n,m}=\varnothing,
 \qquad \Sigma_{n,m}\nearrow\Sigma_n.
\end{equation}
Indeed, once a finite induced prefix reaches $Y_n$, its finite list avoids
$0$; for sufficiently large $m$ it avoids $W_{n+m}$, and after entry into
$Y_n$ no physical return to $W_n$ occurs.

Let $\cB_0$ be the union of all finite puzzle boundaries and all defined
finite preimages of $\R$.  This is a countable planar null set: finite puzzle boundaries are quasicircles, hence null, and finite holomorphic preimages of the real line are unions of analytic arcs and isolated critical points.  All tagged
decompositions below are understood modulo $\cB_0$.

\begin{lemma}[A point avoiding the next level]\label{lem:nonlanding}
For every $n\ge1$ there is $y\in K_n\subset W_n$ whose full physical
orbit avoids $W_{n+1}$.
\end{lemma}
\begin{proof}
The central map $C_n:\,h_nU^0\to W_n$ of $\mathscr R_n$ has a fixed
point. Here is an elementary existence argument. In a Riemann coordinate
on $W_n$, its source is compactly contained in the unit disk and the map
is proper of degree $\ell$ onto that disk. On a sufficiently high regular
level of its modulus, the modulus of the map exceeds that of the source
coordinate. Rouch\'e's theorem applied on the sublevel domain shows that
the map minus the identity has $\ell$ zeros counted with multiplicity.
Let $y$ be one. Then $y\in K_n$. It cannot belong to $W_{n+1}$, because
\[
 \mathscr R_n(W_{n+1})=h_nU^1,
\]
whose closure is disjoint from the central source containing $y$.
Every central return block starts and ends at $y$, and its proper
prefixes avoid $W_n$ by Lemma~\ref{lem:first-return}. Thus the physical
orbit avoids $W_{n+1}$.
\end{proof}

For an eligible tagged copy $Q\to P_n^{\sigma(Q)}$ with physical time $k_Q$
and full hull $\wh Q\subset W_1$, let $\psi_Q$ be its inverse and set
\begin{equation}\label{eq:paired-sets}
 \cX_Q=\psi_Q(X_{n,m}\cap P_n^{\sigma(Q)}),
 \qquad
 \cE_Q=\psi_Q(\Sigma_{n,m}\cap P_n^{\sigma(Q)}).
\end{equation}
Let $\cS_{n,m}$ contain the two identity tags $P_n^\pm$ and every nontrivial
eligible tag with $\wh Q\subset W_1$ and
\begin{equation}\label{eq:not-all-land}
 \wh Q\not\subset\widetilde X_{n+m}.
\end{equation}

\begin{lemma}[Exclusion of an induced transition]\label{lem:transition}
For every nonidentity $Q\in\mathscr S_{n,m}$,
\begin{equation}\label{eq:transition-exclusion}
 \widehat Q\cap(X_{n,m}\cup\Sigma_{n,m})=\varnothing.
\end{equation}
\end{lemma}
\begin{proof}
Write $C=\widehat Q$ and $k=k_Q>0$. If $C$ meets either set, it meets
$W_n$ and hence $C\subset W_n$ by depth and nesting. It is then
nonprimitive, so $C\cap K_{n-1}=\varnothing$ by
Lemma~\ref{lem:copy-sync}.

Suppose first that $x\in C\cap X_{n,m}$. Let $T$ be the full component
of an admissible parent-word pullback of $W_{n+m}$ through $x$.
It contains a preimage of $0$ in $K_{n-1}$, because $0\in K_{n-1}$
and finite admissible parent preimages preserve $K_{n-1}$.
Both $T$ and $C$ are full physical puzzle pieces. If $T\subset C$,
the core point contradicts $C\cap K_{n-1}=\varnothing$; if
$C\subset T$, then $C\subset\widetilde X_{n+m}$, contradicting
\eqref{eq:not-all-land}.

Next suppose $x\in C\cap\Sigma_{n,m}$. Let an admissible parent word
take $x$ to $Y_n$ and let $r$ be its physical length. Necessarily
$r\ge k$: otherwise the defined value $D^{\circ k}(x)\in W_n$
would be a positive-time return after the entry into $Y_n$.
Let $T$ be the full component of this parent-word pullback of $W_n$
through $x$. Its depth is at least that of $C$, so $T\subset C$.
Again $T$ contains an admissible parent preimage of $0$, which is a
point of $K_{n-1}$. This is the same contradiction.
\end{proof}

\begin{lemma}[Synchronized pairing]\label{lem:pairing}
Modulo $\mathscr B_0$,
\begin{align}
 \widetilde X_{n+m}&=\bigcup_{Q\in\mathscr S_{n,m}}\mathscr X_Q,
       \label{eq:representation}\\
 \{\mathscr X_Q,\mathscr E_Q:Q\in\mathscr S_{n,m}\}
       &\text{ is pairwise disjoint},\label{eq:multiplicity-one}\\
 \bigcup_{Q\in\mathscr S_{n,m}}\mathscr E_Q
       &\subset\widetilde X_n.\label{eq:escape-packing}
\end{align}
\end{lemma}
\begin{proof}
\emph{Coverage.} Fix $x\in\widetilde X_{n+m}\setminus\mathscr B_0$.
Let $j$ be the least physical time for which
$D^{\circ j}(x)\in X_{n,m}$. Such a time exists because the first
landing in $W_{n+m}$ has this property. If $j=0$, the appropriate
identity tag covers $x$.

Assume $j>0$ and let $C$ be the full pullback component of $W_n$
at time $j$ containing $x$. Depth and nesting give $C\subset W_1$.
Select the half-tag $Q$ through $x$. If it is nonprimitive, let
$r=\max\{i<j:D^i(C)\cap W_n\ne\varnothing\}$ and write $C_i=D^i(C)$. Nesting gives $C_r\subset W_n$, so $D^r(x)\in W_n$. If that hull met $K_{n-1}$,
Lemma~\ref{lem:copy-sync} would identify the suffix to time $j$
with an admissible parent word. Since its endpoint is in $X_{n,m}$,
the point at time $r$ would also be in $X_{n,m}$, contradicting
minimality of $j$. Hence $Q$ is eligible.

We must also check $C\not\subset\widetilde X_{n+m}$; this condition
cannot simply be assumed. By Corollary~\ref{cor:full-copy}, each
positive-length tail hull is univalent over $W_n$ and does not contain
$0$. If a tail $C_r$, $r<j$, met $W_{n+m}$, nesting would give
either $W_{n+m}\subset C_r$, impossible because $0\notin C_r$, or
$C_r\subset W_{n+m}$, impossible because it would put
$D^{\circ r}(x)$ in $X_{n,m}$ before time $j$. Thus every earlier
tail hull avoids $W_{n+m}$.
Pull the point $y$ of Lemma~\ref{lem:nonlanding} back by the inverse
on the full $W_n$. Its preimage in $C$ cannot land in $W_{n+m}$
before time $j$, by the preceding avoidance, or afterwards, by the
choice of $y$. This proves $Q\in\mathscr S_{n,m}$ and
$x\in\mathscr X_Q$. The reverse inclusion is immediate.

\emph{Disjointness.} Within one full inverse component the two
half-tags are disjoint, and within either tag the landing and loss sets
are disjoint. Distinct full hulls of the same time are disjoint.
If $k_1<k_2$ and their paired sets intersect, nesting forces
$C_2\subsetneq C_1$. Push $C_2$ forward by $k_1$. The resulting
positive-length full hull $C'$ lies in $W_n$, is eligible by tail
stability, and is not contained in $\widetilde X_{n+m}$: otherwise
its physical preimage $C_2$ would be contained there. It therefore
satisfies Lemma~\ref{lem:transition}. The larger copy's paired sets
map into $X_{n,m}\cup\Sigma_{n,m}$, whereas the smaller copy's paired
sets map into $C'$, a contradiction. This also handles comparison
with an identity tag.

\emph{Packing.} Every point of $\mathscr E_Q$ lies in $W_1$ and has
a defined physical image in $W_n$. Thus it belongs to
$\widetilde X_n$. Countability follows from the finite number of
words and sheets at each finite time.
\end{proof}

% END sections/05_pairing.tex

% BEGIN sections/06_recurrence.tex
\section{Weighted area recurrence and the zero-area criterion}\label{sec:criterion}
\subsection{Distortion and normalization}
For $z,w\in\Omega_*$, define
\begin{equation}\label{eq:Cstar}
 C_*(z,w)=\sup_\phi\frac{|\phi'(z)|}{|\phi'(w)|},
\end{equation}
where the supremum is over univalent $\phi:\Omega_*\to\C$.
Koebe distortion makes this finite at each pair and uniformly bounded
on compact subsets. The supremum is lower semicontinuous and hence Borel.
Let $\eps_n=(-1)^n$, identifying the signs with $\{+1,-1\}$, and normalize
\begin{align}
 A_{n,m}^\tau&=h_n^{-1}(\Sigma_{n,m}\cap P_n^{\eps_n\tau}),\label{eq:A-normalized}\\
 B_{n,m}^\tau&=h_n^{-1}(X_{n,m}\cap P_n^{\eps_n\tau}).\label{eq:B-normalized}
\end{align}
By \eqref{eq:stationary-X}, $B_{n,m}^\tau=L_{m+1}^\tau$ and
$\area(B_{n,m}^\tau)=a_{m+1}$. The first index in $A_{n,m}^\tau$
remains relevant: $Y_n$ is defined using ordinary physical iterates of
$D$, not a rescaled map. In particular $h^{-2}Y_3$ is not $Y_1$.
The parent conjugacy normalizes the domain and checkpoints but does not
remove the physical nonreturn condition.

% END sections/06_recurrence.tex

% BEGIN sections/06_support.tex
\subsection{Keeping the distortion integral on its actual support}
Let $\overline C(z,w)\ge C_*(z,w)$ be any finite Borel derivative-ratio
majorant on $\Omega_*\times\Omega_*$; Section~\ref{sec:kernel} gives the explicit
choice used in the computation. For a measurable set $E\subset P_0^\tau$
of positive area, put
\begin{equation}\label{eq:gE}
                  g_E(z)=\int_E\overline C(z,w)^{-2}\,dA(w).
\end{equation}
It is positive on $\Omega_*$. An upper majorant for distortion gives a lower
kernel and therefore a conservative area estimate.

\begin{lemma}[Integrated copy estimate]\label{lem:integrated}
If $\phi:\Omega_*\to\C$ is univalent and $B\subset P_0^\tau$ is
measurable, then
\begin{equation}\label{eq:single-weighted}
 \area(\phi(B))\le\area(\phi(E))\int_B\frac{dA(z)}{g_E(z)}.
\end{equation}
\end{lemma}
\begin{proof}
For each $z$,
$\area(\phi(E))\ge |\phi'(z)|^2g_E(z)$, by the change-of-variables
formula and the derivative-ratio bound. Divide by $g_E(z)>0$ and
integrate on $B$. An infinite integral makes the inequality vacuous.
\end{proof}

Fix $n\ge2$, $m_0\ge1$, and certified positive-area sets
$E^\tau\subset A_{n,m_0}^\tau$ from \eqref{eq:A-normalized}.
Loss monotonicity permits these same sets for every $m\ge m_0$.
The exact landing normalization is
\begin{equation}\label{eq:B-is-L}
 B_{n,m}^\tau=h_n^{-1}(X_{n,m}\cap P_n^{\eps_n\tau})
                    =L_{m+1}^\tau.
\end{equation}

\begin{proposition}[Weighted recursive estimate]\label{prop:weighted}
For every $m\ge m_0$,
\begin{equation}\label{eq:weighted-recursion}
 a_{n+m}\le a_n\max_{\tau\in\{+,-\}}
                 \int_{L_{m+1}^\tau}\frac{dA(z)}{g_{E^\tau}(z)}.
\end{equation}
\end{proposition}
\begin{proof}
For each $Q\in\mathscr S_{n,m}$ put $\tau_Q=\eps_n\sigma(Q)$,
$\phi_Q=\psi_Q\circ h_n$, and $F_Q=\phi_Q(E^{\tau_Q})$.
Proposition~\ref{prop:main} makes $\phi_Q$ univalent on $\Omega_*$;
this also holds for the identity tags. These are the paired landing and
loss copies in concrete terms:
\[
 \phi_Q(L_{m+1}^{\tau_Q})=\mathscr X_Q,
 \qquad F_Q\subset\mathscr E_Q.
\]
They use the \emph{same inverse branch and the same normalization}.
Apply \eqref{eq:single-weighted} to each copy. By
Lemma~\ref{lem:pairing}, the landing copies cover
$\widetilde X_{n+m}$ and the $F_Q$ are disjoint and packed in
$\widetilde X_n$. Therefore
\begin{align*}
 \area(\widetilde X_{n+m})
 &\le\sum_Q\area(\mathscr X_Q)\\
 &\le\left(\max_\tau\int_{L_{m+1}^\tau}g_{E^\tau}^{-1}\,dA\right)
                     \sum_Q\area(F_Q)\\
 &\le\left(\max_\tau\int_{L_{m+1}^\tau}g_{E^\tau}^{-1}\,dA\right)
                     \area(\widetilde X_n).
\end{align*}
Finally $\area(\widetilde X_r)=2\alpha^4a_r$, so the common factor
cancels. There is no additional half-plane or scaling factor.
\end{proof}

\begin{theorem}[Support-restricted uniform contraction]\label{thm:supported}
Fix integers $n\ge2$, $m_0\ge1$, $t\ge1$, and measurable sets
$E^\tau\subset A_{n,m_0}^\tau$ of positive area. Suppose there are
measurable sets $C_t^\tau\subset\Omega_*\cap H^\tau$ and constants
$A_n^+,\gamma>0$ such that, for both signs,
\[
 \area(L_t^\tau\setminus C_t^\tau)=0,\qquad
 \mathop{\rm ess\,inf}_{z\in C_t^\tau}g_{E^\tau}(z)\ge\gamma,
 \qquad a_n\le A_n^+<\gamma.
\]
Put $\theta=A_n^+/\gamma$ and $k_0=\max\{m_0+1,t\}$. Then
\begin{equation}\label{eq:supported-recursion}
 a_{k+n-1}\le\theta a_k\qquad(k\ge k_0),
\end{equation}
and
\[
 a_{k_0+r+q(n-1)}\le\theta^q a_{k_0+r}
 \qquad(0\le r<n-1,\ q\ge0).
\]
In particular $a_k\to0$ exponentially and $\area K(D)=\area J(D)=0$.
\end{theorem}
\begin{proof}
For $m\ge\max\{m_0,t-1\}$, the set $E^\tau$ remains an admissible
loss subset and $L_{m+1}^\tau\subset L_t^\tau\subset C_t^\tau$
modulo null sets. Hence the integral in
\eqref{eq:weighted-recursion} is at most $a_{m+1}/\gamma$.
Set $k=m+1$ to obtain \eqref{eq:supported-recursion}.
Iteration gives geometric decay on each residue class modulo $n-1$;
equivalently, monotonicity fills any intervening indices. The last
assertion follows from Lemma~\ref{lem:null-detection} below.
\end{proof}

No forward invariance of $C_t^\tau$ is required. Its role is only to
contain every sufficiently deep landing set. An integral smaller than
one for a \emph{single} pair $(n,m)$ would give only a finite-step
improvement; the lower bound on this fixed support is what makes
\eqref{eq:supported-recursion} uniform in all later indices.

\subsection{From landing-area decay to a null Julia set}
The following elementary lemma makes the noncritical-survivor argument
self-contained; see also the noncritical puzzle analysis in \cite{Smania}.
\begin{lemma}[Survival without a critical branch]\label{lem:noncritical-survivor}
Let $V_1,\ldots,V_N$ be bounded pairwise disjoint topological disks.
Let $U_1,\ldots,U_M$ be pairwise disjoint topological disks with
$U_j\Subset V_{s(j)}$. Suppose $g|_{U_j}:U_j\to V_{t(j)}$ is
biholomorphic. Then the set $K(g)$ of points with an infinite forward
orbit in $\bigcup_j U_j$ has area zero.
\end{lemma}
\begin{proof}
Choose Riemann maps $\chi_i:V_i\to\D$. Finiteness and compact
containment give a common $\kappa<1$ with
$\chi_{s(j)}(U_j)\subset\D_\kappa$. Each inverse branch is
$\kappa$-Lipschitz for the respective hyperbolic metrics: in disk
coordinates it has the form $\kappa u$, $u:\D\to\D$, and
Schwarz--Pick applies. Choose $\kappa<r<1$ and put
$V_i'=\chi_i^{-1}(\D_r)$. Every survivor in $V_i$ lies in $V_i'$.
Fix a disk $B_i$ compactly contained in
$\chi_i^{-1}(\D_r\setminus\overline{\D_\kappa})$, which avoids all
sources in $V_i$.

A compatible composition $b$ of $q$ inverse branches is univalent on its
terminal $V_i$. Koebe distortion on the fixed disk $\D_r$ gives a common
$C\ge1$, independent of $q,b$, such that
\[
 \area(b(B_i))\ge C^{-1}\area(b(V_i')),\qquad
 \operatorname{diam}(b(V_i'))^2\le C\area(b(V_i')).
\]
The fixed coordinate maps and the finitely many $B_i$ are absorbed in
$C$. Hyperbolic contraction also implies that these diameters tend
uniformly to zero. Every point of $b(B_i)$ escapes after its $q$
prescribed steps, so $b(B_i)\cap K(g)=\varnothing$.
For $x\in K(g)$, take the inverse $b_q$ along its itinerary. Then
$x\in b_q(V_i')$. The disk about $x$ of radius
$\operatorname{diam}(b_q(V_i'))$ contains $b_q(B_i)$, and the two
inequalities give a uniform positive proportion outside $K(g)$.
The radii tend to zero, so no survivor is a Lebesgue density point.
The density theorem proves the claim.
\end{proof}

\begin{lemma}\label{lem:null-detection}
If $a_r\to0$, then $\area K(D)=\area J(D)=0$.
\end{lemma}
\begin{proof}
The ordinary critical pieces $C_d$ shrink to $0$: they are nested and
$C_{S_{2j}-2}=W_j=h^jV$ has diameter tending to zero. At depth $d$,
restrict $D$ from the noncritical pieces of depth $d+1$ to the pieces of
depth $d$, calling the resulting map $g_d$. There are finitely many
branches. Each is biholomorphic, since its proper source does not
contain $0$, and each source is compactly contained in its depth-$d$
container. Lemma~\ref{lem:noncritical-survivor} gives $\area K(g_d)=0$.
If $x\in K(D)$ and $0\notin\omega_D(x)$, its orbit eventually avoids
some $C_d$, and hence has an iterate in $K(g_d)$. Thus
\[
 \{x\in K(D):0\notin\omega_D(x)\}
 \subset\bigcup_{d,N\ge0}D^{-N}(K(g_d)).
\]
Finite preimages of null sets under a nonconstant holomorphic map are
null, by local conformal inverses away from isolated critical points.
This proves nullity of the displayed set.

Set $S=\bigcap_{r\ge1}\widetilde X_r$. Continuity of measure from above,
boundedness of $W_1$, and $\area(\widetilde X_r)=2\alpha^4a_r$ give
$\area S=0$. Every orbit accumulating at $0$ has a forward image in
$W_1$ which subsequently enters every $W_r$; that image belongs to $S$.
The remaining part of $K(D)$ lies in $\bigcup_{N\ge0}D^{-N}(S)$, also
null. Finally $J(D)\subset K(D)$.
\end{proof}

% END sections/06_support.tex

\part*{\Large{Part II. Computer-assisted quartic application}}
\addcontentsline{toc}{part}{Part II. Quartic application}

% BEGIN sections/07_validation.tex
\section{Validation of the analytic solution}\label{sec:validation}
\subsection{The Banach-space equation}
For numerical work it is convenient to use a quotient function with
the normalization
\begin{equation}\label{eq:quotient-normalization}
 f(z)=\HH(z^4/\alpha^4),\qquad
 f_\alpha(z)=\HH(z^4),\qquad \HH(0)=1.
\end{equation}
The functional equation becomes
\begin{equation}\label{eq:banach-equation}
 \Phi(\alpha,\HH)(x)
 =\HH(x)+\alpha^{-2}\HH\bigl(\HH(\alpha^4x)^4\bigr)=0.
\end{equation}
Let $\mathcal A_R$ be the real Taylor algebra with norm
$\|\sum_{j\ge0}c_jx^j\|_R=\sum_{j\ge0}|c_j|R^j$,
where $R=13/5$. The unknowns are $\alpha$ and the coefficients of
$\HH-1$, with product norm $|\delta\alpha|+\|\delta\HH\|_R$.
The residual lies in $\mathcal A_R$; its constant coefficient is an
equation, not a free unknown.

The center is the exact terminating-decimal parameter
\[
 \overline\alpha=0.4724174456070521437409162653300465328574
\]
and the degree-20 polynomial $\overline H$ listed in
Appendix~\ref{app:coefficients}. All higher center coefficients are
zero, but the unknown function has an unrestricted, controlled tail.
Use the radius $r=3\cdot10^{-13}$.

\subsection{Finding the finite Newton candidate}
The search uses the vector $v=(\alpha,c_1,\ldots,c_{20})$, with $c_0=1$,
and the first 21 coefficients of \eqref{eq:banach-equation}. The constant
residual is an equation; it is not dropped. Writing
$p(x)=H(\alpha^4x)$ and $s(x)=p(x)^4$, the analytic Jacobian is
\begin{align*}
 \partial_\alpha\Phi&=-2\alpha^{-3}H(s)
       +16\alpha xH'(s)p^3H'(\alpha^4x),\\
 D_H\Phi[u]&=u(x)+\alpha^{-2}
       \bigl(u(s)+4H'(s)p^3u(\alpha^4x)\bigr),\qquad u(0)=0.
\end{align*}
Products and compositions are truncated at degree 20 for this candidate
search only. At each step solve $D\Phi_{20}(v)\delta=-\Phi_{20}(v)$
and replace $v$ by $v+\delta$; the code uses pivoted elimination.
The accompanying \texttt{newton\_candidate\_demo.py} starts reproducibly
from
\[
 \alpha_0=\sqrt{(\sqrt2-1)/2},\qquad c_1=-1-\alpha_0^2,
 \qquad c_2=\cdots=c_{20}=0.
\]
This is the nonconstant degree-one truncated solution. With 80 decimal
digits it converges to the degree-20 candidate; the full residual and step
trace are written to \texttt{runs/newton\_candidate\_demo.json}.
This fresh search demonstrates how to obtain the candidate. It is not a
historical search log, and its floating-point residual is not an existence
proof. The preserved terminating-decimal center is then validated by
Proposition~\ref{prop:validated}, using the infinite-dimensional operator
and its complete tail estimates.

\subsection{The full contraction certificate}
\begin{proposition}[Validated solution]\label{prop:validated}
There is a unique real solution of \eqref{eq:banach-equation} in
\begin{equation}\label{eq:banach-ball}
 |\alpha-\overline\alpha|+\|\HH-\overline H\|_{13/5}\le r.
\end{equation}
It satisfies $0<\alpha<1$ and $\HH'(0)<0$, with
$|\HH'(0)|>1.3274224713696$.
\end{proposition}
\begin{proof}
Let $B$ be the computed approximate inverse of the finite derivative,
with columns indexed by $\alpha,c_1,\ldots,c_{20}$ and rows indexed
by residual coefficients $0,\ldots,20$. The coefficient weights are
included in this matrix. Define the full operator
$A=B\oplus I_{>20}$ and the Newton map $\mathcal T=I-A\Phi$.
On the complete ball \eqref{eq:banach-ball}, the directed calculation
gives the following outward bounds:
\begin{center}
\begin{tabular}{@{}ll@{}}
\toprule
Quantity & Certified upper bound\\\midrule
$\|A\Phi(\overline\alpha,\overline H)\|$ & $Y=1.604384\cdot10^{-17}$\\
$\|I-A D\Phi(\overline\alpha,\overline H)\|$ & $Z=0.017348053$\\
$\sup\|A D^2\Phi\|$ on the ball & $L=7435.70$\\
$\|I-BD_{20}\|$ & $3.203\cdot10^{-97}$\\
$\|B\|$ & $1.887694$\\\bottomrule
\end{tabular}
\end{center}
In particular
\[
 Y+Zr+\tfrac12Lr^2<r,\qquad Z+Lr<0.017348056<1.
\]
The composition arguments remain strictly inside $|x|<R$ throughout
the ball; the majorants in Appendix~\ref{app:majorants} verify this
and the operator bounds, including the infinite tail.
Thus $\mathcal T$ maps the closed ball into itself and is a contraction.
The finite inverse-defect bound makes $BD_{20}$ invertible, hence
$B$ invertible; adjoining the identity tail makes $A$ injective.
The unique fixed point therefore satisfies $\Phi=0$, rather than
merely a projected equation. The real Banach space gives reality.
Finally $|\delta\HH'(0)|\le r/R$, which proves the derivative claim;
the parameter interval is contained in $(0,1)$.
\end{proof}

\subsection{Rigorous evaluation and connected continuation}
For $|x|\le s<R$, the coefficient-ball bound gives
$|\HH(x)-\overline H(x)|\le r$. Derivative errors follow, for example,
from
\begin{equation}\label{eq:derivative-tail}
 |(\HH-\overline H)^{(d)}(x)|
 \le\frac{d!\,r}{R^d(1-s/R)^{d+1}}\qquad(d\ge1).
\end{equation}
Polynomial evaluation, the interval for $\alpha$, and these tails give
uniform rectangular enclosures. The scalar $\alpha$ is never replaced
by its midpoint in a proof calculation.

The formula
\[
 x\longmapsto-\alpha^{-2}
         \HH\bigl(\HH(\alpha^4x)^4\bigr)
\]
provides a connected continuation whenever both nested arguments have
strictly valid enclosures. A coefficient majorant first proves this
on the entire disk $|x|<19/5$; the remaining outer-argument margin
in the local $13/5$ disk exceeds $0.119$. A second majorant extends
the connected chart to $|x|<11/2$: the inner evaluation is local and
the outer argument has modulus at most $3.538651<19/5$.
Each extension agrees with the preceding germ on an open overlap.

Other finite charts use explicit rational corridors with nonempty
open overlaps and a finite acyclic nesting of the same equation.
After the identification in Section~\ref{sec:identification}, Buff's maximal-extension assertion \cite[Definition 16 and Theorem B]{Buff}
identifies these connected continuations after lifting by $z\mapsto z^4$.
An isolated numerical evaluation of the functional expression, without
such connected agreement, is not used to recognize a physical branch.

\subsection{From local validation to canonical identification}
Proposition~\ref{prop:validated} proves existence and uniqueness inside a
specified Banach ball. It does not by itself identify the zero with the
canonical solution, nor does a derivative sign on a disk establish the
anti-Herglotz inverse property in P1. Section~\ref{sec:identification}
constructs proper domains independently, proves the complete Fibonacci
combinatorics by induction, and only then applies cycle uniqueness.

% END sections/07_validation.tex

% BEGIN sections/07_identification.tex
\section{Canonical identification through an independently constructed Fibonacci map}\label{sec:identification}
We now prove that the solution of Proposition~\ref{prop:validated} induces
a Fibonacci map. In this section only, $\alpha,\HH,f,q,G,D,B,V,U^j$
denote objects constructed from that validated pair. Their identification
with the canonical objects of Section~\ref{sec:preliminaries} is the
conclusion, not a premise. The local input is exactly the full Taylor
ball \eqref{eq:banach-ball} and equation \eqref{eq:banach-equation},
including the unknown tail. The independent finite checks use only these data.
Put
\begin{equation}\label{id:eq:maps}
 q(z)=\HH(z^4),\quad G(z)=q(\alpha z),\quad
 f(z)=q(z/\alpha),\quad h=-\alpha^2.
\end{equation}
Initially these are coherent germs and their validated analytic charts.
The same residual gives
\begin{equation}\label{id:eq:local-identities}
 \begin{gathered}
 q\circ G=hq,\qquad q(0)=1,\qquad q(1)=h,\\
 f(z)=-\alpha^{-2}f\bigl(\alpha f(\alpha z)\bigr).
 \end{gathered}
\end{equation}
Write $\disk s=\{z\in\C:|z|<s\}$.

\begin{theorem}[Canonical identification]\label{thm:canonical-identification}
The validated pair in \eqref{eq:banach-ball} is the canonical quartic pair:
\[
 \alpha=\alpha_{\rm can},\qquad
 f|_{\disk\epsilon}=f_{\rm can}|_{\disk\epsilon}
 \quad\text{for some }\epsilon>0.
\]
Equivalently, if $f_{\rm can}(z)=F_{\rm can}(z^4)$, then
\[
 F_{\rm val}(x):=\HH(x/\alpha^4)=F_{\rm can}(x),\qquad
 \HH(x)=F_{\rm can}(\alpha^4x)
\]
as germs at zero. In particular, $\HH$ and $F_{\rm can}$ have different
input normalizations. The germ induces a real-symmetric generalized polynomial-like map
$D:U^0\sqcup U^1\to V$, of branch degrees four and one, with
Fibonacci combinatorics and $\mathcal R^{\circ2}[D]=[D]$.
The two renormalization phases are distinct.
\end{theorem}

The uniqueness theorem used at the end is for real Fibonacci
renormalization cycles, as stated in \cite[Section 4.2, Theorem 3]{Buff}. We do not
apply uniqueness to an arbitrary zero of \eqref{eq:banach-equation}. The proof
first puts this particular zero in the dynamical class required by that
theorem. The canonical slit-domain geometry is not used to construct the domains
in this section. Buff explicitly attributes his cycle theorem to
van Strien--Nowicki \cite[Theorem 7.1, as cited in Buff's Theorem 3]{NS};
we use that result with the normalization and phase described in P4.

\subsection{Independent finite analytic inequalities}\label{id:sec:finite}
The programs \texttt{rigorous\_contours.py} and \texttt{basin\_checks.py}
prove the inequalities in this subsection with exact rational interval
arithmetic. The files \texttt{G\_disk\_16\_5\_certificate.json} and
\texttt{basin\_checks\_certificate.json} record the 1,024 contour boxes
and eight basin checks, respectively. The displayed decimals are weaker rational bounds selected
for readability; the certificates contain exact fractions.

\begin{proposition}[The new finite certificate]\label{id:prop:finite}
The coherent continuation of $\HH$ is holomorphic on $\disk{11/2}$ and
has no derivative zero there. Set $R=16/5$ and $r=8/5$. Then $G$ is
holomorphic on a neighborhood of $\overline{\disk R}$. Indeed,
$\alpha<0.4725$ and $(\alpha R)^4<1.512^4=5.226454388736<11/2$
place the entire closed disk in the continuation chart. We have
\begin{align}
 |G(z)|^2-R^2&>1 &&(|z|=R),\label{id:eq:boundary}\\
 |G(z)|&<1.445<r &&(|z|\le r).\label{id:eq:trap}
\end{align}
Moreover,
\begin{gather}
 \alpha R<1.512<r,\qquad \alpha R>1,\qquad
 \alpha^2R<0.715<1,\label{id:eq:scaling}\\
 G([0.9,1])\subset(0.934,0.957),\qquad
 \operatorname{Re}\HH'(x)<-0.7359\quad(|x|\le2.59).
 \label{id:eq:realbounds}
\end{gather}
\end{proposition}
\begin{proof}
Write $p=\overline H$, $\varepsilon=3\cdot10^{-13}$ and
$M(s)=\sum_{j=0}^{20}|p_j|s^j$. The local estimate gives
$|\HH(x)-p(x)|\le\varepsilon$ for $|x|<13/5$. On
$|x|\le d=259/100$,
\[
 |\HH'(x)-p_1|
 \le\sum_{j=2}^{20}j|p_j|d^{j-1}
       +\frac{\varepsilon}{(13/5)(1-d/(13/5))^2}.
\]
Subtracting the right side from $|p_1|$ gives a bound greater than
$0.73595$. Since $p_1<0$, this also bounds $-\operatorname{Re}\HH'$.

For an extension radius $s$, set $u=\alpha_+^4s$, where $\alpha_+$
is the upper parameter endpoint. The expression on the right of
\eqref{eq:banach-equation} is defined throughout $\disk s$ whenever
\[
 u<13/5,\qquad (M(u)+\varepsilon)^4<s_{\rm previous}.
\]
The exact scalar calculations give these successive extensions:
\begin{center}\small
\begin{tabular}{@{}llll@{}}\toprule
New radius & Previous radius & Outer argument bound & Inner $|\HH|$ bound\\\midrule
$19/5$ & $13/5$ & $<2.481<2.59$ & $>0.744$\\
$11/2$ & $19/5$ & $<3.539<3.8$ & $>0.628$\\\bottomrule
\end{tabular}
\end{center}
Each expression agrees with the preceding germ near zero by
\eqref{eq:banach-equation}, hence throughout its overlapping disk. Differentiation gives
\[
 \HH'(x)=-4\alpha^2\HH'\bigl(\HH(\alpha^4x)^4\bigr)
              \HH(\alpha^4x)^3\HH'(\alpha^4x).
\]
All factors are nonzero by the local derivative bound, the inner lower
bounds and induction through the extensions. This also covers the
annulus between radii $2.59$ and $2.6$.

Since $\alpha_+^4R^4<5.223<11/2$, $G$ is analytic on a neighborhood
of the closed disk, and
\begin{equation}\label{id:eq:Gderivative}
 G'(z)=4\alpha^4z^3\HH'(\alpha^4z^4)
\end{equation}
vanishes precisely at zero, to order three.

For the full boundary, use
\[
 z(t)=R\left(\frac{1-t^2}{1+t^2}
                   +i\frac{2t}{1+t^2}\right),\qquad0\le t\le1.
\]
Partition $[0,1]$ into 1,024 equal intervals. Monotonicity of each
coordinate gives exact rectangles containing the entire corresponding
arcs. Quartic symmetry $G(iz)=G(z)$ covers the other quadrants.
For each box $Z_j$, outward Horner evaluation, the full Taylor error,
and the connected functional expressions give
$\mathbf G_j\supset G(Z_j)$. The exact minimum over all rows of
\[
 \inf\{a^2:a\in\operatorname{Re}\mathbf G_j\}
 +\inf\{b^2:b\in\operatorname{Im}\mathbf G_j\}-R^2
\]
is greater than $1.0295784546957$. This proves \eqref{id:eq:boundary}.

For the whole trapping disk the scalar bound is
\[
 |G(z)|\le M(\alpha_+^4r^4)+\varepsilon<1.444525\qquad(|z|\le r).
\]
The parameter interval gives \eqref{id:eq:scaling}. A whole real-interval
evaluation gives an enclosure inside $(0.934054,0.956792)$ for
$G([0.9,1])$. These last checks are in the separate basin-check block.
\end{proof}

\begin{figure}[htbp]
\centering\includegraphics[width=1\textwidth]{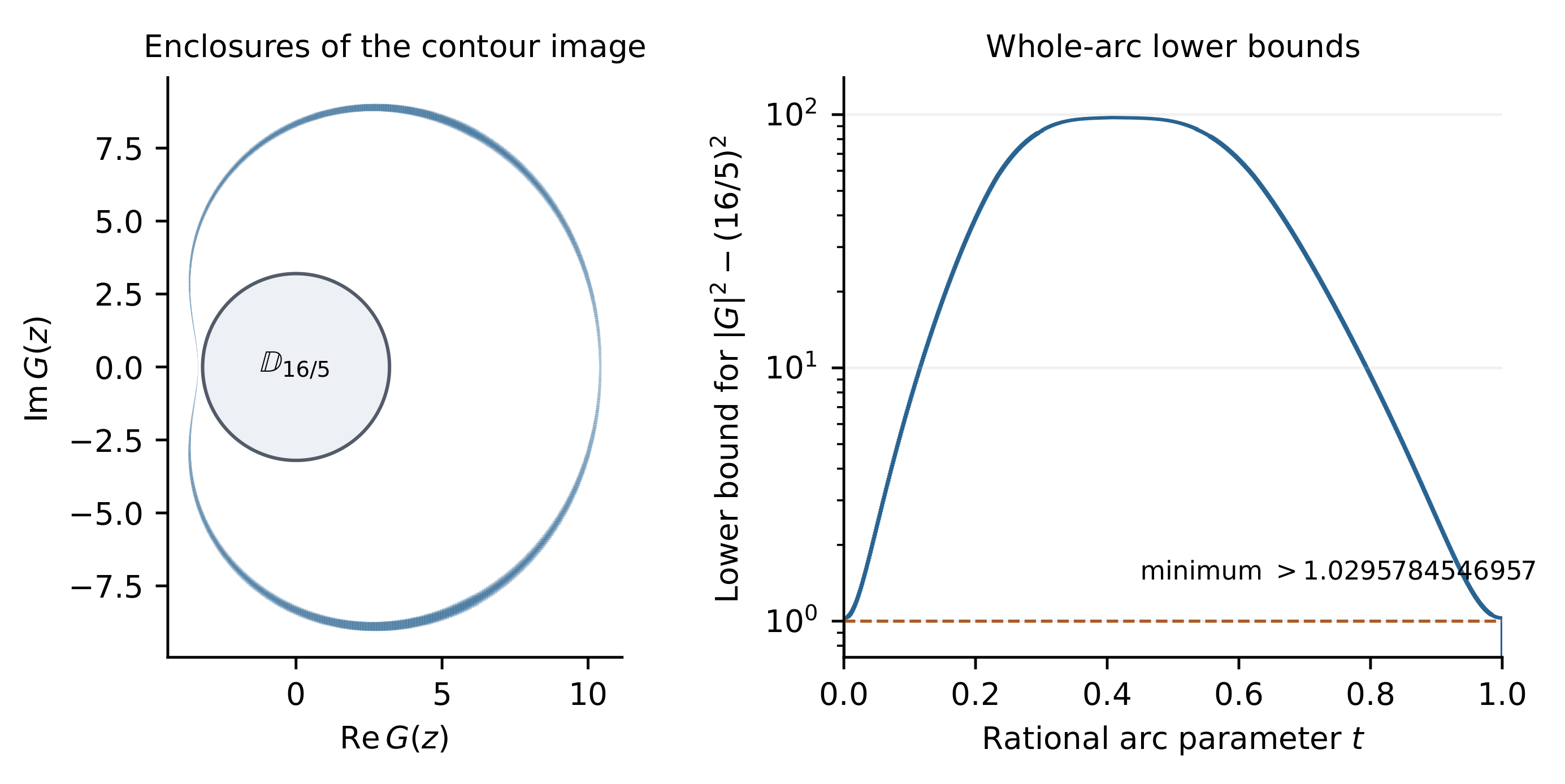}
\caption{The complete contour certificate for $G$. Each rectangle encloses the image of one full parameter interval; the 1,024 intervals cover the quarter-circle, and quartic symmetry covers the rest. The margins are lower bounds for $|G|^2-(16/5)^2$. Plot coordinates are decimal displays; the certificate comparisons use exact fractions.}\label{fig:contours}
\end{figure}
\subsection{The attracting basin and its linearizing function}
\begin{lemma}[A polynomial-like restriction]\label{id:lem:poly}
Let $\mathcal U$ be the component containing zero of
$\{z\in\disk R:|G(z)|<R\}$. Then $G:\mathcal U\to\disk R$
is polynomial-like of degree four, with only critical point zero.
Its immediate attracting basin $B$ has a fixed point $\rho\in(0.9,1)$,
and
\begin{equation}\label{id:eq:Bbounds}
 \disk r\subset B,\quad \overline B\Subset\mathcal U\Subset\disk R,
 \quad \alpha\overline B\Subset B,\quad [0,1/\alpha]\subset B.
\end{equation}
The domain $B$ is a quasidisk and $G:B\to B$ is proper of degree four.
\end{lemma}
\begin{proof}
The boundary inequality makes each sublevel component relatively
compact in $\disk R$ and its restriction proper onto $\disk R$.
For any Jordan curve in such a component, the maximum modulus principle
puts its entire interior in $|G|<R$ and hence in the same component.
The curve and its interior are in $\disk R$, where $G$ is holomorphic.
The component is therefore simply connected. Equation
\eqref{id:eq:Gderivative} and Riemann--Hurwitz now give degree four for
the component through zero. Neither a winding-number certificate nor
global injectivity of $\HH$ is needed.

The disk $\disk r$ lies in $\mathcal U$ and is mapped compactly into
itself. Denjoy--Wolff gives an interior attracting fixed point and
convergence throughout the disk. The real-interval inclusion locates
this point in $(0.9,1)$. Thus the only critical point is attracted to
it. Straightening \cite{DH} and the attracting-fixed-point case for unicritical
polynomials imply that $B$ is a quasidisk, $\overline B=K(G)$, and
$G:B\to B$ is proper of degree four. These are the standard
polynomial-like facts used in \cite[Section 2.1, Proposition 1; Section 4.3, paragraph following
Lemma 6]{Buff}. Here their required proper map and attracting critical
orbit have been established independently; the inverse-class theorem
is not a premise of this application.

The derivative bound makes $\HH$ strictly decreasing on $[0,1]$.
The residual at zero gives $\HH(1)=-\alpha^2$, so
\[
 G([0,1/\alpha])=\HH([0,1])=[-\alpha^2,1]\subset\disk r.
\]
This interval is in $\disk R$, by \eqref{id:eq:scaling}. It is a connected
set of basin points attached to zero, so belongs to $B$.
Finally $\alpha\overline B\subset\disk r$ follows from $\alpha R<r$.
Quartic symmetry of $G$ and $0\in B$ give $B=-B$ and $iB=B$.
\end{proof}

The original local identity $q\circ G=hq$ holds on neighborhoods of
the critical orbit from zero into the real invariant interval in
\eqref{id:eq:realbounds}. At the fixed point it gives $q(\rho)=0$.
The local derivative bound gives $q'(\rho)\ne0$, and differentiation
then yields $G'(\rho)=h=-\alpha^2$.

Choose a small forward-invariant neighborhood $N$ of $\rho$ where this
identity holds. For $z\in B$ and sufficiently large $k$ set
\begin{equation}\label{id:eq:Qextension}
 Q(z)=h^{-k}q(G^{\circ k}(z)),\qquad G^{\circ k}(z)\in N.
\end{equation}
The identity makes this independent of $k$ and consistent on overlaps,
so $Q$ is holomorphic on $B$ and
\begin{equation}\label{id:eq:Qrelation}
                         Q\circ G=hQ\quad\text{on }B.
\end{equation}
It agrees with the original germ at zero: the finite segment of the
critical orbit into $N$ uses the original local charts, and the same
identities hold on neighborhoods of all its points. This gives agreement
on an open neighborhood, not just on a real sequence. Reality and quartic
rotational symmetry also extend to $Q$.

Put $V=\alpha B$. The functions $Q(z)$ and $G(z/\alpha)$ are holomorphic
on this connected domain and agree near zero. Thus, writing $q=Q$ on
$B$ and $f(z)=Q(z/\alpha)$ on $V$, we obtain
\begin{equation}\label{id:eq:qproper}
 q:V\longrightarrow B\text{ proper of degree four},\qquad
 q(z)=G(z/\alpha)\ (z\in V),\qquad V\Subset B.
\end{equation}
Its only critical point is zero, and $1\in V$ by \eqref{id:eq:Bbounds}.
This proper map has been constructed without the canonical slit inverse.

\subsection{The two physical branch domains}
Set $Z=\alpha^2B=\alpha V$ and
\begin{equation}\label{id:eq:Uone}
 U^1=\Comp_1\bigl((q|_V)^{-1}(Z)\bigr).
\end{equation}
The mark is valid because $q(1)=-\alpha^2\in Z$. The domain
$Z\Subset B$ avoids the only critical value $1$ of $q|_V$, by
$\alpha^2R<1$. Properness and simple connectivity give
\begin{equation}\label{id:eq:Uoneproper}
                 U^1\Subset V,\qquad q:U^1\longrightarrow Z
                 \text{ conformal and onto}.
\end{equation}
It is also the full component through $1$ of $Q^{-1}(Z)$ in $B$:
its closure lies in $V\Subset B$, and its boundary maps into
$\partial Z$, so that component cannot continue across the boundary.

Now $Q(G(V))=hQ(V)=hB=Z$. The connected set $G(V)$ contains $1=G(0)$,
so lies in the full component $U^1$. Conformality of $q|_{U^1}$ and
surjectivity of the last image give
\begin{equation}\label{id:eq:two-components}
                       G(V)=U^1,\qquad q(Z)=U^1,\qquad q(U^1)=Z.
\end{equation}
Here $q(\alpha V)=G(V)$ proves the middle equality; that restriction
has degree four.

The exact equation and evenness give
\begin{equation}\label{id:eq:q2}
                         q^{\circ2}(hz)=hG(z)\qquad(z\in B).
\end{equation}
All intermediate points belong to $Z\cup U^1\subset V$, by
\eqref{id:eq:two-components}. Thus $q^{\circ2}|_Z$ is conjugate to $G|_B$.
Its iterates converge to $p_0=h\rho<0$ throughout $Z$; throughout
$U^1=q(Z)$ they converge to $p_1=q(p_0)>0$. The last sign follows
from the independently verified bound
\[
 |q(h\rho)-1|\le M(\alpha_+^8)-1+\varepsilon<0.003294,
\]
since $0<\rho<1$. Distinct limits imply $Z\cap U^1=\varnothing$.
Both domains consist of nonescaping points of $q:V\to B$.
Their closures may meet; disjointness of their boundaries is not needed.

Define auxiliary domains using only the displayed proper maps:
\begin{equation}\label{id:eq:auxiliary}
 \begin{gathered}
 \Ak_{-1}=\Jk_{-1}=V,\\
 \Ak_k=(G|_B)^{-(k+1)}V,\qquad
 \Jk_k=(q|_V)^{-(k+1)}V\quad(k\ge0).
 \end{gathered}
\end{equation}
These full inverse images are connected quasidisks. Every successive
range contains the critical value $1$, whose critical preimage accounts
for the full degree. Properness yields
\begin{equation}\label{id:eq:auxiliary-inclusions}
 \begin{gathered}
 \Ak_k\Subset\Ak_{k+1}\Subset B\quad(k\ge-1),\\
 \Jk_{k+1}\Subset\Jk_k\Subset V\quad(k\ge0),\qquad
 K(q)\Subset\Jk_k,\qquad \Jk_0=\alpha\Ak_0.
 \end{gathered}
\end{equation}
The first nesting uses $G(\overline V)=\overline{U^1}\Subset V$;
the second uses $V\Subset B$ and successive proper preimages under $q$.
Scaling $q(\alpha z)=G(z)$ proves the equality. In particular,
\begin{equation}\label{id:eq:useful-inclusion}
 \begin{gathered}
 U^1\subset K(q)\subset\Jk_0=\alpha\Ak_0\subset\alpha\Ak_k,\\
 \Jk_{2k}\subset\alpha\Ak_k\qquad(k\ge0).
 \end{gathered}
\end{equation}

Finally set
\begin{equation}\label{id:eq:D}
 U^0=\alpha^2\Ak_0,\qquad
 D:U^0\sqcup U^1\to V,\qquad
 D(z)=\begin{cases}f(z),&z\in U^0,\\
                  \alpha^{-1}q(z),&z\in U^1.
       \end{cases}
\end{equation}
Since $f(\alpha^2z)=G(z)$ on $B$, the central branch is proper of
degree four onto $V$. The lateral branch is conformal onto $V$.
Also $U^0\Subset Z$ and $U^1\Subset V$, with
$\overline{U^0}\cap\overline{U^1}=\varnothing$: the latter follows
from $\overline{U^0}\subset Z$ and $U^1\cap Z=\varnothing$.
All domains and maps are real symmetric. Nonescape of the critical
orbit of $D$ is established by the next induction.

\begin{figure}[htbp]
\centering
\begin{tikzpicture}[manuscript diagram]
  \node (Z) at (-3,1.15) {$Z=\alpha^2B$};
  \node (U1) at (3,1.15) {$U^1$};
  \node (U0) at (-3,-1.15) {$U^0=\alpha^2\mathcal A_0\Subset Z$};
  \node (V) at (3,-1.15) {$V=\alpha B$};
  \draw[map arrow] (Z.north east) to[bend left=19]
    node[map label,above=3pt] {$q$\quad (degree $4$)} (U1.north west);
  \draw[map arrow] (U1.south west) to[bend left=19]
    node[map label,below=3pt] {$q$\quad (degree $1$)} (Z.south east);
  \draw[map arrow] (U1.south) --
    node[map label,right=5pt] {$D_1=\alpha^{-1}q$\\(degree $1$)} (V.north);
  \draw[map arrow] (U0.east) --
    node[map label,below=4pt] {$D_0=f$\quad (degree $4$)} (V.west);
\end{tikzpicture}
\caption{The independently constructed proper domains and their mapping degrees. The upper domains lie in distinct attracting basins for $q^2$. The source $U^0\Subset Z$ and $U^1$ have disjoint closures in $V$. This is a diagram of proper maps, not of Euclidean basin shapes.}\label{fig:independent-domains}
\end{figure}
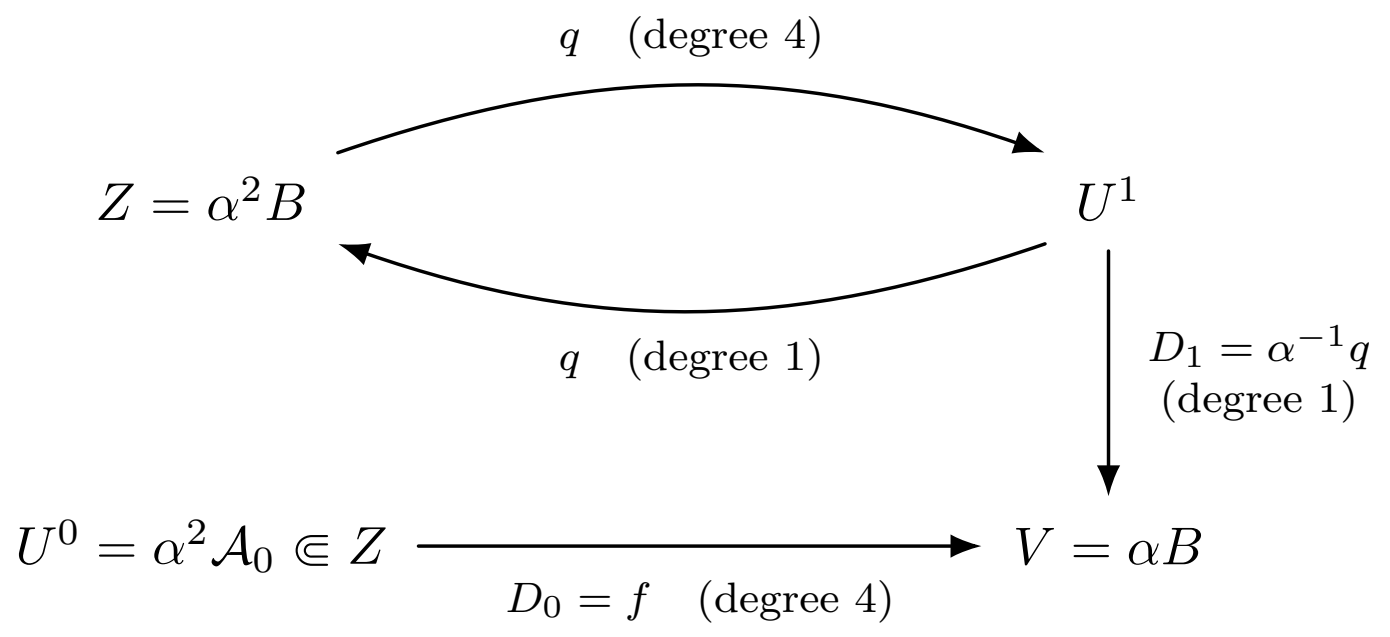
\subsection{The infinite Fibonacci induction on actual domains}
\label{id:sec:cycle-induction}

We now use the independently constructed domains to prove that $D$ is a
Fibonacci map and that its renormalization germ has period two. No
Epstein-class assumption enters this step. Throughout this section,
$h=-\alpha^2$ denotes a scalar, so $hE=\{hz:z\in E\}$.
The domains $B,V,\mathcal A_k,\mathcal J_k$ are invariant under
$z\mapsto-z$. We use the inclusions
\eqref{id:eq:auxiliary-inclusions}--\eqref{id:eq:useful-inclusion}; in particular,
\begin{equation}\label{id:eq:cycle-nesting-used}
 \begin{split}
 &V\subset\mathcal A_k\subset B,\qquad
   \mathcal J_{k+1}\subset\mathcal J_k\subset V,\qquad
   \mathcal J_0=\alpha\mathcal A_0,\\
 &U^1\subset K(q)\subset\mathcal J_0\subset\alpha\mathcal A_k,
   \qquad
   \mathcal J_{2k}\subset\mathcal J_0\subset\alpha\mathcal A_k
   \quad(k\ge0).
 \end{split}
\end{equation}
The notation $K(q)$ refers to the filled set of the polynomial-like
restriction $q:V\to B$. An iterate of $D$ below always uses the two
physical domains in \eqref{id:eq:D}.

\begin{lemma}[The initial return diagram]\label{id:lem:cycle-base-return}
The restrictions
\[
 D^{\circ3}:hU^0\longrightarrow hV,
 \qquad D^{\circ2}:hU^1\longrightarrow hV
\]
are proper first-return branches of degrees $4$ and $1$, respectively.
They satisfy
\begin{equation}\label{id:eq:cycle-base-identities}
 D^{\circ3}(hz)=h f(z)\quad(z\in U^0),\qquad
 D^{\circ2}(hz)=h\alpha^{-1}q(z)\quad(z\in U^1).
\end{equation}
\end{lemma}
\begin{proof}
Since $\alpha^2V\subset U^0$ and
$f(\alpha^2V)=q(\alpha V)=U^1$, the two successive restrictions
\[
 D:\alpha^2V\longrightarrow U^1,
 \qquad D:U^1\longrightarrow V
\]
are proper of degrees $4$ and $1$. Thus their composition is a
physical iterate, and the functional equation gives
\begin{equation}\label{id:eq:cycle-base-two}
 D^{\circ2}(z)=-\alpha f(z/\alpha)
                      \qquad(z\in\alpha^2V).
\end{equation}
Indeed the equality holds near $0$ by the functional identity and
then throughout this connected domain by the identity theorem.

Both $hU^0$ and $hU^1$ lie in $hV=\alpha^2V$. Their first images
lie in $U^1$, which is disjoint from $hV\subset U^0$. For the
central branch, \eqref{id:eq:cycle-base-two} gives
\[
 D^{\circ2}(hU^0)=-\alpha q(U^0)\subset-\alpha U^1
                                      \subset U^0.
\]
The first inclusion uses $U^0\subset Z=\alpha^2B$ and
$q(Z)=U^1$. The second follows from
$U^1\subset\mathcal J_0=\alpha\mathcal A_0$ and
$U^0=\alpha^2\mathcal A_0$.
Also $U^1\cap\alpha V=\varnothing$, because these are the two
open basin domains in \eqref{id:eq:two-components}.
Consequently $(-\alpha U^1)\cap hV=\varnothing$.
The third central iterate is therefore physically defined, and
\[
 D^{\circ3}(hz)=f(-\alpha q(z))
               =f(\alpha q(z))=-\alpha^2f(z).
\]
For the lateral branch, \eqref{id:eq:cycle-base-two} already gives
$D^{\circ2}(hz)=-\alpha q(z)=h\alpha^{-1}q(z)$.
The two identities conjugate the original proper branches of $D$
to the displayed return restrictions. The intermediate-domain
exclusions prove the first-return assertion.
\end{proof}

Set
\[
 S_0=1,\qquad S_1=2,\qquad S_{n+1}=S_n+S_{n-1},
 \qquad e_n=(-1)^{\lfloor(n+1)/2\rfloor},
\]
and, wherever the formula is defined, put
\[
                         \Theta_n(z)=e_n\alpha^n f(z/\alpha^n).
\]

\begin{lemma}[All Fibonacci iterates]\label{id:lem:all-fibonacci-iterates}
For every $k\ge0$, the following are actual physical iterates,
proper of degree $4$, with sole critical point $0$:
\begin{align}
 D^{\circ S_{2k}}=\Theta_{2k}:\quad
 \alpha^{2k+2}\mathcal A_k
   &\longrightarrow\alpha^{2k}\mathcal A_{k-1},
                         \label{id:eq:even-physical-iterate}\\
 D^{\circ S_{2k+1}}=\Theta_{2k+1}:\quad
 \alpha^{2k+2}V
   &\longrightarrow\alpha^{2k}V.
                         \label{id:eq:odd-physical-iterate}
\end{align}
In particular, the critical orbit of $D$ is nonescaping and
\begin{equation}\label{id:eq:all-critical-values}
                         D^{\circ S_n}(0)=e_n\alpha^n
                         \qquad(n\ge0).
\end{equation}
\end{lemma}
\begin{proof}
For $k=0$, \eqref{id:eq:even-physical-iterate} is the definition
$f:U^0=\alpha^2\mathcal A_0\to V$, and
\eqref{id:eq:odd-physical-iterate} is
\eqref{id:eq:cycle-base-two}. Suppose the two assertions are known
through index $k-1$, where $k\ge1$.

\emph{The even step.} Let $X=\alpha^{2k+2}\mathcal A_k$.
Since $\mathcal A_k\subset B=\alpha^{-1}V$ and
$\alpha V\subset V$, we have
\[
 X\subset\alpha^{2k+1}V\subset\alpha^{2k}V.
\]
The previous odd iterate is therefore defined on $X$. Its image
satisfies
\begin{align*}
 \Theta_{2k-1}(X)
 &=(-1)^k\alpha^{2k-1}q(\alpha^2\mathcal A_k)\\
 &\subset(-1)^k\alpha^{2k-1}U^1
 \subset\alpha^{2k}\mathcal A_{k-1}.
\end{align*}
The last inclusion is \eqref{id:eq:cycle-nesting-used} and symmetry.
Thus the previous even iterate is defined on the entire image.
The functional equation and evenness give, near $0$,
\[
                 \Theta_{2k-2}\circ\Theta_{2k-1}=\Theta_{2k}.
\]
Both sides are holomorphic on the connected domain $X$, whose
physical itinerary has just been checked. The identity theorem
therefore gives the equality on $X$. Finally,
\[
 \Theta_{2k}(\alpha^{2k+2}z)=(-1)^k\alpha^{2k}G(z)
                                  \qquad(z\in\mathcal A_k).
\]
Since $G:\mathcal A_k\to\mathcal A_{k-1}$ is proper of degree
$4$ with sole critical point $0$, this proves
\eqref{id:eq:even-physical-iterate}, including properness and degree.

\emph{The odd step.} Let $Y=\alpha^{2k+2}V$. It lies in the
just-established even domain. Moreover,
\[
 \Theta_{2k}(Y)=(-1)^k\alpha^{2k}f(\alpha^2V)
              =(-1)^k\alpha^{2k}U^1
              \subset\alpha^{2k}V.
\]
The previous odd iterate is defined on this image. The functional
equation again gives a germ identity
\[
                 \Theta_{2k-1}\circ\Theta_{2k}=\Theta_{2k+1},
\]
which extends to $Y$ by the same argument. For $z\in V$,
\[
 \Theta_{2k+1}(\alpha^{2k+2}z)
                    =(-1)^{k+1}\alpha^{2k+1}q(z).
\]
The proper map $q:V\to B=\alpha^{-1}V$ proves
\eqref{id:eq:odd-physical-iterate} with the asserted degree and
critical point.

This completes the induction. Since $S_n\to\infty$ and each
iterate $D^{\circ S_n}$ is physically defined on a neighborhood
of $0$, every finite iterate of $0$ is defined. Thus the critical
orbit never escapes the physical domains. Evaluating at $0$ gives
\eqref{id:eq:all-critical-values}.
\end{proof}

We next identify the actual puzzle depths. For $d\ge0$, let $C_d$
be the component of $\operatorname{Dom}(D^{\circ(d+1)})$ containing
$0$, and set $C_{-1}=V$. If a point $x\in K(D)$ is specified,
$P_d(x)$ denotes its puzzle component at depth $d$.

\begin{lemma}[Critical puzzle domains]\label{id:lem:cycle-puzzle-domains}
For every $k\ge0$ and $n\ge0$,
\begin{align}
 C_{S_{2k}-2}&=\alpha^{2k}V,
 & C_{S_{2k+1}-2}&=\alpha^{2k+2}\mathcal A_k,
                                     \label{id:eq:cycle-puzzle-parities}\\
 C_{S_{n+1}-3}&=\alpha^n\mathcal J_{n-1}.
                                     \label{id:eq:cycle-puzzle-adjacent}
\end{align}
\end{lemma}
\begin{proof}
A connected physical pullback that maps properly onto its prescribed
range is a full component of the inverse image: it is open there,
and properness makes it relatively closed. This observation applies
to the domains of Lemma~\ref{id:lem:all-fibonacci-iterates}.

For \eqref{id:eq:cycle-puzzle-parities}, the initial domains are
$C_{-1}=V$ and $C_0=\alpha^2\mathcal A_0$.
The odd iterate in \eqref{id:eq:odd-physical-iterate} pulls
$\alpha^{2k}V$ back to $\alpha^{2k+2}V$ through $0$, increasing
its depth by $S_{2k+1}$. The even iterate with index $k+1$ pulls
$\alpha^{2k+2}\mathcal A_k$ back to
$\alpha^{2k+4}\mathcal A_{k+1}$, increasing its depth by
$S_{2k+2}$. The Fibonacci identities give precisely the indices
in \eqref{id:eq:cycle-puzzle-parities}.

For \eqref{id:eq:cycle-puzzle-adjacent}, the base case is
$C_{S_1-3}=C_{-1}=\mathcal J_{-1}=V$.
For $z\in\mathcal J_n$ the formula is
\[
 \Theta_n(\alpha^{n+1}z)=e_n\alpha^n q(z).
\]
This is a physically defined restriction. If $n=2k+1$, its domain
$\alpha^{2k+2}\mathcal J_{2k+1}$ is contained in
$\alpha^{2k+2}V$, the domain of
\eqref{id:eq:odd-physical-iterate}. If $n=2k$, the inclusion
$\mathcal J_{2k}\subset\alpha\mathcal A_k$ in
\eqref{id:eq:cycle-nesting-used} places its domain in that of
\eqref{id:eq:even-physical-iterate}.
Since $q:\mathcal J_n\to\mathcal J_{n-1}$ is proper of degree
$4$, the restriction maps properly onto
$\alpha^n\mathcal J_{n-1}$. Pulling back through $0$ increases
the depth $S_{n+1}-3$ by $S_n$, giving $S_{n+2}-3$.
This proves the induction.
\end{proof}

\begin{proposition}[Fibonacci combinatorics at every depth]
\label{id:prop:cycle-fibonacci-combinatorics}
The map $D$ is a real symmetric Fibonacci map. More precisely,
\begin{equation}\label{id:eq:full-fibonacci-grid}
 D^{\circ S_n}(0)\in
 C_{S_{n+1}-3}\setminus C_{S_{n+1}-2}
                           \qquad\text{for every }n\ge0.
\end{equation}
\end{proposition}
\begin{proof}
Because $1\in U^1\subset K(q)\subset\mathcal J_{n-1}$,
symmetry and \eqref{id:eq:all-critical-values} give
\[
 D^{\circ S_n}(0)=e_n\alpha^n
                         \in\alpha^n\mathcal J_{n-1}
                         =C_{S_{n+1}-3}.
\]
For the exclusion, $1\notin\alpha V=Z$, because $1\in U^1$
and the two open basin domains $U^1,Z$ are disjoint.
If $n=2k$, then
\[
 \alpha^{-n}C_{S_{n+1}-2}
                       =\alpha^2\mathcal A_k\subset\alpha V.
\]
If $n=2k+1$, the rescaled critical piece is exactly $\alpha V$.
Symmetry therefore excludes both $1$ and $-1$ in the required case,
proving \eqref{id:eq:full-fibonacci-grid}.

The proper-domain hypotheses and nonescape required for a generalized
polynomial-like map were established above. Condition
\eqref{id:eq:full-fibonacci-grid} is the Fibonacci marked-grid
characterization used in \cite[Section 4.2, discussion of
Proposition 7; Section 4.4, Step 4]{Buff}. It gives the stated
Fibonacci combinatorics. Here the condition has been proved at every
index, rather than inferred from a finite critical itinerary.
\end{proof}

\subsection{Renormalization germs and canonical identification}
\label{id:sec:cycle-identification}

The preceding proposition allows the canonical Fibonacci renormalization
to be applied. We verify equality of the resulting germs by comparing
both the formulas and the filled sets.

\begin{lemma}[The normalized restriction has the same filled set]
\label{id:lem:cycle-same-germ}
Let $R$ denote normalized Fibonacci renormalization, and let
$\mathfrak F$ change the sign of the lateral branch. Then
\[
                            [\mathfrak F R(D)]=[D]
\]
in the space of generalized polynomial-like germs.
\end{lemma}
\begin{proof}
By the canonical return construction
\cite[Section 4.2, Proposition 6]{Buff}, the first pre-renormalization
$\Psi$ has range $U^0$, central domain
$C_2=\alpha^2\mathcal J_1$, and lateral domain
$P_1(c_2)$, where $c_2=D^{\circ2}(0)=-\alpha$.
Its two formulas are
\[
 \Psi|_{C_2}(z)=D^{\circ2}(z)=-\alpha f(z/\alpha),
 \qquad \Psi|_{P_1(c_2)}(z)=D(z)=f(z).
\]
Normalize by the affine map $\lambda(z)=-\alpha z$.
The central formula of $R(D)=\lambda^{-1}\circ\Psi\circ\lambda$
is $f$, and its lateral formula is $-\alpha^{-1}q$.
After applying $\mathfrak F$, these are the two formulas of $D$.

The normalized range is
\[
 \widehat V=\lambda^{-1}(U^0)
           =\alpha\mathcal A_0=\mathcal J_0\Subset V.
\]
The normalized central domain is
\[
 \widehat U^0=\alpha\mathcal J_1
       \subset\alpha\mathcal J_0=\alpha^2\mathcal A_0=U^0.
\]
Its proper degree-four map onto $\mathcal J_0$ is the full
central pullback under $D|_{U^0}$.
For the normalized lateral domain,
\[
 \widehat U^1=\lambda^{-1}(P_1(c_2))
                    \subset\lambda^{-1}(U^0)
                    =\mathcal J_0\subset V\subset B.
\]
It is connected, contains $1$, and its formula after the flip is
$\alpha^{-1}q$, mapping properly onto $\mathcal J_0\subset V$.
The domain $U^1$ is the full component through $1$ of the inverse
image of $V$ under this formula in $B$. Proper-component recognition
therefore gives $\widehat U^1\subset U^1$; it is precisely the
pullback of $\mathcal J_0$ under the conformal branch
$D|_{U^1}$. Thus $\mathfrak F R(D)$ is the restriction of $D$
obtained by replacing its range with $\mathcal J_0$.

Finally, $U^0\subset Z$ and $U^1$ belong to the two invariant basin domains
for $q^{\circ2}$. Hence
\begin{equation}\label{id:eq:cycle-filled-set-containment}
                  U^0\cup U^1\subset K(q)\Subset\mathcal J_0.
\end{equation}
Every orbit in $K(D)$ stays in $U^0\cup U^1$ and consequently in
$\mathcal J_0$. Every such orbit is retained by the restriction.
Conversely, an orbit of the restriction is an orbit of $D$.
Their filled sets are therefore equal. Since the maps also coincide
on their common branch domains, they define the same germ. This
establishes the claim beyond equality of the formulas near the
marked points.
\end{proof}

\begin{proposition}[The actual period-two cycle]
\label{id:prop:independent-two-cycle}
The independently constructed map satisfies
\[
                              R^{\circ2}[D]=[D],
\]
and $R[D]\ne[D]$.
\end{proposition}
\begin{proof}
For real symmetric Fibonacci maps, the lateral flip is an involution
and commutes with normalized renormalization
\cite[Section 4.2, Lemmas 4--5]{Buff}. Therefore
\[
             (\mathfrak F R)^2=\mathfrak F^2R^2=R^2.
\]
Lemma~\ref{id:lem:cycle-same-germ} proves $R^2[D]=[D]$.
The second critical value of $D$ is $-\alpha$, whereas that of
$R(D)$ is $+\alpha$, since $q(1)=-\alpha^2$ and the latter map
has lateral formula $-\alpha^{-1}q$. These values are distinct.
Thus the period of the cycle is exactly two.
\end{proof}

\begin{proof}[Completion of the proof of Theorem~\ref{thm:canonical-identification}]
We have constructed a real symmetric degree-four Fibonacci map from
the validated germ, with proper physical domains, nonescaping critical
orbit, and a genuine period-two renormalization germ. Uniqueness of
that cycle in the space of real symmetric Fibonacci germs
\cite[Section 4.2, Theorem 3, citing \cite{NS}, Theorem 7.1]{Buff} now applies to this constructed
map. Its phase is the one with negative second critical value.
Equality of the normalized Fibonacci germs gives equality on their
filled set. The critical values $e_n\alpha^n$ in that set accumulate at
$0$, so the identity theorem identifies the central holomorphic germs.
The parameter is also identified, because
\[
                                  \alpha=-D^{\circ2}(0)
\]
in this normalization. This proves the identification claimed in
Theorem~\ref{thm:canonical-identification}. The argument uses the independently
certified polynomial-like attracting map and the exact functional
identity; it does not infer a global inverse-class property from a
local derivative sign.
\end{proof}

% END sections/07_identification.tex

% BEGIN sections/08_geometry.tex
\section{Finite enclosures of the canonical domains}\label{sec:geometry}
Write $E_\alpha=W_\alpha$, and retain $W_j$ exclusively for the
homothetic physical tower. Let
\[
 q_\alpha(z)=f(\alpha z)=\HH(z^4),\qquad
 G(z)=f(\alpha^2z)=\HH(\alpha^4z^4).
\]
Buff's canonical attracting basin $B$ of $G$ satisfies
$V=\alpha B$ and
\begin{equation}\label{eq:sandwich}
              \alpha E_\alpha\subset V\Subset E_\alpha.
\end{equation}
In the $V$ coordinate the basin map is
\begin{equation}\label{eq:FV}
                    F_V(v)=\alpha\HH(v^4),\qquad F_V(V)\subset V.
\end{equation}
This is conjugate to $G$ by multiplication by $\alpha$.

\begin{figure}[htbp]
\centering
\begin{tikzpicture}[manuscript diagram]
  % Nested domain boundaries are complete, with space for every label.
  \draw[diagramblue,fill=diagramblue!3,line width=.75pt]
    (-3.3,0) ellipse[x radius=2.2cm,y radius=1.45cm];
  \draw[diagramteal,fill=diagramteal!4,line width=.7pt]
    (-3.3,0) ellipse[x radius=1.65cm,y radius=1.03cm];
  \draw[diagramorange,fill=diagramorange!11,line width=.7pt]
    (-3.3,0) ellipse[x radius=1.03cm,y radius=.60cm];
  \node at (-3.3,1.23) {$E_\alpha$};
  \node at (-3.3,.82) {$V=\alpha B$};
  \node at (-3.3,.35) {$\alpha E_\alpha$};
  \fill (-3.3,-.08) circle[radius=1pt];
  \node[font=\scriptsize,below right=2pt] at (-3.3,-.08) {$0$};
  \node[font=\footnotesize] at (-3.3,-1.86) {Canonical inclusions};
  \node (B1) at (1.10,1.05) {$B$};
  \node (B2) at (4.45,1.05) {$B$};
  \node (V1) at (1.10,-1.05) {$V$};
  \node (V2) at (4.45,-1.05) {$V$};
  \draw[map arrow] (B1) -- node[map label,above=3pt] {$G$} (B2);
  \draw[map arrow] (V1) -- node[map label,below=3pt] {$F_V$} (V2);
  \draw[map arrow] (B1) -- node[map label,left=3pt] {$z\mapsto\alpha z$} (V1);
  \draw[map arrow] (B2) -- node[map label,right=3pt] {$z\mapsto\alpha z$} (V2);
  \node[font=\footnotesize] at (2.78,-1.86) {Change of coordinate};
\end{tikzpicture}
\caption{Canonical domain inclusions and the exact change of coordinate $v=\alpha z$, conjugating $G$ to $F_V$. The inclusions are schematic; Theorem~\ref{thm:canonical-identification} identifies these objects for the validated germ.}\label{fig:basin-conjugacy}
\end{figure}
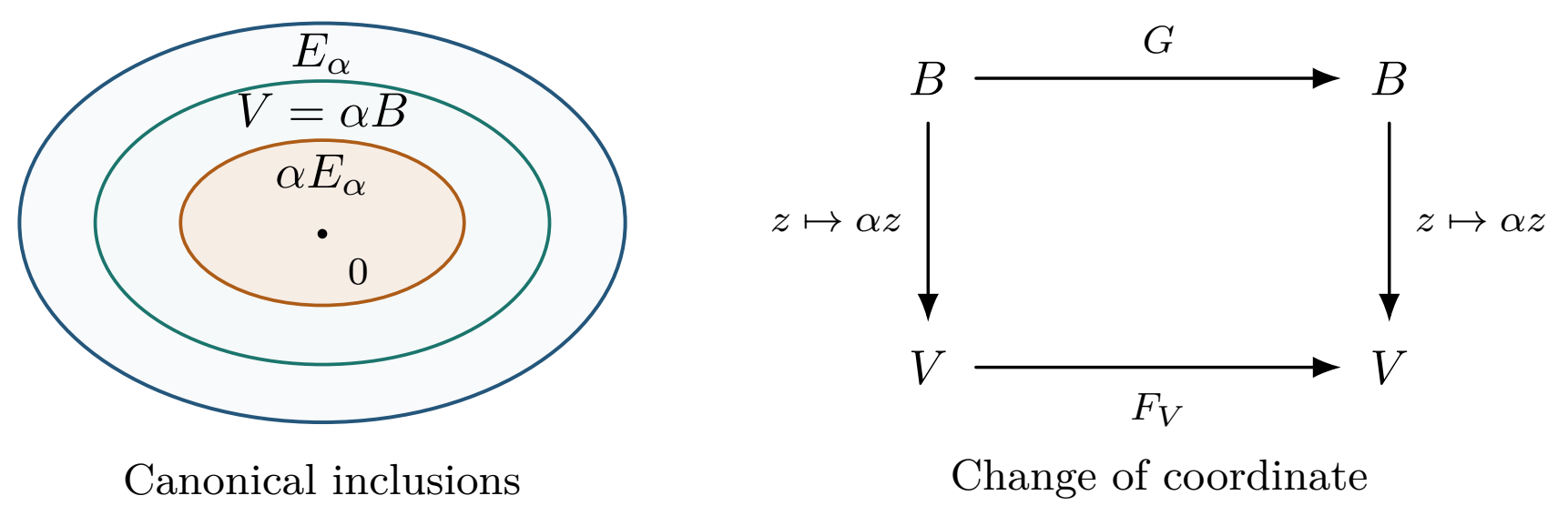
\subsection{An outer source polygon by an interior crossing argument}
The exact functional identity is
\begin{equation}\label{eq:linearizer}
                       q_\alpha(G(z))=-\alpha^2q_\alpha(z).
\end{equation}
Let $T(z)=\overline{G(z)}$ and let
$S=\{x+iy:0\le y\le x\}$. The Epstein inverse gives a unique tile
$P\subset\operatorname{int}S$ on which
$q_\alpha:P\to\mathbb H^-$ is conformal. Its inverse will be denoted
$\psi$. The sign follows from $\HH'(0)<0$ and the fourth-root lift.
The eight rotations and reflections of this tile are dense in
$E_\alpha$.

We have $G(E_\alpha)\subset E_\alpha$. Indeed
$E_\alpha\subset B$ and $G(B)\subset B$, while
$-\alpha^2\Delta_0\subset\Delta_0$. Properness makes $E_\alpha$
the component of $q_\alpha^{-1}(\Delta_0)$ in $B$ containing $0$
and $1$. The connected image $G(E_\alpha)$ lies in that same
component by \eqref{eq:linearizer} and $G(0)=1$.

It follows that $T$ maps $P$ onto itself and
\begin{equation}\label{eq:tile-T}
 T=\psi\circ\bigl(w\mapsto-\alpha^2\overline w\bigr)\circ q_\alpha
                                      \quad\text{on }P.
\end{equation}
To identify the image tile, use the expansion near the origin in $S$:
$T(z)=1+\alpha^4\HH'(0)\overline z^4+O(|z|^8)$ lies in the first
octant. Surjectivity follows because $q_\alpha(P)=\mathbb H^-$ and
the middle map in \eqref{eq:tile-T} is an automorphism of that half-plane.
The inverse $\psi$ extends analytically at $w=0$ to the positive
regular zero $\rho$ of $q_\alpha$.

\begin{lemma}[Interior exit-gate lemma]\label{lem:outer-gate}
Let $K\subset S$ be a closed polygon with $\rho$ in its relative
interior along the real side. Let $\Gamma$ be its outer edges, and
suppose every path in $S$ from a sufficiently small relative
neighborhood of $\rho$ to $S\setminus K$ crosses $\Gamma$.
If $T$ is defined on $\Gamma$ and $T(\Gamma)\cap K=\varnothing$,
then $P\subset K$.
\end{lemma}
\begin{proof}
If $P\not\subset K$, the curves
\[
              u\longmapsto\psi(uq_\alpha(z)),\qquad0<u\le1,
\]
from $\rho$ to points $z\in P\setminus K$ show that
$\Gamma\cap P\ne\varnothing$. Put
$d=\inf\{|q_\alpha(z)|:z\in\Gamma\cap P\}$.
Then $d>0$: otherwise $q_\alpha(z_j)\to0$ would imply
$z_j=\psi(q_\alpha(z_j))\to\rho$, contradicting compactness of
$\Gamma$ and $\rho\notin\Gamma$.
Choose $z\in\Gamma\cap P$ with $|q_\alpha(z)|<d/\alpha^2$.
By the gate hypothesis and \eqref{eq:tile-T}, $Tz\in P\setminus K$
and $|q_\alpha(Tz)|=\alpha^2|q_\alpha(z)|<d$.
The radial-coordinate curve from $\rho$ to $Tz$ crosses the gate
at a point $z'\in\Gamma\cap P$ with
$|q_\alpha(z')|\le|q_\alpha(Tz)|<d$, a contradiction.
No assertion about a parametrization of the entire geometric boundary
of $P$ is needed.
\end{proof}

Apply this lemma to
\begin{equation}\label{eq:OABC}
 K=\operatorname{conv}\left\{(0,0),\left(\frac{41}{20},0\right),
       \left(\frac{39}{20},\frac{21}{20}\right),
       \left(\frac{71}{50},\frac{71}{50}\right)\right\}.
\end{equation}

\begin{figure}[htbp]
\centering
\begin{tikzpicture}[manuscript diagram,x=2.35cm,y=2.35cm]
  \coordinate (O) at (0,0); \coordinate (A) at (2.05,0);
  \coordinate (B) at (1.95,1.05); \coordinate (C) at (1.42,1.42);
  \draw[->,gray!60,line width=.45pt] (-.13,0) -- (2.45,0)
    node[right,font=\scriptsize] {$\Re z$};
  \draw[gray!55,line width=.45pt] (O) -- (1.78,1.78);
  \node[font=\scriptsize,text=gray!75] at (1.63,1.92) {$\arg z=\pi/4$};
  \fill[diagramblue!6] (O)--(A)--(B)--(C)--cycle;
  \draw[diagramblue,line width=.6pt] (A)--(O)--(C);
  \draw[diagramorange,line width=1pt] (A)--(B)--(C);
  \node[text=diagramorange] at (2.12,.66) {$\Gamma$};
  \node at (1.48,.48) {$K$};
  \foreach \point/\placement in {O/below left,A/below right,B/right,C/above left}
    {\fill (\point) circle[radius=.65pt];\node[\placement=2pt] at (\point) {$\point$};}
  \fill (.72,0) circle[radius=.7pt];
  \node[below=3pt] at (.72,0) {$\rho$};
  \draw[map arrow,diagramteal] (.72,0)
    .. controls (1.12,.35) and (1.45,.93) .. (1.71,1.35);
  \node[font=\scriptsize,text=diagramteal] at (.68,.90) {$\psi(uw)$};
\end{tikzpicture}
\caption{The polygon $K$ and its tested edges $\Gamma$, with the exact vertices in \eqref{eq:OABC}. The path and $\rho$ are schematic. Lemma~\ref{lem:outer-gate} rules out an exit when $T(\Gamma)\cap K=\varnothing$. Its crossing point lies in $\Gamma\cap P$.}\label{fig:exit-gate}
\end{figure}
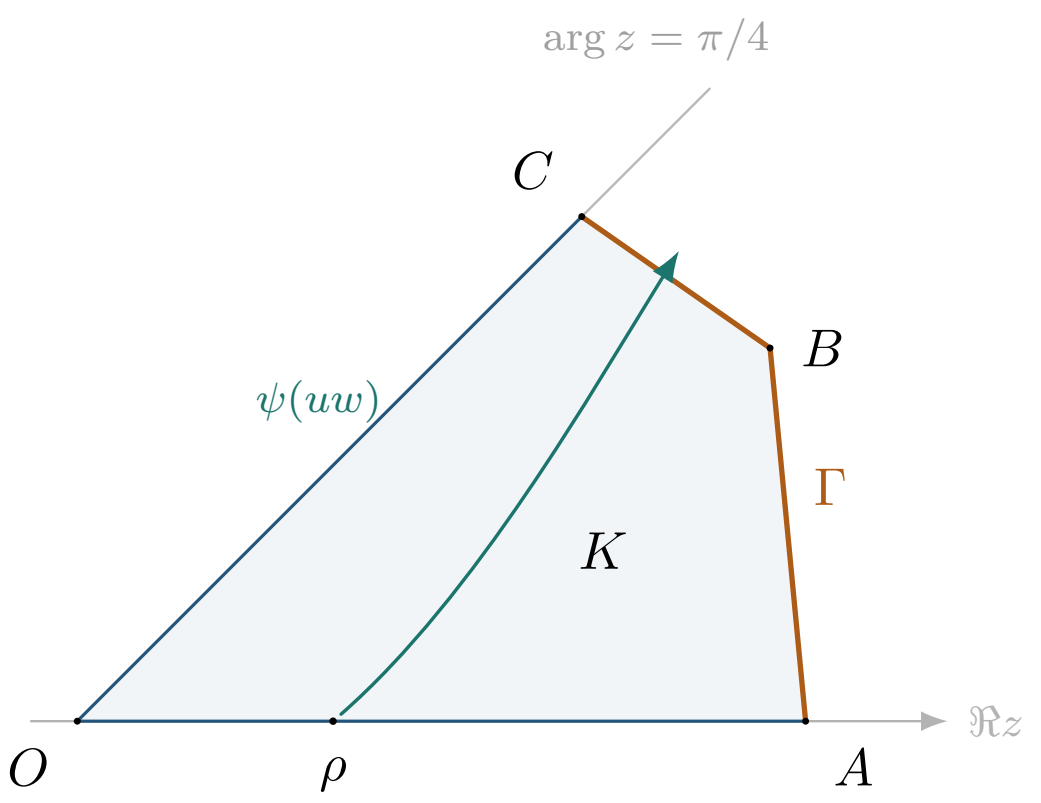
The two outer edges are the gate; the other two lie on the sector
boundary. Each outer edge is divided into $2^{10}$ segments.
The exact certificate evaluates $T$ on all 2,048 containing boxes
and gives a strictly violated half-plane inequality of $K$ for
each entire image. Its minimum oriented margin exceeds
$0.00075127607697$. All evaluations use strictly valid local arguments.
Since $0<\rho<1$, the start of each radial-coordinate path lies in
$K$ relative to $S$. Lemma~\ref{lem:outer-gate} and symmetry yield
\begin{equation}\label{eq:E-outer}
 \overline{E_\alpha}\subset E^{\rm out},
 \qquad E^{\rm out}=\bigcup\{\text{eight dihedral copies of }K\}.
\end{equation}
The exact upper-half area of this initial enclosure is $6861/1000$.

\subsection{Invariant pruning of the outer basin cover}
Start with the closed polygonal cover $E^{\rm out}$ from
\eqref{eq:E-outer}. Partition it by dyadic squares, clipping the
boundary pieces exactly. Suppose a current closed cover $C$ contains
$\overline V$. A square $Q$ can be removed if a valid complete
evaluation gives
\begin{equation}\label{eq:outer-pruning}
                     F_V^{\circ j}(Q)\cap C=\varnothing
\end{equation}
for some finite $j$. A hypothetical point in $Q\cap\overline V$
would remain in $\overline V\subset C$, contrary to the enclosure.
Analytic guards and continuity justify this also on the closures.
Unknown evaluations retain the square.

After each round retain the origin component of the closed-cell
adjacency graph. Eight-neighbor adjacency includes corner contacts;
ignoring polygon clipping can only introduce additional connections.
Connectedness of $V$ ensures that no part of $V$ is lost. This
finite pruning can be repeated with the improved cover.
The final mesh-$2^{-8}$ cover has upper-half area
\begin{equation}\label{eq:Vouter-area}
                    \area(P_0^+)\le\frac{73671}{32768}.
\end{equation}
The 37,708 initial image enclosures and seven pruning rounds have
complete exact rational replay records.

\subsection{Inner source and basin enclosures}
A connected closed set through the mark $0$, with complete
$q_\alpha$ image in $\Delta_0$, belongs to the marked source
$E_\alpha$ by proper-component recognition. To certify a source disk,
the computation covers a full circle in the quotient $x$ coordinate
and proves that $\Re\HH(x)$ lies strictly between conservative slit
endpoints. The harmonic maximum and minimum principles prove the
same strict inequalities in the disk. For radius $1.39$ in the $z$
coordinate, 956 quotient boundary boxes suffice. Additional connected
grid cells use complete slit-image checks. Their union subsequently
contains the closed disk of radius $3/2$.

Equation \eqref{eq:sandwich} transfers these source subsets to
inner subsets of $V$. A new whole square is an entry cell if all
its pre-entry $F_V$ enclosures lie in the proved physical source
and a finite terminal enclosure is in an already known inner subset
of $V$. This proves membership in the full attracting basin. Only
the connected collection attached to the known immediate component
is retained; shared closed edges lie in the open basin and justify
the attachment.

The final mesh-$2^{-9}$ inner grid and its exact area give
\begin{equation}\label{eq:Vinner-area}
 \frac{70987}{32768}\le\area(P_0^+)
                   \le\frac{73671}{32768}.
\end{equation}
It contains the closed centered disk of radius $197259/200000$.
The disk claim is an exact grid inclusion: the nearest omitted
first-octant cell has squared minimum radius $255008/512^2$,
strictly larger than $(197259/200000)^2$.

\subsection{Physical branch recognition and outer isolation}
\begin{lemma}[Proper-source recognition]\label{lem:physical-recognition}
Let $F_j$ be a holomorphic continuation of the physical proper map
$D_j:U^j\to V$. If a connected set $A$ lies in its continuation
domain, $F_j(A)\subset V$, and $A$ meets $U^j$, then $A\subset U^j$.
\end{lemma}
\begin{proof}
Properness makes $U^j$ relatively closed as well as open in the
component of $F_j^{-1}(V)$ containing its mark. It is that component.
The connected set $A$ cannot leave it.
\end{proof}

In the finite calculation $A$ is the rectangular hull of a query
box and a point in the authenticated source grid. Either complete
image containment in the inner $V$ grid is proved directly, or a
finite subdivision proves that every image point enters $V$ under
$F_V$. In the latter case the connected image $F_j(A)$ lies in the
full basin and meets its immediate component at the anchor image;
it therefore lies in $V$. Lemma~\ref{lem:physical-recognition} applies.
Different children may have different finite entry times.

Outer source polygons are proved by a first-exit argument. A
holomorphic branch formula is evaluated on the full boundary of
a polygon containing its physical mark. If every boundary image
misses a closed outer cover of $V$, no path in $U^j$ can cross that
boundary. Continuity also excludes contact of $\overline{U^j}$.
The final central polygon is the 16-gon generated by rotating
\[
 (0.64,0),\quad(0.59,0.245),\quad(0.45,0.45),\quad(0.245,0.59),
\]
and the lateral polygon is $[0.52,1.15]+i[-0.30,0.30]$.
The 5,120 boundary checks have exact replay. These are outer
isolators; membership in them does not certify a physical branch.
The inner source masks at mesh $2^{-10}$ have full-plane areas
$149883/262144$ and $17785/262144$, respectively.

% END sections/08_geometry.tex

% BEGIN sections/09_algorithms.tex
\section{Whole-set landing and loss certification}\label{sec:algorithms}
Theorem~\ref{thm:canonical-identification} and the preceding enclosures
supply a fixed canonical map and tests on its physical domains. LAND
constructs an upper bound by retaining every unresolved cell. LOSS
constructs a lower bound by accepting only whole cells with a witness.

Both bounds use two complementary methods. In the Koebe method a
terminal disk in the physical hole is pulled back along a certified center
itinerary; a source cell is certified when it fits inside that inverse
disk. Direct propagation advances an enclosure of the entire cell until
a terminal hole test succeeds. The counted objects remain dyadic cells:
the area of a containing ball is never substituted for the cell's area.
LAND first forms a complete retained cover by Koebe exclusions and then
refines that cover or selected descendants by direct propagation.
LOSS first forms a Koebe-certified union $E_K$ and then adds direct
witnesses by geometric union, removing all intersections before adding
area. Thus Corollary~\ref{cor:koebe-land} is used in the actual LAND
dependency chain; the later direct algorithms reuse its retained cover.

\subsection{Quantities required by the uniform criterion}
Recall $L_r^\tau=h^{-1}\widetilde X_r\cap P_0^\tau$,
$a_r=\area(L_r^+)$ and $\area(\widetilde X_r)=2\alpha^4a_r$.
All numerical source cells use these normalized coordinates. The loss
normalization retains its first index. Explicitly, put
\[
 \mathcal Y_n=h_n^{-1}Y_n
 =\{z\in V:D^k(h_nz)\notin h_nV\text{ for every defined }k\ge1\}.
\]
Parent conjugacy gives
\begin{equation}\label{eq:loss-normalized-event}
 \begin{split}
 A_{n,m}^\tau=\{z\in P_0^\tau:\ &D^q(hz)\in h\mathcal Y_n
                  \text{ for some defined }q\ge0,\\
 &D^j(hz)\notin W_{m+1}\text{ for every defined }j\ge0\}.
 \end{split}
\end{equation}
Indeed $\cR_{n-1}^j(h_nz)=h_{n-1}D^j(hz)$ and
$h_{n-1}^{-1}Y_n=h\mathcal Y_n$. The initial parent itinerary is
normalized, but the terminal event still concerns ordinary iterates at
level $n$. For $(n,m_0,t)=(3,2,5)$ the four required outputs are
\[
 a_3\le A_3^+,\quad E^\tau\subset A_{3,2}^\tau,\quad
 \area(L_5^\tau\setminus C_5^\tau)=0,\quad
 g_{E^\tau}(z)\ge\gamma^->0\text{ a.e. on }C_5^\tau.
\]
The support may extend outside $V$, provided every required kernel row
lies in its slit domain. No forward invariance of the support is needed.

\subsection{Exact algorithms and finite search parameters}
Even with exact set operations, the algorithm uses finite subdivision
and a finite orbit horizon. An unresolved LAND cell remains in the upper
cover, whereas an unresolved LOSS cell is absent from the accepted union.
Let $L$ be the mesh level, so a square has side $2^{-L}$; it is distinct
from dynamical depths $n,m,t$ and from Fibonacci times $S_r$.
The production LAND horizon $K_{\rm phys}=500$ counts ordinary iterates
of $D$. If block $i$ has duration $t_i$, its cumulative clock is
$k_i=\sum_{j<i}t_j$; a block is admitted only if
$k_i+t_i\le K_{\rm phys}$. A bound on the number of blocks would be a
different parameter and is not the meaning of 500 in these records.
The LOSS bounds $K_1,K_2$ count parent iterates and ordinary tail iterates,
respectively. Auxiliary iterations of the basin map $F_V$, used only in
recognition, have separate limits. No one of these finite cutoffs proves
escape without a terminal witness.

\subsection{Exact cemetery-hole pullbacks}
Set $U=U^0\sqcup U^1$ and
$\mathcal O=V\setminus(\overline{U^0}\cup\overline{U^1})$.
An orbit entering this open hole has no further physical iterate.

\begin{lemma}[Cemetery-hole nesting]\label{lem:cemetery}
Let $\Delta\Subset\mathcal O$ be a disk and $A$ a component of
$\Dom(D^{\circ k})\cap D^{-k}(\Delta)$. For any critical puzzle
piece $C_d$, either $A\cap C_d=\varnothing$ or $A\subset C_d$.
Every selected physical pullback over $\Delta$ is conformal.
\end{lemma}
\begin{proof}
The case $k=0$ is immediate. If $a\in A\cap C_d$ and $1\le k\le d$,
the depth-$d$ itinerary gives $D^{\circ k}(a)\in U$, contradicting
membership in $\Delta$. Thus $k\ge d+1$. The component of
$\Dom(D^{\circ k})$ containing $A$ has depth $k-1\ge d$ and meets
$C_d$, so it is contained in $C_d$ by the Markov property.
Finally every critical value of a defined physical composite is
postcritical and belongs to $K(D)\subset U$. Hence $\Delta$ contains
no such value. Properness and simple connectivity make its selected
unbranched pullbacks conformal.
\end{proof}

\begin{corollary}[Koebe landing exclusion]\label{cor:koebe-land}
Let $x,Dx,\ldots,D^{\circ k}x$ be a physical orbit whose points
avoid $W_n$, and suppose $B_R(D^{\circ k}x)\Subset\mathcal O$.
Then the disk
\[
                  B_{R/(4|(D^{\circ k})'(x)|)}(x)
\]
is disjoint from the global physical landing set of $W_n$.
For $k=0$ the conclusion with the smaller radius $R/4$ is also valid.
\end{corollary}
\begin{proof}
Pull the terminal disk backward through the physical itinerary.
If any suffix pullback met $W_n$, Lemma~\ref{lem:cemetery} would
put the whole suffix, including its center, in $W_n$, a contradiction.
No point in the terminal disk can later return. The inverse is
univalent, and Koebe's one-quarter theorem gives the stated disk.
When $R$ is a distance supremum, first use every radius smaller
than $R$ and then take the union.
\end{proof}

For a normalized square centered at $c$, apply the corollary with
$x=hc$. If $r_Q$ is its circumradius, a sufficient exclusion test is
\begin{equation}\label{eq:koebe-test}
 R>4r_QJ,\qquad J\ge |h|\,|(D^{\circ k})'(hc)|.
\end{equation}
The center $c$ need not belong to $V$: the global exclusion disk
excludes its intersection with the true normalized source.
Its physical itinerary must nevertheless be certified. The same
disk may be reused for any smaller cell whose farthest corner is
strictly inside it. Exact dyadic subdivision and disk reuse preserve
every cell not covered by a proved exclusion.

\subsection{Exact landing algorithm}
Let $\mathcal Q_0$ be a finite family with disjoint interiors covering
$P_0^+$ modulo boundaries. For $Q\in\mathcal Q_0$ put
$S_0=h(Q\cap V)$, intersecting with the true domain. A direct exclusion
witness consists of exact sets $S_0,\ldots,S_N$, integers $t_i\ge1$ and
physical holomorphic blocks $F_i$ such that
\begin{gather*}
 F_i=D^{t_i}\text{ near }S_i,\quad S_{i+1}=F_i(S_i)\quad(0\le i<N),\\
 S_i\cap W_n=\varnothing\quad(0\le i\le N),\\
 D^r(S_i)\cap W_n=\varnothing\quad(0\le i<N,\ 0<r<t_i),
 \qquad S_N\subset\mathcal O.
\end{gather*}
Every displayed preterminal iterate must be defined on the entire exact
set. These conditions imply $Q\cap L_n^+=\varnothing$.
Equivalently, starting with $(S,k)=(S_0,0)$:
\begin{enumerate}[label=(L\arabic*)]
\item If target avoidance cannot be proved, retain the cell or subdivide it.
\item If $S\subset\mathcal O$, exclude the cell: all previous images
avoided the target and the entire orbit has terminated.
\item Choose a common physical branch, including every intermediate-time
avoidance, advance to its whole-set image, and add its duration to $k$.
A chosen block must fit the remaining physical budget.
\item If no common branch is certified or a work limit is reached, retain
the cell, or replace it by all children when further subdivision is allowed.
\end{enumerate}
Ordinary branches have $t_i=1$. An accelerated block
$F_i=h_{r_i}D_{j_i}h_{r_i}^{-1}$ requires
$S_i\subset h_{r_i}U^{j_i}$ with $0\le r_i\le n-1$.
For $r_i\ge1$ its duration is $S_{2r_i}$ or $S_{2r_i-1}$;
Lemma~\ref{lem:first-return} proves all proper intermediate exclusions
because $W_n\subset W_{r_i}$. For $r_i=0$ its duration is one.
Thus $t_i$ is a block duration, while $k_i=t_0+\cdots+t_{i-1}$ is the
cumulative physical time and $S_i=D^{k_i}(S_0)$. Different source cells
generally select different sequences. Setting every $t_i=1$ and a maximum
ordinary time $M$ is valid, but avoidance for $k<M$ without a terminal
hole witness does not prove nonlanding after $M$. Acceleration reduces
repeated evaluations and growth of interval enclosures.

At the finite horizon a source cell is excluded only by this direct
witness or Corollary~\ref{cor:koebe-land}. A cell proved to land is
retained as well; such a test is optional for an upper bound. If
$\mathcal Q$ is the retained family, then
\[
 L_n^+\subset\bigcup_{Q\in\mathcal Q}Q\text{ modulo null boundaries},
 \qquad a_n\le A_n^+=\area\Bigl(\bigcup_{Q\in\mathcal Q}Q\Bigr).
\]
Every unselected root remains in the union. The retained area of a
selected batch alone is not an upper bound for $a_n$.
If only outer source information is available, retain every possible
physical branch choice. A branch of this calculation follows the exact
subset $S\cap U^j$, using an enclosing image if necessary. A path is
eliminated only by proved impossibility or a proved terminating exclusion;
a failed analytic chart retains it. A cell is excluded only after all
possible physical paths have been excluded without a target visit.
Analytic continuation bounds physical subsets but does not introduce new
physical orbits.

\subsection{Exact parent-saturated loss algorithm}
For $Q\subset P_0^+$, start from $x=h_nz$ with $z\in Q$ and use
the parent coordinates
\begin{equation}\label{eq:loss-coordinates}
 u_0=hz,\qquad u_j=D^{\circ j}(u_0),\qquad
 x_j=h_{n-1}u_j=\mathscr R_{n-1}^{\circ j}(x).
\end{equation}
For a whole cell $Q\subset P_0^+$ write
$\mathcal U_j=D^j(hQ)$ for $0\le j\le q$. The exact first stage requires
\[
 \mathcal U_j\cap W_{m+1}=\varnothing\ (0\le j\le q),\qquad
 \mathcal U_j\subset U^{\epsilon_j}\ (0\le j<q),\qquad
 \mathcal U_q\subset W_1.
\]
For $T_0=h_{n-1}\mathcal U_q\subset W_n$, the second stage requires
\[
 T_r=D^r(T_0),\quad T_r\cap W_n=\varnothing\ (1\le r\le N_{\rm tail}),
 \quad T_{N_{\rm tail}}\subset\mathcal O.
\]
Every preterminal step is physical on the whole set. Choose finite
$q\le K_1$, $1\le N_{\rm tail}\le K_2$. These conditions give
$T_0\subset Y_n$, and hence $Q\subset A_{n,m}^+$. At $q=0$ the
initial exclusion $Q\cap W_m=\varnothing$ remains necessary.
An excluded LAND cell can be a candidate, but LOSS also requires the
whole cell in $P_0^+$ and the tail beginning at
$h_{n-1}D^q(hQ)$, with the stated positive-time avoidance. A LAND
exclusion alone supplies neither condition.

\begin{figure}[htbp]
\centering
\begin{tikzpicture}[manuscript diagram]
  \node[cell box,minimum width=2.4cm] (start) at (0,0) {$\mathcal U_0=hQ$};
  \node[cell box,minimum width=2.7cm] (parent) at (7.6,0) {$\mathcal U_q\subset W_1$};
  \node[cell box,minimum width=2.7cm] (tail) at (7.6,-2.65)
    {$T_0=h_{n-1}\mathcal U_q$\\$\subset W_n$};
  \node[cell box,minimum width=2.4cm,draw=diagramorange,fill=diagramorange!5]
    (hole) at (0,-2.65) {$T_{N_{\rm tail}}\subset\mathcal O$};
  \draw[map arrow] (start) --
    node[map label,above=5pt] {$q\le K_1$ iterates of $D$}
    node[map label,below=5pt] {Avoid $W_{m+1}$\\at every parent checkpoint} (parent);
  \draw[map arrow] (parent) -- node[map label,right=5pt] {$h_{n-1}$} (tail);
  \draw[map arrow] (tail) --
    node[map label,above=5pt] {$N_{\rm tail}\le K_2$ ordinary iterates}
    node[map label,below=5pt] {Avoid $W_n$ at every\\strictly positive tail time} (hole);
\end{tikzpicture}
\caption{The two parts of a LOSS witness. The upper row is in normalized parent coordinates; the lower row uses ordinary physical iterates. Every box denotes the image of a whole source cell. The parent and tail clocks have separate finite bounds.}\label{fig:loss-witness}
\end{figure}
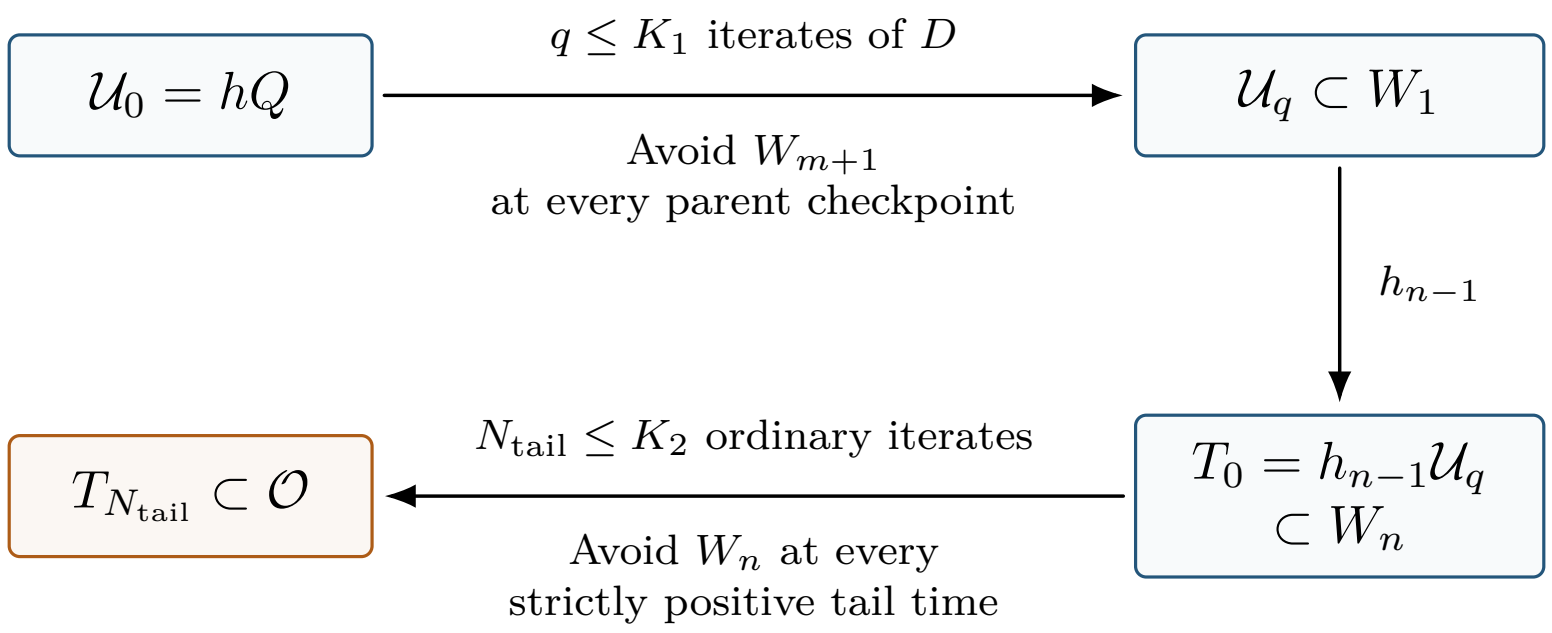
The same two stages in point coordinates are as follows.

\begin{enumerate}
\item Follow a common admissible parent itinerary through index $q$,
requiring $u_j\notin W_{m+1}$ for every $0\le j\le q$.
At a candidate index require $u_q\in W_1$, equivalently $x_q\in W_n$.
\item Starting at $w=x_q$, prove that every strictly positive
physical time through a finite tail length $N_{\rm tail}$ avoids $W_n$
and that its terminal image is in $\mathcal O$.
\end{enumerate}
Then $w\in Y_n$. Earlier parent checkpoints avoided $W_{n+m}$;
no later parent checkpoint can reach it, because this would be a
positive-time physical return to $W_n$ after $w$. The entire square
is therefore contained in $A_{n,m}^+$ after the prescribed sign
normalization. Accepted loss squares are required to lie wholly
inside the true source $P_0^+$.

Both stages can be certified by direct whole-set propagation.
They can also be certified from a center orbit and a terminal disk,
using the following consequence of Lemma~\ref{lem:cemetery}.

\begin{lemma}[Koebe loss inclusion]\label{lem:koebe-loss}
Suppose the center orbit satisfies the two-stage conditions above
and ends in a disk $\Delta=B_R(D^{\circ N_{\rm tail}}w)\Subset\mathcal O$.
Let
\[
 J\ge |h_n|\prod_{j<q}|D'(u_j)|
                     \prod_{r<N_{\rm tail}}|D'(D^{\circ r}w)|.
\]
The normalized disk of radius $R/(4J)$ about the source center
belongs to the normalized loss event. Consequently a whole square
in $P_0^+$ satisfying $R>4r_QJ$ is an accepted loss square.
\end{lemma}
\begin{proof}
Pull $\Delta$ backward along the complete physical word.
At the candidate checkpoint the suffix contains $w\in W_n$,
so Lemma~\ref{lem:cemetery} puts it wholly in $W_n$. At each
strictly positive tail time it avoids $W_n$, since an intersection
would put its center there. Thus this entire suffix lies in $Y_n$.
At every parent checkpoint the remaining pullback is contained in
$W_{n-1}$, because its center is there. Pulling back one more parent
block therefore uses the declared full physical return component.
It avoids $W_{n+m}$: an intersection would, by the same lemma, put
its center in that forbidden piece. Hence the entire chained
preimage has the asserted parent-saturated loss property.
The terminal disk contains no physical critical value, so its
selected inverse is conformal. Koebe and the chain rule prove
the radius bound. In the derivative of each parent conjugate the
two linear factors cancel. There is no additional factor
$|h_{n-1}|$ in $J$.
\end{proof}

The index $q$ counts parent iterates; $N_{\rm tail}$ counts ordinary physical
iterates. For $n=3$ the parent branch times are $8$ and $5$.
At $q=0$, $u_0=hz$ is automatically in $W_1\subset U^0$.
Its first parent central block has endpoint $h_{n-1}D_0(hz)$ and
time $S_{2n-2}$. All its proper positive prefixes avoid $W_n$.
The initial forbidden test is $hz\notin W_{m+1}$, or
$z\notin W_m$; it must not be omitted.

\subsection{Interval implementation of the exact tests}
The mesh describes the source cells, not their image enclosures. A square
of side $2^{-L}$ lies in a disk of radius $2^{-L}/\sqrt2$. A whole-ball
certificate therefore certifies the square. Nonlinear images are generally
neither squares nor disks; intersections of ball and rectangle enclosures
can be much tighter than either alone. All exact orbit points, including
rescaled dyadic centers, are
represented by interval rectangles. A dyadic center is exact, but
$h_nc$, subsequent values of $D$, and the parameter $\alpha$ are
not binary rational point data. Every derivative and set image
uses the complete current enclosure and the analytic tail bound.

For the final LAND and additional LOSS computations the state
consists of a rectangle $R$ containing the actual image $S$, a
rectangle $Z$ containing an exact analytic reference point $c$,
and a bound $|z-c|\le\delta$ for all $z\in S$.
Replace $R$ by $R\cap(Z+[-\delta,\delta]^2)$.

\begin{lemma}[Hybrid propagation]\label{lem:hybrid}
Let $S\subset R$ be the exact current orbit set and let
$F=h_rD_jh_r^{-1}$ agree with a physical block on all of $S$. Suppose its analytic extension is defined
on the convex hull of $R\cup Z$ and $|F'|\le K$ there.
Neither $R$ nor the reference rectangle $Z$ must lie wholly in the
physical source; the exact set $S$ must follow the declared branch. Then
\[
 F(S)\subset F(R)\cap B_{K\delta}(F(c)).
\]
If $\mathbf J$ encloses $F'$ on that hull and $Z'$ encloses $F(c)$,
the further rectangular intersection
\begin{equation}\label{eq:centered-image}
                    F(R)\cap\bigl(Z'+\mathbf J(R-Z)\bigr)
\end{equation}
also contains $F(S)$.
\end{lemma}
\begin{proof}
For $z\in S$, the segment from $c$ to $z$ lies in the derivative
hull. Integrating $F'$ gives $|F(z)-F(c)|\le K\delta$.
The mean derivative on the segment belongs to the rectangular,
convex enclosure $\mathbf J$, giving \eqref{eq:centered-image}.
The reference point is used analytically and need not be physical;
the whole actual set follows the certified physical block.
\end{proof}

Initially $S=h(Q\cap V)\subset W_1\subset U^0$. Thus a direct
central fallback is valid for the actual set even if the analytic
reference center has not been recognized in $V$; its analytic
derivative hull must still be valid. Every admitted block has range
in some $W_r\subset V$, so the invariant $S\subset V$ persists.
This permits a terminal test in which the containing rectangle is
disjoint from both source closures even though the rectangle itself
is not contained in the available inner $V$ grid. The actual set
is then in $\mathcal O$. The separate disk test requires a proved
disk in $\mathcal O$ and the strict inequality $R_{\rm hole}>\delta$.

Source-closure exclusion is either geometric, using the outer
isolators, or analytic. If a complete $D_j$ image has a finite
$F_V$ iterate disjoint from the closed outer cover of $V$, it
misses $\overline V$. Continuity and $D_j(U^j)\subset V$ then
exclude $\overline{U^j}$ from the original box. Every finite
analytic evaluation has strict domain guards. The final LAND
implementation uses at most six such basin iterations and at most
ten basin-entry iterations in source recognition.

\paragraph{Lower distance bounds in disk tests.}
For a forbidden closed set $F$ and an endpoint enclosure $Z\ni z_{\rm ex}$,
a lower bound $d^-$ must satisfy
\[
 0\le d^-\le\inf_{z\in Z}\operatorname{dist}(z,F)
             \le\operatorname{dist}(z_{\rm ex},F).
\]
For a terminal disk in $\mathcal O$, $F$ must contain the complement of a
proved inner subset of $\mathcal O$. Distance to the source domains alone
does not establish containment in $V$. The direct whole-set terminal test
has the separate exact invariant $S\subset V$.

\subsection{Exact area bookkeeping and finite stopping rules}
A square of side $2^{-L}$ contributes $2^{-2L}$; a diagonal
first-octant triangle contributes half that amount. If $u$ is the
integer number of retained half-square units at a common finest
level $L$, the upper-half area is exactly
\begin{equation}\label{eq:dyadic-area}
                              \frac{2u}{4^L}.
\end{equation}
Mixed-depth cells are converted to these units by integer powers
of four. Each selected root is replaced by all of its children,
and every unselected root and unresolved child is retained.
Intersections of valid outer covers are computed geometrically;
unions of loss cells remove overlaps exactly before summing areas.

More explicitly, for a complete initial LAND cover $C^{(0)}$ and the
certified excluded union $B_s$ in pass $s$,
\[
 C^{(s)}=C^{(s-1)}\setminus B_s.
\]
Each tested descendant in the production chain is inside a retained root
of $C^{(s-1)}$. Only previously retained area can be removed. The residual
is already a union of dyadic cells and diagonal half-cells, so no new
outer approximation or boundary-area error is introduced between passes.

For LOSS, put $E^{(0)}=E_K$, and let $B_s$ be the raw accepted union of
pass $s$. The new accepted set is
\begin{equation}\label{eq:loss-union-update}
 N_s=B_s\setminus E^{(s-1)},\qquad
 E^{(s)}=E^{(s-1)}\sqcup N_s.
\end{equation}
Thus the $N_s$ are mutually disjoint and disjoint from $E_K$.
Within a batch, duplicate coordinates are removed and the dyadic union is
coalesced into leaves with disjoint interiors. For a partially covered
candidate it is valid to certify the whole cell and add only its uncovered
part; a fully covered cell can be skipped. All operations use exact dyadic
coordinates. Raw accepted areas from different passes are never added
without these differences.

The fourfold upper-half symmetry used for these dyadic areas is
valid for the actual normalized events: $V$ and $W_j$ have quartic
symmetry, and the first actual point lies in the central source,
where $D_0(iz)=D_0(z)$. All subsequent physical images agree.
This is not a symmetry assertion for arbitrary continued lateral
itineraries or for the two-point distortion kernel.

The final complete LAND refinement uses the ordinary physical horizon
$K_{\rm phys}=500$: each accelerated block must fit the remaining budget,
and the clock is advanced by its duration $t_i$, not by one.
The value $6$ is the subdivision depth of the \emph{physical-source
recognition routine}; it is not the depth of LAND cell refinement.
The anchor search radius is 16 grid steps. The final source meshes are
$2^{-14}$ and $2^{-15}$ after the earlier $2^{-13}$ passes.
The auxiliary basin tests have separate limits: at most six $F_V$
iterations for source-closure exclusion and ten for source entry.
These counters neither replace nor increase the physical horizon.
If any budget or chart guard is exhausted, an unfinished LAND cell
remains in the upper cover and an unfinished LOSS cell is not accepted.
These are finite sufficient tests, with no claim of convergence of the
classification algorithm.

% END sections/09_algorithms.tex

% BEGIN sections/10_certificate.tex
\section{Kernel evaluation and the final certificate}\label{sec:certificate}
The preceding algorithms supply a complete upper landing cover and a
finite lower loss union. This section evaluates the distortion kernel on
those sets. Their canonical interpretation is supplied by
Theorem~\ref{thm:canonical-identification}.
\begin{proposition}[Certified inputs for the quartic criterion]\label{prop:certificate-intro}
For $(n,m_0,t)=(3,2,5)$ the constructions give measurable sets
$E\subset A_{3,2}^+$ and $C_5^+\subset\Omega_*\cap H^+$ such that
$\area(E)>0$, $\area(L_5^+\setminus C_5^+)=0$, and
\begin{align}
 a_3&\le A_3^+=\frac{23788665}{536870912},\label{eq:cert-A}\\
 \mathop{\rm ess\,inf}_{z\in C_5^+}g_E(z)&\ge\gamma^-
 =\frac{6481603931109663}{144115188075855872}.\label{eq:cert-gamma}
\end{align}
Conjugation gives the lower-half-plane sets and bounds. In particular,
\begin{equation}\label{eq:theta}
 \theta=\frac{A_3^+}{\gamma^-}
 =\frac{2128573712302080}{2160534643703221}<1,\qquad
 \gamma^- -A_3^+=\frac{95882794203423}{144115188075855872}>0.
\end{equation}
\end{proposition}

\subsection{Which cells enter a kernel sum?}
The landing computation supplies a dyadic support $C_t^+$, while the loss
computation supplies a disjoint dyadic union $E$. Kernel evaluation may
use coarser rectangles. A row $Z_i$ covers all source positions at which
one common lower bound is computed. A loss group is the actual measurable
set $E_j\subset E$, with exact area $e_j$; its rectangular hull $Y_j$
contains it. The sets $E_j$ partition $E$, although their hulls $Y_j$ may
overlap. Only $e_j$, not $\area(Y_j)$, is used as the weight.
Using the fine original cells for every pair would also be valid. Grouping
reduces the number of pairs at the cost of a weaker whole-box bound.
The support-capacity test $b_i=\area(C_t^+\cap Z_i)>0$ keeps every row
meeting the actual support in positive area, including partial rows.

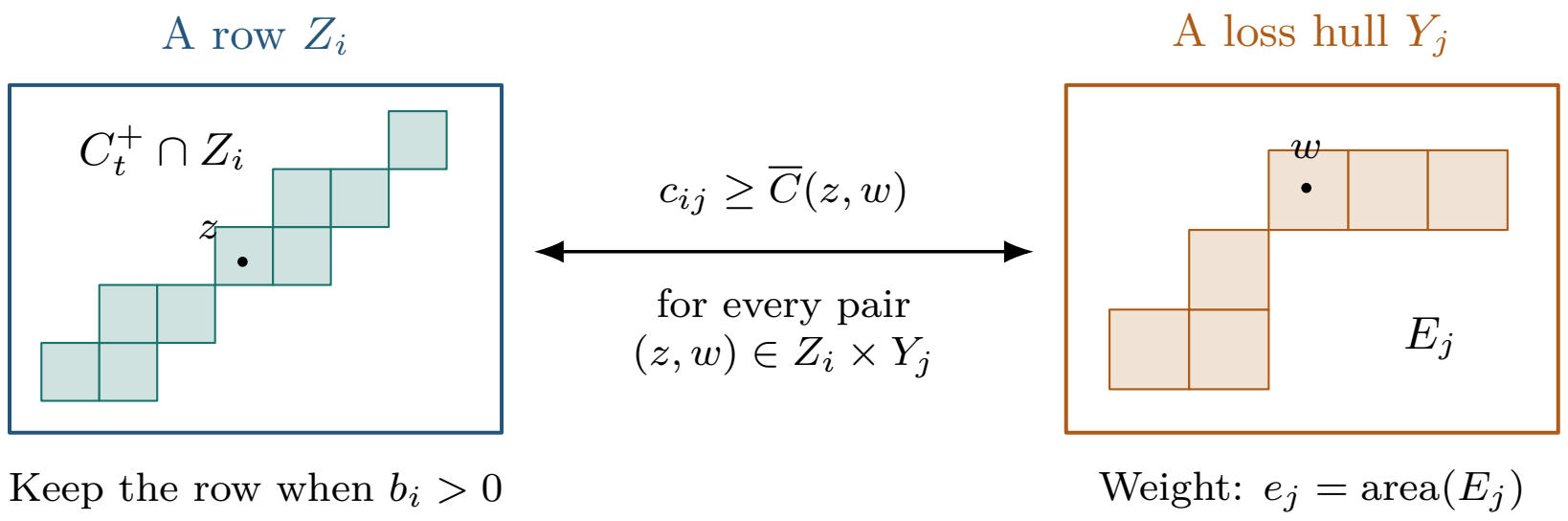
\begin{figure}[htbp]
\centering
\begin{tikzpicture}[manuscript diagram]
  \draw[diagramblue,line width=.75pt] (0,0) rectangle (3.4,2.4);
  \node[font=\footnotesize,text=diagramblue] at (1.7,2.74) {A row $Z_i$};
  \foreach \x/\y in {.22/.22,.62/.22,.62/.62,1.02/.62,1.42/1.02,1.82/1.02,1.82/1.42,2.22/1.42,2.62/1.82}
    {\filldraw[fill=diagramteal!21,draw=diagramteal,line width=.4pt]
      (\x,\y) rectangle ++(.4,.4);}
  \node[font=\footnotesize] at (1.06,1.95) {$C_t^+\cap Z_i$};
  \fill (1.61,1.18) circle[radius=1pt];
  \node[font=\scriptsize,above left=1pt] at (1.61,1.18) {$z$};
  \draw[diagramorange,line width=.75pt] (7.3,0) rectangle (10.7,2.4);
  \node[font=\footnotesize,text=diagramorange] at (9,2.74) {A loss hull $Y_j$};
  \foreach \x/\y in {7.6/.3,8.15/.3,8.15/.85,8.7/1.40,9.25/1.40,9.8/1.40}
    {\filldraw[fill=diagramorange!18,draw=diagramorange,line width=.4pt]
      (\x,\y) rectangle ++(.55,.55);}
  \node[font=\footnotesize] at (9.82,.65) {$E_j$};
  \fill (8.96,1.69) circle[radius=1pt];
  \node[font=\scriptsize,above=2pt] at (8.96,1.69) {$w$};
  \draw[<->,line width=.65pt,shorten <=2pt,shorten >=2pt] (3.55,1.25) -- (7.15,1.25)
    node[midway,map label,above=5pt] {$c_{ij}\ge\overline C(z,w)$}
    node[midway,map label,below=5pt] {for every pair\\$(z,w)\in Z_i\times Y_j$};
  \node[font=\scriptsize] at (1.7,-.40) {Keep the row when $b_i>0$};
  \node[font=\scriptsize] at (9,-.40) {Weight: $e_j=\area(E_j)$};
\end{tikzpicture}
\caption{One kernel contribution $e_j/c_{ij}^2$. The teal rectangles are retained support cells and the orange rectangles are accepted loss cells. A row $Z_i$ and a hull $Y_j$ are larger rectangles used for interval bounds. The weight is $e_j=\area(E_j)$, not $\area(Y_j)$. Locations are schematic.}\label{fig:kernel-groups}
\end{figure}
\subsection{An explicit slit-domain majorant}\label{sec:kernel}
Use exact rational endpoints $a,b$ with
$a_*\le a<b\le b_*$ and a certified complete outer cover of
$V$ compactly contained in
\[
 \Omega(a,b)=\C\setminus(({-\infty},a]\cup[b,{+\infty}))
                                      \subset\Omega_*.
\]
The majorant below is used on $\Omega(a,b)\times\Omega(a,b)$. If a majorant on all of $\Omega_*^2$ is desired in Section~\ref{sec:criterion}, extend it by $C_*$ outside this product. All computed source rows and loss hulls lie in the smaller slit. Thus every copy map in Proposition~\ref{prop:main} restricts
univalently to this smaller domain. The production endpoints are
\begin{align*}
 a&=-\frac{1265032011577179519724227}{302231454903657293676544},\\
 b&=\frac{298811595760410736573365}{151115727451828646838272}.
\end{align*}
They are inward choices for the intrinsic postcritical slit,
not the slit defining $E_\alpha$.

Define, with the square root having positive real part,
\begin{align}
 \chi(z)&=\sqrt{\frac{z-a}{b-z}},&
 T_0(z)&=|\chi(z)|\,|b-z|^2\Re \chi(z),\label{eq:kernel-chart}\\
 N(z,w)&=|\chi(z)-\chi(w)|^2,&
 \rho(z,w)^2&=\frac{N(z,w)}{N(z,w)+4\Re \chi(z)\Re \chi(w)}.
 \label{eq:rho}
\end{align}
Here $\chi$ maps the slit domain conformally to the right half-plane
and $0\le\rho<1$. A valid derivative-ratio majorant is
\begin{equation}\label{eq:explicit-C}
 \overline C(z,w)=\max\left\{1,
       \left(\frac{1+\rho(z,w)}{1-\rho(z,w)}\right)^2
                       \frac{T_0(w)}{T_0(z)}\right\}.
\end{equation}
To verify the direction of the quotient, set
$\zeta_w(z)=(\chi(z)-\chi(w))/(\chi(z)+\overline{\chi(w)})$.
Koebe's disk derivative theorem applied to
$\phi\circ\zeta_w^{-1}$ gives
\[
 \frac{|\phi'(z)|}{|\phi'(w)|}
 \le\frac{1+\rho}{(1-\rho)^3}
                     \frac{|\zeta_w'(z)|}{|\zeta_w'(w)|}.
\]
Use
\[
 \frac{|\zeta_w'(z)|}{|\zeta_w'(w)|}
 =(1-\rho^2)\frac{\Re \chi(w)}{\Re \chi(z)}
                    \frac{|\chi'(z)|}{|\chi'(w)|},\qquad
 \chi'(z)=\frac{b-a}{2\chi(z)(b-z)^2}.
\]
This simplifies to the second entry of \eqref{eq:explicit-C}.
Taking its maximum with one is conservative.

\subsection{Whole-box lower sums}
Partition the accepted loss set into disjoint measurable pieces $E_j$,
each with an exact dyadic area $e_j$ and a rectangular hull $Y_j$.
Partition the required source cover into rows $Z_i$. If
$c_{ij}$ bounds \eqref{eq:explicit-C} above on the whole product
$Z_i\times Y_j$, then
\begin{equation}\label{eq:kernel-row}
                 g_E(z)\ge\ell_i:=\sum_j\frac{e_j}{c_{ij}^2}
                                      \qquad(z\in Z_i).
\end{equation}
Rectangular hulls of loss groups may overlap; their actual accepted
subsets are disjoint, and it is their exact areas that enter the sum.

Every chart encloses the entire rectangle. For example, a direct
interval square-root image may be intersected with a center image
inflated by a derivative bound on the rectangle, using
\[
                  |\chi'(z)|=\frac{b-a}{2\sqrt{|z-a|\,|b-z|^3}}.
\]
This improves a whole-box enclosure; it does not replace it by a sample.
To evaluate \eqref{eq:rho}, use its monotonicity in
$N\ge0$ and $p=4\Re \chi(z)\Re \chi(w)>0$:
\[
                      \rho^2\le\frac{N^+}{N^++p^-}.
\]
The same $N^+$ is kept in numerator and denominator. The code
bounds distortion upward, each contribution downward, and each
accumulation downward. A failed pair contributes zero without
removing its row. Positive branch, denominator, and $\rho<1$
guards are mandatory.

For the support calculation, compute exact intersection capacities
$b_i=\area(C_5^+\cap Z_i)$ on a disjoint row partition. Only rows
with $b_i=0$ can be omitted, and their union meets the support in
a null set. Every positive-capacity row remains. Exact clipping
and area conservation verify coverage after each split.
The final lower bound is
\begin{equation}\label{eq:gamma-rows}
                    \gamma^-=\min_{i:b_i>0}\ell_i.
\end{equation}
Conjugation gives the same kernel and accepted set in the lower
half-plane. Quarter-turn symmetry of the kernel is not assumed;
all four upper-half octants are evaluated separately.

\subsection{The accepted loss set and complete landing covers}
For $(n,m_0)=(3,2)$ the union is stored as 439,868 coalesced dyadic
squares in the first octant. The four upper-half symmetry copies give
1,759,472 squares and
\begin{equation}\label{eq:E-area}
 \area(E)=\frac{5543321}{4194304}.
\end{equation}
Every cell has a whole-source inclusion and a finite loss witness.
The initial level-12 Koebe pass gives $2635251/4194304$; its adaptive
level-13 extension gives $E_K$ in Table~\ref{tab:loss-stages}.
\begin{table}[htbp]\centering\small
\begin{tabular}{@{}lrrrr@{}}\toprule
Pass & Source level & Raw units & Overlap units & New units\\\midrule
Koebe union $E_K$ & 12, 13 & 15,665,201 & 0 & 15,665,201\\
Direct ball & 12 & 5,972,100 & 4,121,687 & 1,850,413\\
Direct parent search & 11 & 4,409,424 & 216,412 & 4,193,012\\
Direct hybrid strip & 12 & 466,684 & 2,026 & 464,658\\\bottomrule
\end{tabular}
\caption{LOSS construction in units $2^{-24}$ of upper-half area.
Overlap is against the entire preceding union. The new units sum to
22,173,284, giving \eqref{eq:E-area}. Coalescing may store larger squares
without changing the union.}\label{tab:loss-stages}
\end{table}
For example, the parent-search overlap contains 115,067 units in $E_K$
and 101,345 in the preceding direct addition outside $E_K$; both are
removed. The final strip tests 298,384 squares and accepts 116,671;
its net addition is $464658/2^{24}=232329/8388608$.
The first direct pass historically selected candidates from unresolved
records of the coarse Koebe pass, so it repeated work: 913,898 accepted
level-12 squares were already wholly covered by the final adaptive $E_K$.
The later two direct passes skip cells wholly covered by the current union.
This scheduling affects runtime, not the exact union
\eqref{eq:loss-union-update}. The accepted set is fixed as $m$ grows.

\begin{table}[htbp]\centering\small
\begin{tabularx}{\textwidth}{@{}cXcl@{}}\toprule
Stage & Method and retained cover & Mesh level & Upper bound for $a_3$\\\midrule
1 & Intersection of complete Koebe covers & 13 & $360967/2097152$\\
2 & Direct ball propagation on stage 1 & 13 & $1746263/16777216$\\
3 & Direct hybrid propagation on stage 2 & 13 & $1315579/16777216$\\
4 & Direct refinement of all stage-3 roots & 14 & $3121929/67108864$\\
5 & Direct refinement of selected stage-4 roots & 14, 15 & $23788665/536870912$\\\bottomrule
\end{tabularx}
\caption{Complete LAND chain for target $W_3$. A level-$L$ square has
side $2^{-L}$ in normalized coordinates and side $\alpha^2 2^{-L}$
after applying $h$. Stage 5 keeps every unselected level-14 root and
unresolved level-15 child.}\label{tab:land-stages}
\end{table}
The numbers of retained first-octant cells after stages 1--4 are
2,888,818, 1,747,238, 1,316,413 and 3,123,471. Stage 2 removes 1,141,580
cells and stage 3 removes 430,825. Stage 4 enumerates 5,263,984 admissible
children and excludes 2,140,513, with diagonal half-square weights.
Stage 5 selects 148,767 roots and keeps all 2,974,704 unselected roots.
It enumerates 595,045 admissible children, excludes 593,404 and keeps
1,641. Its removed area is $1186767/536870912$, so
\[
 \frac{3121929}{67108864}-\frac{1186767}{536870912}
 =\frac{23788665}{536870912}.
\]
The selected roots are unresolved cells ranked by the ratio
$R_{\rm hole}/\delta$ in positive-hole-radius failure records. Selection
is a work-allocation decision, not a membership certificate: every
excluded child still requires its own witness.
Thus $3121929/67108864$ already incorporates Koebe exclusions and several
direct passes. The last subtraction operates only inside that retained
cover, so none of the earlier removed area is removed twice.

The fixed $C_5^+$ cover was obtained separately by applying target-5
exclusions to a complete target-3 cover. Nesting permits all
untested cells of that base cover to remain. The cover has exact
counted upper-half area $333543081/4294967296$ and may be larger
than a subsequently improved cover for $L_3^+$. This causes no
problem: $C_5^+$ is a valid geometric support, whereas $A_3^+$ is
the separate recurrence normalizer. The kernel result is rebound
to the improved normalizer only after map identity and target
indices have been checked.

The support computation has 5,905 final positive-capacity rows and 8,200 loss groups, giving 48,421,000 whole-box pairs.
The complete directed row replay, exact support-area conservation,
and independent disjointness audit give precisely
\eqref{eq:cert-gamma}. Together with the LAND area recount,
this proves Proposition~\ref{prop:certificate-intro}.

\subsection{Exact comparison and conclusion}
All final comparisons use rational numbers. Equivalently, the
strict comparison of the two fractions in
Proposition~\ref{prop:certificate-intro} is the integer inequality
\[
 3428307930548529927290880
             <3479784613717629946822656.
\]
It gives the coefficient and positive margin in \eqref{eq:theta}.
For reference, the magnitudes are
\begin{center}
\begin{tabular}{@{}lr@{}}
\toprule
Certified quantity & Decimal display\\\midrule
$A_3^+$ & $0.044309841468930244$\\
$\gamma^-$ on $C_5^+$ & $0.04497516200511797$\\
$A_3^+/\gamma^-$ & $0.985206934082594$\\
Full-domain kernel lower bound & $0.018149579875847383$\\
Corresponding full-domain ratio upper bound & $2.441370090769744$\\
\bottomrule
\end{tabular}
\end{center}
The decimal entries are displays of the exact certificate values,
not replacements for their directed bounds. The full-domain ratio
being above one means that this particular sufficient bound failed;
it does not give a lower bound on the optimal full-domain coefficient.

\begin{proof}[Proof of Theorem~\ref{thm:main}]
Theorem~\ref{thm:canonical-identification} identifies the validated germ with
Buff's canonical quartic germ. Its canonical domains, physical
Fibonacci tower, and intrinsic slit are those in
Sections~\ref{sec:dynamics}--\ref{sec:extension}.
Proposition~\ref{prop:main} and Lemma~\ref{lem:pairing} establish
the weighted recurrence. The finite certificates give
$E\subset A_{3,2}^+$, $L_5^+\subset C_5^+$ modulo a null set,
and the bounds of Proposition~\ref{prop:certificate-intro}.
Apply Theorem~\ref{thm:supported} with $(n,m_0,t)=(3,2,5)$.
It yields $a_{j+2}\le\theta a_j$ for every $j\ge5$, with the
fixed rational $\theta<1$ in \eqref{eq:theta}.
Thus, for example,
\[
 a_j\le a_5\,\theta^{\lfloor(j-5)/2\rfloor}\qquad(j\ge5).
\]
Lemma~\ref{lem:null-detection} proves
$\area K(D)=\area J(D)=0$.
\end{proof}

% END sections/10_certificate.tex

% BEGIN sections/11_reproducibility.tex
\section{Reproducibility and scope of the finite verification}\label{sec:reproducibility}
The readers can get access to the numerical verification code via \url{https://github.com/RunzeZhangMATH/quartic-fibonacci-verification/releases/tag/v1.0.0}.
The release contains a numerical verification package combining the authenticated production tree, the independent
identification supplement, the review's additional replay programs and
records, and a relocatable driver.

\subsection{Inputs, arithmetic and failure conventions}
The finite calculation refers to the exact degree-20 center in
Appendix~\ref{app:coefficients} and the full local certificate, whose
SHA-256 digest is
\begin{center}\small\ttfamily
7eafbd53637c8f5e95b3e2bb7930a6a8b86b\\
d205db9239fa70574fea95e1b7ee
\end{center}
All LAND, LOSS, geometry and kernel dependencies must identify this same
validated germ. Section~\ref{sec:identification} supplies its canonical
interpretation. Historical receipts that say ``conditional canonical
identification'' are preserved as records of their original computational
scope; their analytic dependency is now discharged by
Theorem~\ref{thm:canonical-identification}.

The production geometric validators use exact rational interval arithmetic.
The intensive forward and kernel loops use outward binary64 enclosures
under their checked IEEE-754 host assumptions. C++ builds require
\texttt{-frounding-math -ffp-contract=off}; fast-math is not allowed.
The independent local implementation uses integer intervals with 400
fractional bits, while its geometric and kernel charts use 120 fractional
bits. The independent kernel accumulator checks that \texttt{long double}
has at least 64 significand bits and uses directed lower summation.
A different valid implementation need not reproduce the same last bound
bit; it must prove the stated inequality in the conservative direction.

Run Python without \texttt{-O} or \texttt{-OO}: inherited verification
modules use assertions as proof gates. The top-level driver rejects an
optimized interpreter. Dependencies are Python 3 with NumPy and mpmath,
a C++17 compiler with OpenMP, and the usual Linux process environment.
The proof programs do not require the plotting dependencies SciPy and
Matplotlib. The driver raises the stack allowance to 128 MiB for generated
geometry tables, rejects a smaller hard limit, records every subprocess
and exit code, and stops at the first failure. A failed check is never
replaced by a saved success receipt.

\subsection{What each kind of check establishes}
The following distinctions are part of the verification contract.
\begin{enumerate}
\item \emph{Integrity.} SHA-256 and size comparisons establish that files
have the recorded bytes. They do not prove a numerical inequality or a
trajectory. \texttt{RELEASE\_MANIFEST.json} inventories immutable release
inputs; the original production inventory authenticates 1,682 frozen files.
\item \emph{Scalar and local validation.} The Taylor validator recomputes
the residual, derivative defect and unrestricted tail estimates. The
identification driver recomputes every contour arc and basin inequality.
The scalar area contract compares exact fractions, map identities and
indices; by itself it does not establish dynamical set membership.
\item \emph{Fresh orbit replay.} The LAND drivers rebuild the physical
engine and regenerate center or whole-cell trajectories for the required
passes. The full LOSS wrapper tests every cell of the final accepted
union, subdividing unresolved candidates to level 14 and accepting only
whole-cell witnesses. Its final gate requires zero unresolved area and
exact area $5543321/4194304$.
\item \emph{Terminal checks.} The terminal replay recomputes every accepted
hole inequality from its saved interval endpoint. This is separate from
recomputing the path to that endpoint, which the orbit replay supplies.
\item \emph{Set accounting.} Integer dyadic masks verify complete retained
covers, selected-root replacement, diagonal weights, within-batch
deduplication, and differences against the entire preceding loss union.
Kernel weights are checked against the actual accepted subsets in all
8,200 groups, not their rectangular hull areas.
\item \emph{Kernel inequalities.} Whole-box chart and row computations
re-evaluate all 48,421,000 pairs and take the minimum over all 5,905
positive-capacity rows. Exact coverage and positive capacities are checked
separately. The final comparison uses rational arithmetic.
\end{enumerate}

The review programs for local Taylor validation, geometry and kernel
bounds are separate implementations. The additional LAND and full LOSS
replay drivers reuse the production physical-orbit engine; their fresh
execution and complete-set checks are not claimed to be an independent
implementation of that engine. Mathematical continuation, disjointness,
all-depth combinatorics and the canonical uniqueness application remain
proved in the manuscript, rather than inferred from any program's return
code.

\subsection{One ordered release recipe}
Unpack the verification ZIP, enter its top-level directory, and run
\begin{Verbatim}[fontsize=\small,breaklines=true,breakanywhere=true]
python3 reproduce.py all --threads 4
\end{Verbatim}
This is the complete recipe. It runs the following phases in order;
\texttt{python3 reproduce.py all --list} prints every underlying command
and its working directory without executing anything.
\begin{Verbatim}[fontsize=\small,breaklines=true,breakanywhere=true]
python3 reproduce.py integrity
python3 reproduce.py local
python3 reproduce.py geometry
python3 reproduce.py orbits --threads 4
python3 reproduce.py kernel
python3 reproduce.py accounting
\end{Verbatim}
The full run can take hours and produce several gigabytes of intermediate
streams. The output directory \texttt{release\_runs/} contains dated logs
and a machine-readable summary. Running an individual phase reports only
that phase; it does not report a complete replay.

The local phase demonstrates the finite Newton search, reruns the full
independent Taylor validation, and replays the independent identification
contours and basin bounds. The geometry phase recomputes the legacy and
upgraded inner sets, all outer-cover pruning stages, physical isolators,
the initial cover intersection and slit endpoints. Path adapters relocate
historical filesystem locators while keeping the authenticated JSON
unchanged and comparing every numerical field exactly.

The orbit phase first rebuilds the interval executables. It then runs
\path{replay_original_stage1.py}, which reproduces the initial cover,
early Koebe exclusions, dense center enumeration, adaptive Koebe passes,
disk reuse and their exact intersection. Next
\texttt{replay\_earlier.py} regenerates LAND stages 2 and 3 and the target-5
support refinement. The complete level-14 run is
\begin{Verbatim}[fontsize=\small,breaklines=true,breakanywhere=true]
production/quartic_fibonacci_numerics/build/land \
  production/quartic_fibonacci_numerics/frozen/rigorous_l4/threshold_work/land/full_L13_vw9_retained.bin \
  replay_L14 13 14 3 500 6 8192
\end{Verbatim}
Here 13 and 14 are source mesh levels, 3 is the target depth, 500 is the
physical horizon, 6 is the source-recognition subdivision depth, and 8192
is a processing batch size. The driver subsequently checks every accepted
terminal, reruns the final selected level-15 pass and the final LOSS strip,
and links the fresh classifications to the consumed cover masks using
\texttt{check\_fresh\_cover\_chain.py}. It also regenerates terminal checks
for the target-5 support. Finally it compiles and executes
\texttt{replay\_full\_loss.cpp} on all 439,868 first-octant roots of $E$;
\texttt{verify\_full\_loss.py} enforces the zero-unresolved-area gate.

The kernel phase runs both the production whole-box calculation and the
independent chart/accumulator pair. The independent lower bounds are
written in hexadecimal floating-point notation and converted exactly to
rationals by \texttt{verify\_kernel.py}; all rows must be present and their
minimum must be at least the manuscript's $\gamma^-$.
The accounting phase checks all five LAND stages, reconstructs the full
LOSS union and each kernel weight, verifies support coverage, and checks
$A_3^+<\gamma^-$ with the correct map and target indices. The local
recount entry points differ from the original ZIP-reading scripts only
in their storage access; their mathematical tests are unchanged.

The production inventory lists 79 omitted, regenerable historical streams.
They are not assumed to have been replayed merely because their hashes
are recorded. The complete recipe regenerates the streams needed for the
cover chain and terminal checks. It directly re-certifies the final LOSS
union, so its inclusion does not require replaying every abandoned
historical search that helped discover that union.

\subsection{Review records and this revision's checks}
The directory \texttt{review\_records/} preserves the supplied review's
results separately from new outputs. The review reported fresh replay of
the complete LAND chain, the target-5 support, the physical geometry and
the full LOSS union. Its latter test used 1,208,916 oracle calls, accepted
1,016,654 leaves and left none unresolved. Its independently computed
kernel minimum was
\[
 \frac{3240802860456799}{72057594037927936}
 >\frac{6481603931109663}{144115188075855872}.
\]
The larger value is corroborating evidence; we retain the original,
weaker certified constant $\gamma^-$ in the theorem.

The dated release audit states exactly which checks were rerun while
preparing this revision. A review record, a newly computed record and a
recipe for a longer future replay are kept distinct. No full-run claim
is made on the basis of file authentication or scalar agreement alone.
The source and verification archives are closed and reopened for complete
ZIP CRC checks before delivery; their file inventories and checksums
accompany the release.

\subsection{The numerical illustration figures}\label{sec:illustrations}
Figure~\ref{fig:buff-domains} retains the numerical contours from the
reviewed illustration, with the nearly coincident $r=0.99$ curve removed.
Here $r$ is the radius in an auxiliary unit-disk coordinate on
$\Delta_\alpha$. The blue curve is the fourth-root lift of the distinguished
inverse applied to $|\zeta|=0.9999$ in that coordinate; the ochre contour
is the twelfth full pullback of $|z|=16/5$ under $G$, scaled by $\alpha$.
The dashed green curve is the $\alpha$-scaled blue curve.

Figure~\ref{fig:julia} has been recomputed at greater inverse depth. For
$w\in V$, the five physical inverse branches are
\[
 w\longmapsto\alpha\zeta\,\iota(w)^{1/4}\quad(\zeta^4=1),
 \qquad w\longmapsto\iota(\alpha w)^{1/4},
\]
where the lateral root is the branch through the positive real source.
The plotting code starts at zero and applies all five branches twelve
times, giving $5^{12}=244,140,625$ points, compared with depth ten in the
previous illustration. To resolve the critical neighborhood, it also
plots $h_mD^{-12}(0)$ for $m=1,2,3$. These are deeper physical preimages:
$\cR_m^{12}(h_mz)=0$, and expansion of the twelve return blocks has
ordinary duration between $12S_{2m-1}$ and $12S_{2m}$, namely
24--36, 60--96, and 156--252. The setwise return identities justify this
way of allocating more inverse iterations near zero.

The computations use the preserved polynomial center in ordinary double
precision. Newton inversion of that polynomial is guarded by residual
$2\cdot10^{-12}$ and the local chart radius; these floating-point guards
are consistency checks, not rigorous approximation bounds for the exact
Julia set. A three-pixel plotting mark improves print legibility. 

% END sections/11_reproducibility.tex

\clearpage
\appendix

% BEGIN sections/appendix_coefficients.tex

\section{Exact center polynomial}\label{app:coefficients}
The center $\overline H(x)=\sum_{j=0}^{20}c_jx^j$ has the following
exact terminating-decimal coefficients. They are input rationals;
the table does not assert that the exact solution is polynomial.
The rigorous unknown tail and coefficient errors are enclosed by
\eqref{eq:banach-ball}.
\begin{center}\small
\begin{tabular}{@{}rl@{}}
\toprule
$j$ & $c_j$\\\midrule
0 & $1.0$\\
1 & $-1.32742247136974996312533$\\
2 & $0.1051064512514101623758647$\\
3 & $-0.000448957995929042314437366$\\
4 & $-0.0004425008264219656946746805$\\
5 & $0.00002980095773353479818868747$\\
6 & $-0.0000005129769460912169422886435$\\
7 & $-5.756645654532490979854463\times10^{-8}$\\
8 & $5.87919818037196276101591\times10^{-9}$\\
9 & $-2.699498847563312514806136\times10^{-10}$\\
10 & $2.66549135991815927663057\times10^{-12}$\\
11 & $6.034561594024143718011813\times10^{-13}$\\
12 & $-5.132560586006654479680559\times10^{-14}$\\
13 & $1.921201584639645386507266\times10^{-15}$\\
14 & $5.863087836210449451751711\times10^{-18}$\\
15 & $-5.219730801754346384407887\times10^{-18}$\\
16 & $3.10465482528773884643936\times10^{-19}$\\
17 & $-6.128548708334479949238719\times10^{-21}$\\
18 & $-3.919798379453434350050041\times10^{-22}$\\
19 & $3.832338909411114073645153\times10^{-23}$\\
20 & $-1.291945190901080194029897\times10^{-24}$\\
\bottomrule
\end{tabular}
\end{center}

\section{Majorants for the full Taylor-algebra validation}\label{app:majorants}
This appendix specifies the finite and infinite parts of the
certificate in Proposition~\ref{prop:validated}. All the displayed
norms are the weighted $\ell^1$ norm at $R=13/5$.
Set, at the center,
\[
 a=\overline\alpha^4,\qquad G=\overline H(ax),\qquad
 S=G^4,\qquad Q=\overline H(S).
\]
The residual $\overline H+\overline\alpha^{-2}Q$ and all finite
derivative columns are polynomial algebra calculations with
directed coefficients; they are not truncated before taking
the tail norms. For a direction $(\mu,u)$, where $u(0)=0$,
\begin{align*}
 \delta G&=u(ax)+4\overline\alpha^3\mu x\overline H'(ax),\\
 \delta S&=4G^3\delta G,\\
 D\Phi(\mu,u)&=u-2\overline\alpha^{-3}\mu Q
     +\overline\alpha^{-2}\bigl(u(S)+\overline H'(S)\delta S\bigr).
\end{align*}
The finite columns use $\mu=1,u=0$ and
$\mu=0,u=x^j/R^j$ for $1\le j\le20$. The matrix $D_{20}$ is their
first 21 weighted output coefficients. A decimal Gauss--Jordan
calculation proposes $B$; the interval computation of
$I-BD_{20}$, not the accuracy of the inversion routine, validates it.

Let $\Pi_{\le20}$ and $\Pi_{>20}$ denote coefficient projections.
The residual bound is
\[
 Y\ge\|B\Pi_{\le20}\Phi(\overline\alpha,\overline H)\|
                   +\|\Pi_{>20}\Phi(\overline\alpha,\overline H)\|.
\]
For each finite unit column, add its finite inverse defect and its
exact polynomial tail norm. Their maximum is the finite-column
bound $Z_{\rm fin}<0.000055087$.

For a unit tail $u=\sum_{j>20}u_jx^j$, put $q_0=\|S\|/R<1$.
Then
\[
 \|u(S)\|\le q_0^{21},\qquad
 \|u(ax)\|\le a^{21}.
\]
The identity part cancels in $I-AD\Phi$ on a tail input, so the
conservative bound
\[
 Z_{\rm tail}\le (\|B\|+1)\overline\alpha^{-2}
       \bigl(q_0^{21}+4\|\overline H'(S)\|\,\|G\|^3a^{21}\bigr)
                         <0.017348053
\]
controls every unknown higher mode. Thus
$Z=\max\{Z_{\rm fin},Z_{\rm tail}\}$ bounds the full operator.

For clarity, the second-derivative estimate can be reproduced with
the following scalar majorants, uniformly on the ball of radius
$r$. Write $\lambda_\pm=\overline\alpha\pm r$,
$a_+=\lambda_+^4$, and $s_a=a_+R$.
For $0\le s<R$, define
\[
 H_d(s)=\sum_{j=d}^{20}\frac{j!}{(j-d)!}|c_j|s^{j-d},
 \qquad
 E_d(s)=\frac{d!\,r}{R^d(1-s/R)^{d+1}}\quad(d\ge1).
\]
Use $H_0(s)+r$ for the order-zero bound. Let
$M_d=H_d(s_a)+E_d(s_a)$ for $d=1,2$ and
$M_0=H_0(s_a)+r$. Set
\begin{align*}
 a_1&=4\lambda_+^3,& a_2&=12\lambda_+^2,\\
 G_1&=1+a_1RM_1,&
 G_2&=2a_1+a_1^2R^2M_2+a_2RM_1,\\
 s_S&=M_0^4,& S_1&=4M_0^3G_1,\\
 S_2&=12M_0^2G_1^2+4M_0^3G_2.
\end{align*}
The term $2a_1$ uses
$\|xu'(a_+x)\|\le1$ for unit $u$, since $a_+<1/2$ implies
$\sup_{j\ge1}j a_+^{j-1}=1$.
The strict inequality $s_S<R$ validates the nested composition
throughout the ball. Put
\begin{align*}
 N_0&=H_0(s_S)+r,& N_d&=H_d(s_S)+E_d(s_S)\quad(d=1,2),\\
 U_1&=\frac{1}{R(1-s_S/R)^2},&
 Q_1&=1+N_1S_1,\\
 Q_2&=2U_1S_1+N_2S_1^2+N_1S_2.
\end{align*}
The product rule for $\alpha^{-2}\HH(S)$ now gives
\[
 \sup\|AD^2\Phi\|
 \le\max\{1,\|B\|\}
       \bigl(\lambda_-^{-2}Q_2+4\lambda_-^{-3}Q_1
                                      +6\lambda_-^{-4}N_0\bigr)
 <7435.70.
\]
All upper expressions are rounded upward; positive gaps in
denominators are rounded downward. The fresh replay agrees with
every proof-bearing field of the preserved local certificate,
including the finite inverse defect, the infinite-tail bound,
and both strict contraction inequalities.

% END sections/appendix_coefficients.tex

% BEGIN sections/references.tex

% END sections/references.tex

\end{document}